\documentclass[12pt]{article}

\usepackage[utf8]{inputenc}

\usepackage{ae}
\usepackage{amsthm,amsfonts,amssymb,amsmath,oldgerm}
\usepackage{caption}
\usepackage[thinlines]{easybmat}
\usepackage{enumerate}
\usepackage{enumitem}
\usepackage{fancyhdr}
\usepackage{float}
\usepackage[T1]{fontenc}     
\usepackage[top=25mm, bottom=25mm, left=25mm, right=25mm]{geometry}
\usepackage{graphicx}
\usepackage{hyperref}
\usepackage{mdframed}
\usepackage[numbers]{natbib}
\usepackage{newtxtext,newtxmath}
\usepackage{numprint}
\usepackage{pstricks}
\usepackage{subcaption}
\usepackage{xcolor}
\usepackage{xpatch}
\usepackage[nottoc]{tocbibind}

\definecolor{DarkGreen}{rgb}{0,0.3,0}
\definecolor{LightGreen}{rgb}{0.9,1,0.9}
\definecolor{DarkPurple}{rgb}{0.5,0,0.5}
\definecolor{LightPurple}{rgb}{1,0.85,1}
\definecolor{LightRed}{rgb}{1,0.9,0.9}
\definecolor{Green}{rgb}{0,0.7,0}
\definecolor{Mblue}{rgb}{0,0.6,0.6}

\newtheorem{Th}{Theorem}

\newtheorem*{remarks}{Remarks}

\newcommand{\dd}{\mathrm{d}}

\newcommand{\Linfty}{L^{\infty}}

\newcommand{\eps}{\varepsilon}

\newcommand{\pt}{\partial}

\newcommand{\avgf}{\langle f \rangle}
\newcommand{\avgdf}{F^{\eps}}

\newcommand{\yeps}{y^{\eps}}

\newcommand{\RR}{\mathbb{R}}

\newcommand{\NN}{\mathbb{N}}
\newcommand{\ZZ}{\mathbb{Z}}

\newcommand{\KK}{\mathbb{K}}

\newcommand{\hs}{\hspace{2mm}}

\newcommand{\lnorm}{\left|\left|}
\newcommand{\rnorm}{\right|\right|}
\newcommand{\pphi}{\varphi}

\begin{document}

    
    \title{\rule{160mm}{0.7mm}\\
		\textbf{Machine Learning for highly oscillatory \\ differential equations}\\
		\rule{160mm}{0.5mm}}
	
	
	\author{Maxime Bouchereau$^1$}
	
	\date{
		\small
		$^1$IRMAR, Université de Rennes\\%
		\hs \\
		\large\today	
	}

	\maketitle
	
	\renewcommand{\abstractname}{\Large Abstract}

 
	\begin{abstract}
		\hs\\
		Highly oscillatory differential equations, commonly encountered in multi-scale problems, are often too complex to solve analytically. However, several numerical methods have been developed to approximate their solutions. 
		Although these methods have shown their efficiency, the first part of the strategy often involves heavy pre-computations from averaging theory. In this paper, we leverage neural networks (machine learning) to approximate the vector fields required by the pre-computations in the first part, and combine this with micro-macro techniques to efficiently solve the oscillatory problem. We illustrate our work by numerical simulations.
        
	\end{abstract}

	\hs

    \section{Introduction}


    Highly oscillatory differential equations are frequently used to model phenomena exhibiting a multiscale behavior with periodic dependence \cite{AVERAGING_chartier2020new, AVERAGING_chartier2010higher}. However, most of these equations cannot be solved analytically. While various numerical methods exist for autonomous differential equations \cite{ANUM_Sanz-Serna, ANUM_Butcher, ANUM_Casas, ANUM_SODE1, ANUM_Reich, ANUM_DOPRI}, they are not well-suited for highly oscillatory equations due to their stiffness.

    Although the methods mentioned above are not suitable for solving highly oscillatory ODEs, the analytical properties of these differential equations \cite{ AVERAGING_chartier2010higher, AVERAGING_chartier2012formal, AVERAGING_perko1969higher} can be leveraged to develop effective numerical techniques. In particular, uniformly accurate methods are particularly powerful, as their error bounds do not depend on the stiffness parameter \cite{AVERAGING_chartier2020new, AVERAGING_chartier2022derivative}.

    These methods rely on transforming the studied dynamical system into a modified form to enable the correct application of numerical techniques. For instance, preliminary computations often involve the slow-fast decomposition, which separates the multiscale dynamics of the system into slow and fast components. However, such a decomposition provides only an approximation of the original problem and may not by solely used to construct uniformly accurate (UA) numerical schemes. Two approaches have been developped to remedy this deficiency. The first is based on micro-macro decomposition techniques which allow for UA schemed although it increases by a factor 2 the dimension of the problem. The second is a based on pullback technique enjoys some geometric properties of the system although it induces some additional computations (inversion of a linear operator) \cite{AVERAGING_chartier2020new}.
    

    Although these methods are highly accurate, they require preliminary computations to transform the original dynamical system. For many systems, these transformations involve formal calculations, which can result in significant computational costs \cite{AVERAGING_chartier2010higher, AVERAGING_chartier2012formal}.

    The main idea of this work is to combine the theory of highly oscillatory differential equations with machine learning techniques and in particular neural networks. Given a highly oscillatory differential equation, the slow-fast decomposition of the equation is first learned by extensive simulations and then approximated by inference from a neural network. Since the slow dynamics is described by an autonomous differential equation, the theory of modified equations \cite{ANUM_GNI} is used to learn the slow dynamics. Once this representation of the slow-fast decomposition has been obtained, it is used to solve the whole equation via the micro-macro decomposition for any initial value prescribed in a \textit{learnt} range of phase space.

    \subsection{Scope of the paper}

    The paper is structured into two main sections: Section \ref{Paper2_section:theory} introduces the technique and presents the associated convergence results, while Section \ref{Paper2_section:numerical_experiments} highlights numerical experiments that demonstrate the properties of the proposed schemes.

    Subsection \ref{subsection:General_strategy} describes the strategy for integrating machine learning with averaging theory to solve highly oscillatory differential equations. The central concept is to emulate the \textit{slow-fast} decomposition using neural networks. This method requires the generation of exact data via a highly accurate yet computationally intensive numerical integrator. The data is then used to train neural networks to learn the decomposition by minimizing a $Loss_{Train}$ error function. This approach employs modified equation theory for modeling slow dynamics and an auto-encoder structure for fast dynamics. Moreover, we propose a method based on the micro-macro decomposition that avoids the need for additional training, leading to uniformly accurate approximations. Lastly, an alternative method designed for highly oscillatory autonomous systems is introduced, which eliminates the dependency on an auto-encoder.
    
    Subsection \ref{subsection:error_analysis} focuses on error analysis, providing bounds between exact and approximated solution. First method gives error bound dominated by learning errors of elements of the \textit{slow-fast} decomposition and exponential error w.r.t. $\eps$, providing a good approximation of the solution for small values of $\eps$. Second method based on micro-macro decomposition gives error bounds dominated by learning errors and step size, providing uniformly accurate bounds. Then, alternative method (reserved for oscillatory autonomous systems) provides similar error bounds to first method).
    
    
    In section \ref{Paper2_section:numerical_experiments}, numerical experiments are made. For first and second method, we learn the modified averaged field by using Forward Euler and midpoint schemes, and illustrate that elements of \textit{slow-fast} and micro-macro decomposition can properly be approximated for different dynamical systems. For autonomous systems, the experiments include a comparison with the classical method based on the \textit{slow-fast} decomposition.

    \subsection{Related work}

    Connections between differential equations, machine learning, and multiscale problems have been investigated in numerous studies. The relationship between differential equations and machine learning primarily involves learning hidden dynamics from collected data using regression techniques \cite{LEARNING_DYNAMICS_brunton2016discovering,LEARNING_DYNAMICS_du2022discovery} or statistical methods \cite{STAT_StatFieldLearning,STAT_raissi2017machineGaussian}. In contrast, the connection with multiscale problems largely focuses on learning solutions \cite{PINN_jin2023asymptotic}.

    \textbf{Modified equation coupling.} The connection to modified equation theory has been utilized to study the learning of hidden dynamics in differential equations \cite{LEARNING_DYNAMICS_du2022discovery,LEARNING_DYNAMICS_zhu2020inverse}. Offen et al. extended this approach to learn Hamiltonian functions \cite{GEO_PROP_offen2022symplectic}.

    \textbf{Neural network structure and geometric properties.} To preserve specific properties of equations, various neural network architectures have been developed. For maintaining asymptotic properties in multiscale equations, Jin et al. employed Physics-Informed Neural Networks (PINNs) \cite{PINN_jin2023asymptotic}. Other architectures have also been proposed for multiscale equations \cite{PINN_leung2022nh,PINN_lu2022solving,PINN_weng2022multiscale}, such as Convolutional Neural Networks (CNNs) \cite{PINN_wu2024capturing}. To ensure the preservation of geometric properties like Hamiltonian vector fields \cite{GEO_PROP_courtes2025neural}, Hamiltonian Neural Networks have been used \cite{GEO_PROP_david2021symplectic,GEO_PROP_greydanus2019hamiltonian}. Additionally, reciprocal mappings can be learned effectively using auto-encoders. For example, Jin et al. applied auto-encoders to learn Poisson systems \cite{GEO_PROP_jin2020learningPoisson}. Moreover, Zhu et al. have developped a method to learn Volume-Preserving vector fields \cite{GEO_PROP_VPNets}.

    \textbf{Approximation by Neural Networks.} Accurately approximating functions using neural networks has led to error estimates. Anastassiou \cite{NN_APPROX_anastassiou2000quantitative,NN_APPROX_anastassiou2023general} established convergence rates for approximating functions that map to finite- and infinite-dimensional vector spaces. These rates depend on the number of parameters and the input dimensionality. By treating neural networks as functional spaces, error bounds have been derived \cite{NN_APPROX_de2021approximation,NN_APPROX_gribonval2022approximation}. In a related study, Bach \cite{NN_APPROX_bach2021learning} provided error estimates by considering neural networks as elements within a Hilbert space. An important challenge highlighted by Mallat \cite{NN_APPROX_mallat2019sciences} is the curse of dimensionality: as the dimensionality of the vector field increases, the rate of convergence slows, requiring more parameters and data to achieve satisfactory learning outcomes in high-dimensional spaces compared to low-dimensional ones.

    \section{Approximate solutions of highly oscillatory differential equations with machine learning}\label{Paper2_section:theory}

    Consider an highly oscillatory differential equation of the form

    \begin{equation}
        \left\{ \begin{array}{c c l l}
            \dot{\yeps}(t) & = & f\left(\frac{t}{\eps} , \yeps(t) \right) \in \RR^d, & t\in [0,T] \\
            \yeps(0) & = & y_0 &
        \end{array} \right. \label{equation:highly_oscillatory}
    \end{equation}
    
    \noindent where $f: \mathbb{R} \times \mathbb{R}^d \longrightarrow \mathbb{R}^d$ is assumed to be sufficiently smooth and $2\pi$-periodic with respect to its first variable. Additionally, the parameter $\varepsilon \in (0,1]$ introduces high oscillations. By the Cauchy-Lipschitz theorem, the existence and uniqueness of a solution are guaranteed for any initial $y_0 \in \mathbb{R}^d$. Our goal is to approximate the solution over the interval $[0,T]$ at discrete times points $t_n = nh$, where $0 \leqslant n \leqslant N$, with $h = \frac{T}{N}$ representing the step size and $N$ the number of discretization points, with a uniform accuracy with respect to $\eps \in ]0,1]$.

    \subsection{General strategy}\label{subsection:General_strategy}

    As outlined in the Introduction, we will use neural networks to approximate the \textit{slow-fast} decomposition for solving \eqref{equation:highly_oscillatory}, particularly in cases where the parameter $\varepsilon$ is small. Then, we use methods based on a micro-macro decomposition to provide approximations that are valid for all $\eps \in ]0,1]$ and that enjoy uniform accuracy (UA) with respect to $\eps$.
    \subsubsection{Slow-fast decomposition}

    Consider $t \mapsto \pphi_t^f(y_0)$, the exact flow corresponding to the equation \eqref{equation:highly_oscillatory}. It is known \cite{AVERAGING_chartier2015higher} that one may construct an approximation of $t \longmapsto \pphi_t^f(y_0)$ using the following decomposition:

    \begin{equation}
        \pphi_t^f(y_0) \approx \phi^\eps_{\frac{t}{\eps}}\left( \psi_t^\eps(y_0) \right) \label{formula:slow_fast_decomposition_formal}
    \end{equation}
    
    \noindent where $\tau \mapsto \phi_\tau^\eps$ is $2\pi$-periodic, $\phi_0^\eps = \text{Id}$\footnote{The assumption $\phi_0^\eps = \text{Id}$ arises from stroboscopic averaging \cite{AVERAGING_chartier2020new}. An alternative convention, standard averaging, considers the average $\langle \phi_\cdot^\eps \rangle = \text{Id}$.} and $\phi^\eps = \text{Id} + \mathcal{O}(\eps)$. $\phi^\eps: (\tau, y) \longmapsto \phi_\tau^\eps(y)$ introduces high oscillations while $t \mapsto \psi_t^\eps$ results in global drift. $t \mapsto \psi_t^\eps(y_0)$ is the solution of an autonomous differential equation with an associated vector field, denoted $F^\eps$, called the \textit{averaged field}. $\phi^\eps$ and $F^\eps$ are usually expressed as formal series expansions in terms of $\eps$:

    \begin{equation}
        \phi^\eps_\tau(y) = y + \sum_{j=1}^{+\infty}\eps^j\tilde{\phi}_j(\tau,y) \quad \text{and} \quad \avgdf(y) = \avgf(y) + \sum_{j=1}^{+\infty}\eps^jF_j(y), \label{variable_change_averaged_field_formal}
    \end{equation}
    
    \noindent where $\avgf$ is the \textit{average field} defined as

    \begin{equation}
        \avgf(y) := \frac{1}{2\pi}\int_{0}^{2\pi}f(\tau,y)\,d\tau,
    \end{equation}
    
    \noindent and the coefficient functions $\tilde{\phi_j}$ and $F_j$ depend on derivatives of $f$ and $\avgf$. However, the formula referenced as equation \eqref{formula:slow_fast_decomposition_formal} may not hold because $\phi^\eps$ and $\avgdf$ are defined through formal series that typically do not converge. Despite the non-convergence of the formal series \eqref{variable_change_averaged_field_formal}, it is possible to construct an approximation of the slow-fast decomposition. Specifically, for all $n \in \mathbb{N}$, let $F^{\eps,[n]}$ and $\phi^{\eps,[n]}$ represent the truncation of $\phi^{\eps}$ and $F^{\eps}$ obtained by neglecting terms of order $\mathcal{O}(\eps^{n+1})$ terms. According to averaging theory (see Chartier et al. \cite{AVERAGING_chartier2012formal,AVERAGING_chartier2015higher}), there exists $\eps_0 > 0$ such that, for all $\eps \in ]0,\eps_0]$, there is an $n_\eps \in \mathbb{N}$ such that, for all $t \in [0,T]$,

    \begin{equation}
        \left| \pphi_t^f(y_0) - \phi^{\eps,[n_\eps]}_{\frac{t}{\eps}}\left(\pphi_t^{F^{\eps,[n_\eps]}}(y_0)\right) \right| \leqslant M e^{-\frac{\beta\eps_0}{\eps}} \label{estimate:exponential_error}
    \end{equation}
    
    \noindent for some constants $M,\beta > 0$ independent of $\eps$. Since the error bound decays very rapidly as $\eps \rightarrow 0$, this decomposition provides an accurate approximation of the solution for small values of $\eps$. However, this approximation does not hold for values of $\eps$ close to $1$. Consequently, we need uniformly accurate methods for all values of $\eps \in ]0,1]$.

    \subsubsection{Micro-Macro decomposition}

    Although the \textit{slow-fast} decomposition is efficient for small values of the parameter $\varepsilon$, it does not provide a uniformly accurate approximation of the solution with respect to $\varepsilon$ due to the exponential error described in the referenced estimate \eqref{estimate:exponential_error}. To address this issue, the micro-macro decomposition, as outlined by Chartier et al. \cite{AVERAGING_chartier2020new} leverages the multiscale structure of the equation. This approach decomposes the solution into a \textit{slow-fast} component and a remainder. By considering $\phi^{[p]}$ and $F^{[p]}$ as truncations of $\phi^\varepsilon$ and $F^\varepsilon$, respectively, where terms of order $\mathcal{O}(\varepsilon^{p+1})$ are neglected, the following result holds for all $t \in [0,T]$:

    \begin{equation}
        \yeps(t) = \phi^{[p]}_{\frac{t}{\epsilon}}(v(t)) + w(t),
    \end{equation}
    
    \noindent where $(v,w)$ is the solution of the micro-macro system, where the second equation is obtained by using \eqref{equation:highly_oscillatory} and chain rule:
    
    \begin{equation}
        \left\{
        \begin{array}{lcl}
            \dot{v}(t) & = & F^{[p]}(v(t)) \\
            \dot{w}(t) & = & f\left(\frac{t}{\varepsilon} , \phi^{[p]}(v(t))+w(t)\right) - \frac{1}{\varepsilon}\left(\frac{\partial}{\partial \tau}\phi^{[p]}_{\frac{t}{\varepsilon}}\right)(v(t)) - \left(\frac{\partial}{\partial y}\phi^{[p]}_{\frac{t}{\varepsilon}}\right)(v(t))F^{[p]}(v(t)),
        \end{array}
        \right.
    \end{equation}
    
    \noindent with $(v,w)(0) = (\yeps(0),0)$. Using a numerical integrator of order $p$ for this system produces an approximate solution of order $p$, denoted by $(v_n,w_n)_{n \in \mathbb{N}}$, with respect to the step size $h$:
    
    \begin{equation}
        \max_{0 \leqslant n \leqslant N}\left| \phi^{[p]}(v(t_n)) + w(t_n) - \phi^{[p]}(v_n) - w_n \right| \leqslant \overline{M}h^p,
    \end{equation}
    
    \noindent where the constant $\overline{M}$ is independent of $\varepsilon$ and of $h$. This ensures a uniformly accurate approximation, making the micro-macro method a \textit{Uniformly Accurate} (UA) method.

    \subsubsection{Machine learning method}

    The primary goal of this paper is to approximate the mapping $(\tau,y,\varepsilon) \mapsto \phi^{\varepsilon,[n_\varepsilon]}_\tau(y)$ and the flow $(t,y, \varepsilon) \mapsto \pphi_t^{F^{\varepsilon,[n_\varepsilon]}}(y)$ using neural networks, with the aim of achieving a structure for the approximated solution analogous to \eqref{formula:slow_fast_decomposition_formal}.\\

    Since the flow $(t,y,\varepsilon) \mapsto \pphi_t^{F^{\varepsilon,[n_\varepsilon]}}(y)$ is associated with following autonomous differential equation

    \begin{equation}
        \left\{
        \begin{array}{lcl}
            \dot{z}(t) & = & F^{\eps,[n_{\eps}]}(z(t)) \\
            z(0) & = & y_0,
        \end{array}
        \right. \label{equation:slow_dynamics}
    \end{equation}

    \noindent backward error analysis (also known as modified equation theory \cite{ANUM_GNI}) can be applied to approximate solutions of equation \eqref{equation:slow_dynamics}. Using a numerical method $\Phi_h$ of order $p$ with a given step size $h$, let $\widetilde{F_h^{\varepsilon,[n_\varepsilon]}}$ denote the modified field $F^{\varepsilon,[n_\varepsilon]}$ corresponding to $\Phi_h$ and $\widetilde{F_h^{\varepsilon,[n_\varepsilon],[q]}}$ represent the truncation of $\widetilde{F_h^{\varepsilon,[n_\varepsilon]}}$ at order $h^{p+q-1}$ \cite{FORMER_PAPERS_bouchereau2023machine,ANUM_GNI}. Then, for any compact set $\KK \subset \RR^d$ and $y \in \KK$, we have:
    
    \begin{equation}
        \left| \pphi_h^{F^{\varepsilon,[n_\varepsilon]}}(y) -\Phi_h^{\widetilde{F_h^{\varepsilon,[n_\varepsilon],[q]}}}(y) \right| \leqslant \overline{C}h^{p+q-1}.
    \end{equation}
    
    Next, we approximate $\widetilde{F_h^{\varepsilon,[n_\varepsilon],[q]}}$ using a neural network, denoted by $F_\theta$, which can be interpreted as a perturbation of the averaged field:
    
    \begin{equation}
        F_\theta(y,h,\varepsilon) = \avgf + R_{\theta,F}(y,h,\varepsilon)
    \end{equation}
    
    \noindent where the perturbation $R_{\theta,F}$ is modeled as a multilayer perceptron (MLP). Structure of $F_\theta$ follows both averaged field $\eqref{variable_change_averaged_field_formal}$ and modified fied structures \cite{FORMER_PAPERS_bouchereau2023machine, ANUM_GNI}.\\
    
    Furthermore, we must account for the high-oscillation generator described by the map $(\tau,y,\varepsilon)\mapsto \phi^\varepsilon_\tau(y)$. To get the autonomous equation from the solution of the initial equation, it is necessary to approximate both $\phi^{\varepsilon,[n_\varepsilon]}_\tau$ and its inverse\footnote{this map is invertible if $\varepsilon$ is small enough, as a perturbation of the identity.}. For this purpose, an auto-encoder consisting of a pair of neural networks, $\left(\phi_{\theta,+}(\tau,y,\varepsilon),\phi_{\theta,-}(\tau, y,\varepsilon)\right)$ is used to approximate $\left(\phi^\varepsilon_\tau(y),(\phi^\varepsilon_\tau)^{-1}(y)\right)$. To preserve the structure of $\phi^\varepsilon$ and its inverse, as required by stroboscopic averaging, both neural networks are designed as perturbations of the identity:

    \begin{equation}
        \phi_{\theta,+}(\tau,y,\eps) = y + \eps\left[ R_{\theta,+}\left(\cos(\tau),\sin(\tau),y,\eps\right) - R_{\theta,+}\left(1,0,y,\eps\right) \right]
    \end{equation}
    
    \noindent and 

    \begin{equation}
        \phi_{\theta,-}(\tau,y,\eps) = y + \eps\left[ R_{\theta,-}\left(\cos(\tau),\sin(\tau),y,\eps\right) - R_{\theta,-}\left(1,0,y,\eps\right) \right]
    \end{equation}
    
    \noindent where both $R_{\theta,+}$ and $R_{\theta,-}$ are modeled as MLPs. Additionally, the mappings $y \longmapsto \phi_{\theta,+}(\tau,\cdot,\varepsilon) \circ \phi_{\theta,-}(\tau,\cdot,\varepsilon)(y)$ and $y \longmapsto \phi_{\theta,-}(\tau,\cdot, \varepsilon) \circ \phi_{\theta,+}(\tau,\cdot,\varepsilon)(y)$ must closely approximate the identity to satisfy the auto-encoder structure \cite{GEO_PROP_jin2020learningPoisson}. Moreover, dependancy w.r.t. $\tau$ with trigonometric functions comes from $2\pi$-periodicity of $\phi^\eps$.\\

    The complete numerical procedure consists of three main steps. First, data are collected by accurately simulating the exact flow at various points in the domain, requiring a large number of simulations and high precision to ensure a reliable approximation of the averaged field and high oscillation generator. Second, the neural networks are trained individually by minimizing a prescribed loss function to optimize their performance. Finally, given the initial data, an approximation of the exact solution is obtained by applying the same numerical scheme to the neural networks as was used during training.

    \begin{enumerate}[label=\textbf{\arabic*.}]
        \item \textbf{Construction of the data set:} $K$ ''initial'' data $y_0^{(k)}$ at time $t_0^{(k)}$ are randomly selected into a compact set $\Omega \subset \mathbb{R}^d$ (where we want to simulate the solution) with uniform distribution. The initial time $t_0^{(k)}$ is randomly chosen in $[0,2\pi]$. Then, for all $0 \leqslant k \leqslant K-1$, we compute a very accurate approximation of the exact flow at times $h^{(k)}$ with initial condition $y_0^{(k)}$, denoted $y_1^{(k)}$. step sizes $h^{(k)}$ and high oscillation parameters $\varepsilon^{(k)}$ are chosen in domains $[h_-,h_+]$ and $[\varepsilon_-, \varepsilon_+]$ respectively (we actually pick the values $\log h^{(k)}$ and $\log \varepsilon^{(k)}$ randomly in the domains $[\log h_-,\log h_+]$ and $[\log \varepsilon_-,\log \varepsilon_+]$ with uniform distribution).

        \item \textbf{Training the neural networks:} We minimize the Mean Squared Error (MSE), denoted $Loss_{Train}$, which measures the difference between the predicted data $\hat{y}_1^{(k)}$ and the ``exact data'' $y_1^{(k)}$, by computing the optimal parameters of the NN over $K_0$ data (where $1 \leqslant K_0 \leqslant K-1$) using a gradient method:

        \footnotesize\begin{eqnarray}
            Loss_{Train} & = & \frac{1}{K_0}\sum_{k=0}^{K_0-1}\Big| \underbrace{  \phi_{\theta,+}\Big(\frac{t_0^{(k)}+h^{(k)}}{\eps^{(k)}} , \cdot , \eps^{(k)}\Big) \circ \Phi^{F_{\theta}(\cdot,h^{(k)},\eps^{(k)})}_{h^{(k)}} \circ \phi_{\theta,-}\Big(\frac{t_0^{(k)}}{\eps^{(k)}} , \cdot , \eps^{(k)}\Big)\big( y_0^{(k)} \big)}_{=\hat{y_1}^{(k)}} - \underbrace{\varphi^f_ {h^{(k)}}\big( y_0^{(k)} \big)}_{=y_1^{(k)}} \Big|^2 \nonumber \\
            & + & \frac{1}{K_0}\sum_{k=0}^{K_0-1}\Big| \phi_{\theta,+}\Big(\frac{t_0^{(k)}}{\eps^{(k)}} , \cdot , \eps^{(k)}\Big) \circ \phi_{\theta,-}\Big(\frac{t_0^{(k)}}{\eps^{(k)}} , \cdot , \eps^{(k)}\Big)\big( y_0^{(k)} \big) - y_0^{(k)} \Big|^2 \label{Loss_train} \\
            & + & \frac{1}{K_0}\sum_{k=0}^{K_0-1}\Big| \phi_{\theta,-}\Big(\frac{t_0^{(k)}}{\eps^{(k)}} , \cdot , \eps^{(k)}\Big) \circ \phi_{\theta,+}\Big(\frac{t_0^{(k)}}{\eps^{(k)}} , \cdot , \eps^{(k)}\Big)\big( y_0^{(k)} \big) - y_0^{(k)} \Big|^2 \nonumber
        \end{eqnarray}

        \normalsize
        The first term of $Loss_{Train}$ corresponds to the structure of the slow-fast decomposition equation, while the second and third terms correspond to the auto-encoder structure for the pair $\left(\phi_{\theta,+},\phi_{\theta,-}\right)$.\\

        Simultaneously, we compute the value of another MSE, denoted as $Loss_{Test}$, which measures the difference between the predicted data $\hat{y}_1^{(k)}$ and the ``exact data'' $y_1^{(k)}$ for a subset of initial values that were not used during the training of the neural networks. The purpose of this step is to evaluate the performance of the training process on ``unknown'' initial values:

	  \footnotesize\begin{eqnarray}
            Loss_{Test} & = & \frac{1}{K-K_0}\sum_{k=K_0}^{K-1}\Big| \underbrace{  \phi_{\theta,+}\Big(\frac{t_0^{(k)}+h^{(k)}}{\eps^{(k)}} , \cdot , \eps^{(k)}\Big) \circ \Phi^{F_{\theta}(\cdot,h^{(k)},\eps^{(k)})}_{h^{(k)}} \circ \phi_{\theta,-}\Big(\frac{t_0^{(k)}}{\eps^{(k)}} , \cdot , \eps^{(k)}\Big)\big( y_0^{(k)} \big)}_{=\hat{y_1}^{(k)}} \nonumber \\
            & - & \underbrace{\varphi^f_ {h^{(k)}}\big( y_0^{(k)} \big)}_{=y_1^{(k)}} \Big|^2 \nonumber \\
            & + & \frac{1}{K-K_0}\sum_{k=K_0}^{K-1}\Big| \phi_{\theta,+}\Big(\frac{t_0^{(k)}}{\eps^{(k)}} , \cdot , \eps^{(k)}\Big) \circ \phi_{\theta,-}\Big(\frac{t_0^{(k)}}{\eps^{(k)}} , \cdot , \eps^{(k)}\Big)\big( y_0^{(k)} \big) - y_0^{(k)} \Big|^2 \\
            & + & \frac{1}{K-K_0}\sum_{k=K_0}^{K-1}\Big| \phi_{\theta,-}\Big(\frac{t_0^{(k)}}{\eps^{(k)}} , \cdot , \eps^{(k)}\Big) \circ \phi_{\theta,+}\Big(\frac{t_0^{(k)}}{\eps^{(k)}} , \cdot , \eps^{(k)}\Big)\big( y_0^{(k)} \big) - y_0^{(k)} \Big|^2 \nonumber
        \end{eqnarray}

        \normalsize If $Loss_{Train}$ and $Loss_{Test}$ exhibit similar decay patterns, it indicates the absence of overfitting. In this case, the neural network model does not merely fit the training data but also maintains its ability to generalize and perform accurately on unknown data, which is its primary objective.

        \item \textbf{Numerical approximation using \textit{slow-fast} decomposition:} At the end of training process, an accurate approximation $F_{\theta}(\cdot,h,\eps)$ of $\widetilde{F_h^{\eps,[n_\eps],[q]}}$ is obtained. This approximation is then used to compute the successive values of $\big(\Phi^{\widetilde{F_h^{\eps,[n_\eps]}}}_{h}\big)^n(y_0)$ for $n=0,\ldots,N$. Additionally, an accurate approximation $\phi_{\theta,+}(\cdot,\cdot,\eps)$ of $\phi^\eps$. As a result, an approximation of the solution can be visualized by plotting:

        \begin{equation}
            \yeps_{\theta,n} := \phi_{\theta,+}\left(\frac{t_n}{\eps} , \left(\Phi_{h}^{F_{\theta}(\cdot,h,\eps)}\right)^n(y_0),\eps\right)
        \end{equation}
        
        for all $n=0,\cdots,N$.

        \item \textbf{Numerical approximation using micro-macro correction:} At the conclusion of the training, an alternative option is to reproduce the numerical integration using the Micro-Macro method by plotting:

        \begin{equation}
            \yeps_{\theta,n} = \phi_{\theta,+}\left(\frac{t_n}{\eps},v_{\theta,n},\eps\right) + w_{\theta,n}, \label{formula:Micro-Macro_ML}
        \end{equation}
        
        \noindent where, for all $n\in\NN$:

        \begin{equation}
            v_{\theta,n} = \left(\Phi_h^{F_\theta(\cdot,h,\eps)}\right)^n(y_0),
        \end{equation}
        
        and
        
        \begin{equation}
            w_{\theta,n+1} = \left(\Phi_{t_n,h}^{g_\theta(\frac{t_n}{\eps},\cdot,v_{\theta,n},\eps)}\right)(w_{\theta,n}),
        \end{equation}
        
        where, for all $(\tau,w,v,\eps) \in [0,2\pi]\times\Omega\times\Omega\times]0,1]$,

        \begin{equation}
            g_\theta\left(\tau,w,v,\eps\right) = f\left(\tau,\phi_{\theta,+}\left(\tau,v,\eps\right)+w\right) - \frac{1}{\eps}\pt_\tau\phi_{\theta,+}\left( \tau , v , \eps \right) - \pt_y\phi_{\theta,+}\left(\tau,v,\eps\right)F_\theta\left(v,0,\eps\right).
        \end{equation}

        The micro-macro correction has the advantage of eliminating the need to learn the entire micro-macro vector field, which would require a dataset with double the dimensions and more data. Furthermore, this method does not necessitate more training than the \textit{slow-fast} decomposition-based approach.
        
    \end{enumerate}

    \subsubsection{An alternative method for autonomous systems}

    Let us consider now autonomous highly oscillatory systems of the form

    \begin{equation}
        \dot{\yeps}(t) = \frac{1}{\eps}A\yeps(t) + g(\yeps(t)), \label{equation:autonomous_highly_oscillatory}
    \end{equation}
    
    \noindent where $A \in \mathcal{M}_d(\RR)$  be a matrix with eigenvalues in $i\ZZ$ and let $g: \RR^d \longrightarrow \RR^d$ be a smooth function. A natural approach to analyze this system involves introducing the change of variables $z^\eps(t) = e^{-\frac{t}{\eps}A}\yeps(t)$ to get the new system
    
    \begin{equation}
        \dot{z^\eps}(t) = e^{-\frac{t}{\eps}A}g\left(e^{\frac{t}{\eps}A}z^\eps(t)\right),
    \end{equation}
    
    \noindent which is a system of the form \eqref{equation:highly_oscillatory}.\\
    
    This system can also be studied directly in its autonomous form \eqref{equation:autonomous_highly_oscillatory}. A normal form theorem \cite{AVERAGING_chartier2016averaging} guarantees the existence of a matrix $A^\eps$ and a vector field $g^\eps$ such that $A^\eps$ generates a periodic flow $\tau \longmapsto \Phi_\tau^\eps$, $g^\eps$ generates a flow $t \longmapsto \pphi_t^{g^\eps}$, the Lie bracket $\left[A^\eps,g^\eps\right]$ vanishes\footnote{$A^\eps$ can be considered as a linear vector field $y \longmapsto A^\eps y$. So we have $\left[A^\eps,g^\eps\right](y) = A^\eps g^\eps(y) - \pt_yg^\eps(y)A^\eps(y)$}. Moreover, for any $T>0$, there exists a positive constant $C_T$ such that for all $t\in[0,T]$ and $\eps \in ]0,1]$, the following holds:
    
    \begin{equation}
        \left|\yeps(t) - \phi_{\frac{t}{\eps}}^\eps\left(\pphi_t^{g^\eps}(\yeps(0))\right)\right| \leqslant C_Te^{-\frac{C_T}{\eps}}.
    \end{equation}
    
    The main idea here is to approximate the flows $\Phi^\eps$ and $\pphi^{g^\eps}$ using neural networks. This approach has the advantage of eliminating the need for an auto-encoder to learn $\phi^\eps$.\\
    
    Specifically, we approximate $\pphi_h^{g^\eps}$ with a neural network, denoted $\pphi_\theta$ modeled as an identity perturbation
    
    \begin{equation}
        \pphi_\theta(y,h,\eps) = y + hR_{\theta,\pphi}(y,h,\eps),
    \end{equation}
    
    \noindent where $R_{\theta,\pphi}$ is a multilayer perceptron.\\
    
    We also approximate the periodic flow $\tau \longmapsto \Phi_\tau^\eps$ using a neural network denoted $\phi_\theta$ modeled as an identity perturbation:
    
    \begin{equation}
        \phi_\theta(\tau,y,\eps) = y + \left[R_{\theta,\phi}(\cos(\tau),\sin(\tau),y,\eps) - R_{\theta,\phi}(1,0,y,\eps)\right].
    \end{equation}
    
    \noindent where $R_{\theta,\phi}$ is also a multilayer perceptron. Structure of $\phi_\theta$ follows structure of $\Phi^\eps$.\\
    
    As in the classical case, the numerical procedure is divided into three main steps: data generation, training the neural networks through \textit{loss} minimization, and numerical integration. Specifically, we adhere to the same methodology as in the classical case:

    \begin{enumerate}[label=\textbf{\arabic*.}]
        \item \textbf{Data set construction:}
        We randomly select $K$ initial conditions $y_0^{(k)}$ at time $t=0$ from a compact set $\Omega \subset \RR^d$ with a uniform distribution. Additionally, we randomly choose $h^{(k)}\in[h_-, h_+]$ and $\eps^{(k)}\in[\eps_-,\eps_+]$ where $\log h^{(k)}$ and $\log \eps^{(k)}$ are uniformly distributed. For each $0 \leqslant k \leqslant K-1$, we compute $y_1^{(k)} = \yeps(h)$ using a highly accurate (and computationally expensive) integrator, providing a very precise approximation of the exact flow.

        \item \textbf{Training the neural networks:} We minimize the MSE loss function, $Loss_{Train}$ function, which quantifies the difference between the "exact" data and the predictions. This optimization is performed using a gradient-based method to find the optimal neural network parameters over $K_0$ data points:

        \begin{eqnarray}
            Loss_{Train}  & := & \frac{1}{K_0}\sum_{k=0}^{K_0-1}\left| y_1^{(k)} - \phi_\theta\left( \frac{h^{(k)}}{\eps^{(k)}} , \pphi_\theta(y_0^{(k)},h^{(k)},\eps^{(k)}) , \eps^{(k)} \right) \right|^2 \nonumber \\
            & + & \frac{1}{K_0}\sum_{k=0}^{K_0-1}\left| \phi_\theta\left(\frac{h^{(k)}}{\eps^{(k)}},\pphi_\theta(y_0^{(k)},h^{(k)},\eps^{(k)}),\eps^{(k)}\right)\right. \\
            & - & \left. \pphi_\theta\left( \phi_\theta\left( \frac{h^{(k)}}{\eps^{(k)}},y_0^{(k)},\eps^{(k)} \right),h^{(k)},\eps^{(k)} \right) \right|^2. \nonumber
        \end{eqnarray}

        The first term of $Loss_{Train}$ captures the structure of the equation involving both flows, while the second term enforces the property of flow commutativity, which is equivalent to the vanishing of the Lie bracket of the associated vector fields.

        Simultaneously, we compute the MSE $Loss_{Test}$ to evaluate the performance of the training on "unknown" data:

        \begin{eqnarray}
            Loss_{Test}  & := & \frac{1}{K-K_0}\sum_{k=K_0}^{K-1}\left| y_1^{(k)} - \phi_\theta\left( \frac{h^{(k)}}{\eps^{(k)}} , \pphi_\theta(y_0^{(k)},h^{(k)},\eps^{(k)}) , \eps^{(k)} \right) \right|^2 \nonumber \\
            & + & \frac{1}{K-K_0}\sum_{k=K_0}^{K-1}\left| \phi_\theta\left(\frac{h^{(k)}}{\eps^{(k)}},\pphi_\theta(y_0^{(k)},h^{(k)},\eps^{(k)}),\eps^{(k)}\right) \right. \\
            & - & \left. \pphi_\theta\left( \phi_\theta\left( \frac{h^{(k)}}{\eps^{(k)}},y_0^{(k)},\eps^{(k)} \right),h^{(k)},\eps^{(k)} \right) \right|^2. \nonumber
        \end{eqnarray}

        \item \textbf{Integration:} At the end of the training, we obtain an accurate approximation of $\Phi^\eps$ and $\pphi^{g^{\eps}}$. We then plot the points $(\yeps_{\theta,n})_{0 \leqslant n \leqslant N} = \left( \phi_\theta\left( \frac{t_n}{\eps} , \pphi_\theta(\cdot,h,\eps)^n(\yeps(0)),\eps \right) \right)_{0 \leqslant n \leqslant N}$ for $n = 0,\cdots,N$.
    \end{enumerate}

    \subsection{Error analysis}\label{subsection:error_analysis}

    In this subsection, we analyze the error arising from the methods described in the previous section. More specifically, we provide estimates of the global error for each of the standard numerical methods.

    \subsubsection{Slow-fast decomposition}

    The \textit{slow-fast} decomposition offers a direct numerical method along with corresponding error bounds.

    \begin{Th} \label{theorem:highly_oscillatory_general_slow-fast}

        Let us denote the following learning errors:

        \begin{enumerate}[label = \textbf{(\roman*)}]
            \item \textbf{Learning error for high oscillation generator:}

            \begin{equation}
                \delta_{\phi,+} := \lnorm \frac{\phi^{\eps,[n_\eps]} - \phi_{\theta,+}}{\eps} \rnorm_{\Linfty\left([0,2\pi] \times \Omega \times [0,\eps_+] \right) }
            \end{equation}

            \item \textbf{Learning error for modified averaged field:}

            \begin{equation}
                \delta_{F} := \lnorm \widetilde{F_h^{\eps,[n_\eps],[q]}} - F_{\theta} \rnorm_{\Linfty\left( \Omega \times [0,h_+] \times [0,\eps_+] \right)},
            \end{equation}
        \end{enumerate}
        
    Let us consider these two hypotheses:

    \begin{enumerate}[label = \textbf{(\roman*)}]
        \item For $f_1,f_2:\Omega \rightarrow \RR^d$ smooth enough, we have, for all $h>0$,

        \begin{equation}
            \lnorm \Phi_h^{f_1} - \Phi_h^{f_2} \rnorm_{\Linfty(\Omega)} \leqslant Ch\lnorm f_1 - f_2 \rnorm_{\Linfty(\Omega)}
        \end{equation}
        
        for some positive constant $C>0$ independent on $f_1$ and $f_2$.
    
        \item For $f: \Omega \rightarrow \RR^d$ sufficiently smooth, there exists $L_f > 0$ s.t. for all $y1,y2\in\Omega$ and for all $h>0$, we have

        \begin{equation}
            \left| \Phi_h^f(y_1) - \Phi_h^f(y_2) \right| \leqslant (1+L_fh)|y1-y2|
        \end{equation}
    \end{enumerate}

    Let $y_{\theta,n}^\eps$ denote the following numerical flow:

    \begin{equation}
        \yeps_{\theta,n} := \phi_{\theta,+}\left( \frac{t_n}{\eps} , \left(\Phi_h^{F_{\theta}(\cdot,h,\eps)}\right)^n(y_0) , \eps \right).
    \end{equation}
    
    Then there exist constants $\lambda , \alpha > 0$ (independent of $h,\eps)$ s.t. for all $h \leqslant h_+$ and $\eps \in ]0,\eps_0]$,
    
    \begin{equation}
        \left| y^{\eps}(t_n) - y^{\eps}_{\theta,n} \right| \leqslant Me^{-\frac{\beta\eps_0}{\eps}} + \delta_{\phi,+}\eps + (1+\alpha\eps)\frac{e^{\lambda T}-1}{\lambda}\left[\overline{C}h^{p+q-1}+C\delta_F\right].
    \end{equation}
    \end{Th}

    \begin{remarks}
        \begin{enumerate}[label=\textbf{(\roman*)}]
            \item The vector fields $\phi^{\eps,[n_{\eps}]}_\cdot$ and $\widetilde{F_h^{\eps,[n_\eps]}}$ are smooth by construction, as shown in the formulas \eqref{sequence_approximation_variable_change} and \eqref{sequence_approximation_averaged_field}). Similarly, $\phi_{\theta,-}$, $\phi_{\theta,+}$ and $F_\theta$ are also smooth because they result from the composition of affine function $A_1,\cdots,A_{L+1}$ and nonlinear activation functions  $\sigma_1,\cdots,\sigma_L$. For $L$ layers, the output of the neural network takes the form $A_{L+1} \circ \Sigma_L \circ A_L \circ \cdots \Sigma_1 \circ A_1$. Therefore, if the activation functions are smooth, so are $\phi_{\theta,-}$, $\phi_{\theta,+}$, and $F_\theta$. This is the case, for example, when the $\Sigma_i$'s are hyperbolic tangent functions.

            \item A similar error estimate holds for a variable step-size implementation of the numerical method $\Phi$: if we use the step sequence $0\leqslant h_{j} \leqslant h_+$, then $T = h_0 + \cdots + h_{N-1}$, $t_n = h_0 + \cdots + h_{n-1}$ and
            
            \begin{equation}
                y_{\theta,n} = \phi_{\theta,+}\left(\frac{t_n}{\eps}, \Phi^{F_{\theta}(\cdot,h_{n-1},\eps)}_{h_{n-1}} \circ \ldots \circ \Phi^{F_{\theta}(\cdot,h_0,\eps)}_{h_0}  ( \yeps(0) ) , \eps \right)
            \end{equation}
            
            where constants $\lambda , \alpha > 0$ (independent of $h,\eps)$ s.t. for all $h \leqslant h_+$ and $\eps \in ]0,1]$:

            \begin{equation}
                \underset{0 \leqslant n \leqslant N}{Max}\left|\yeps(t_n) - y_{\theta,n}\right| \leqslant Me^{-\frac{\beta\eps_0}{\eps}} + \delta_{\phi,+}\eps + (1+\alpha\eps)\frac{e^{\lambda T}-1}{\lambda}\left[\overline{C}h^{p+q-1} + C\delta_F\right].
            \end{equation}
            
        \end{enumerate}
    \end{remarks}

    \subsubsection{Micro-Macro correction}

    By applying the method of \textit{slow-fast} decomposition to derive a numerical technique based on micro-macro decomposition, we obtain a novel numerical method along with its corresponding error bounds for numerical approximation.

    \begin{Th} \label{theorem:highly_oscillatory_general_Micro-Macro}
        Let us consider $p\in\NN^*$ and denote the following learning errors:

        \begin{enumerate}[label=\textbf{(\roman*)}]
            \item \textbf{Learning error for high oscillation generator:}

            \begin{equation}
                \delta_{\phi,+} := \lnorm \frac{\phi^{[p]} - \phi_{\theta,+}}{\eps} \rnorm_{W^{1,\infty}(\Omega \times [0,2\pi]),\Linfty([0,\eps_+])} \\
            \end{equation}
            
            \item \textbf{Learning error for modified averaged field:}

            \begin{equation}
                \delta_{F} := \lnorm \widetilde{F^{[p],[q]}_h} - F_{\theta} \rnorm_{\Linfty([0,2\pi]\times[0,\eps_+])}
            \end{equation}

            \item \textbf{Learning error for field of associated to micro part:}

            \begin{equation}
                \delta_g := \lnorm g - g_{\theta} \rnorm_{\Linfty([0,2\pi]\times\Omega\times\Omega\times[0,\eps_+])}
            \end{equation}    
        \end{enumerate}

        Let $\Phi_{t,h}$ denote a numerical method of order $p$ (depending on time $t$). Consider the following two hypotheses:

        \begin{enumerate}[label = \textbf{(\roman*)}]
            \item For $f_1,f_2:[0,T]\times\Omega \rightarrow \RR^d$ sufficiently smooth and $t\in[0,T]$, we have, for all $h>0$,
    
            \begin{equation}
                \lnorm \Phi_{t,h}^{f_1} - \Phi_{t,h}^{f_2} \rnorm_{\Linfty(\Omega)} \leqslant Ch\lnorm f_1 - f_2 \rnorm_{\Linfty(\Omega)}
            \end{equation}
            
            for some positive constant $C>0$ independent of $f_1$ and $f_2$.
        
            \item For $f: [0,T]\times\Omega \rightarrow \RR^d$ sufficiently smooth and $t\in[0,T]$, there exists $L_f > 0$ s.t. for all $y_1,y_2\in\Omega$ and for all $h>0$, we have
    
            \begin{equation}
                \left| \Phi_{t,h}^f(y_1) - \Phi_{t,h}^f(y_2) \right| \leqslant (1+L_fh)|y_1-y_2|
            \end{equation}
        \end{enumerate}

        Then, there exist constants $\alpha_\phi, \lambda, \mu , M, \beta > 0$ (independent of $h,\eps)$ s.t. for all $h \leqslant h_+$ and $\eps \in ]0,1]$:

        \begin{eqnarray}
            \underset{0 \leqslant n \leqslant N}{Max}\left| y^{\eps}(t_n) - y^{\eps}_{\theta,n} \right| & \leqslant & \delta_{\phi,+}\eps + \alpha_\phi\frac{e^{\lambda T}-1}{\lambda}\left[\overline{C}h^{p+q-1}+C\delta_F\right] \\
            & + & \frac{e^{\mu T}-1}{\mu}\left[M'h^p + \frac{e^{\lambda T}-1}{\lambda}\beta(\overline{C}h^{p+q-1}+C\delta_F) + C\delta_g\right] \nonumber
        \end{eqnarray}
    \end{Th}

    \begin{remarks}
    \begin{enumerate}[label = \textbf{(\roman*)}]
        \item For some nonnegative constants $L_f, \alpha_{\theta,F}$ and $\alpha_\phi$ independent of $h, \eps$, we get:

        \begin{equation}
            \delta_g \leqslant (1+L_f+\alpha_{\theta,F}\eps)\delta_{\phi} + \alpha_\phi\delta_F
        \end{equation}
        
        \item Although the flow can be considered for non-autonomous differential equations and depends on time, the previous estimates concerning the numerical flow are assumed to be independent of time (except for the integration time $T$).
        
        \item The norm $\lnorm \cdot \rnorm_{W^{1,\infty}(\Omega \times [0,2\pi]),\Linfty([0,\eps_+])}$ corresponds to $W^{1,\infty}$-Sobolev norm w.r.t. the variables $y$ and $\tau$ and $\Linfty$-norm w.r.t. the variable $\eps$.
    \end{enumerate}
\end{remarks}

    \subsubsection{Alternative method}

    Alternative method for autonomous case yields to specific numerical error bounds.

    \begin{Th}\label{theorem:autonomous_highly_oscillatory}
        Let consider and denote the following learning errors:

        \begin{enumerate}[label=\textbf{(\roman*)}]
            \item \textbf{Learning error for flow $\phi^\eps$ (high oscillations):}

            \begin{equation}
                \delta_\Phi := \underset{(\tau,y,\eps)\in[0,2\pi]\times\Omega\times[0,\eps_+]}{Max}\left| \Phi_\tau^\eps(y) - \phi_\theta(\tau,y,\eps) \right|.
            \end{equation}

            \item \textbf{Learning error for flow $\pphi^{g^\eps}$:}

            \begin{equation}
                \delta_\pphi := \underset{(y,h,\eps)\in\Omega\times[0,h_+]\times[0,\eps_+]}{Max}\left| \frac{\pphi_h^{g^\eps}(y) - \pphi_\theta(y,h,\eps)}{h} \right|.
            \end{equation}
        \end{enumerate}

        Then there exist positive constants $\lambda, L$ and $C_T$ such that, for all $\eps \in ]0,1]$:

        \begin{equation}
            \underset{0 \leqslant n \leqslant N}{Max}\left| \yeps_{\theta,n} - \yeps(t_n) \right| \leqslant \delta_\Phi + L\frac{e^{\lambda_\theta T }-1}{\lambda}\delta_\pphi + C_Te^{-\frac{C_T}{\eps}}
        \end{equation}
    \end{Th}

    \section{Numerical experiments}\label{Paper2_section:numerical_experiments}

    To illustrate our theoretical results, we tested the method outlined in subsection \ref{subsection:General_strategy} on a simple dynamical system from physics:

    \begin{enumerate}[label = \textbf{\arabic*.}]
        \item \textbf{Inverted Pendulum:} This dynamical system describes the evolution of an unstable pendulum, with its center of gravity above the pivot point, undergoing forced oscillations. It is governed by the following equation:

        \begin{equation}
            \left\{\begin{array}{c c l}
                \dot{\yeps_1}(t) & = & \yeps_2(t) + \sin\left(\frac{t}{\eps}\right)\sin\left(\yeps_1(t)\right) \\
                \dot{\yeps_2}(t) & = & \sin\left(\yeps_1(t)\right) - \frac{1}{2}\sin\left(\frac{t}{\eps}\right)^2\sin\left(2\yeps_1(t)\right) - \sin\left(\frac{t}{\eps}\right)\cos\left(\yeps_1(t)\right)\yeps_2(t)
            \end{array}\right.
        \end{equation}
        
        and the average field associated to this equation is given by
    
        \begin{equation}
            \avgf(y) = \begin{bmatrix}
                y_2 \\
                \sin(y_1) - \frac{1}{4}\sin(2y_1)
            \end{bmatrix}
        \end{equation}

        \item \textbf{Van der Pol oscillator:} This system models an electrical circuit with nonlinear damping. It is governed by the following two-dimensional system:

        \begin{eqnarray}
            \left\{\begin{array}{c}
                 \dot{q} \\
                 \dot{p}
            \end{array} \right. & \begin{array}{c}
                 = \\
                 =
            \end{array} & \begin{array}{c}
                 \frac{1}{\eps}p \\
                 -\frac{1}{\eps}q + \left(\frac{1}{4}-q^2\right)p
            \end{array}
        \end{eqnarray}
        
        by performing the variable change $(y_1,y_2) \longmapsto S\left(\frac{t}{\eps}\right)(q,p)$, where $\tau \longmapsto S(\tau)$ is given by
        
        \begin{equation}
            \tau \longmapsto S(\tau) = \begin{bmatrix}
                \cos(\tau) & -\sin(\tau) \\
                \sin(\tau) & \cos(\tau)
            \end{bmatrix} \label{Variable_change_VDP}
        \end{equation}
        
        we get the system:
        
        \small\begin{equation}
            \left\{ \begin{array}{c c c}
            \dot{y_1}(t) & = & -\sin\left( \frac{t}{\eps}\right)\left[ \frac{1}{4} - \left( y_1(t)\cos\left(\frac{t}{\eps}\right) + y_2(t)\sin\left(\frac{t}{\eps}\right) \right)^2  \right]\left[ -y_1(t)\sin\left(\frac{t}{\eps}\right) + y_2(t)\cos\left( \frac{t}{\eps} \right) \right] \\
            \dot{y_2}(t) & = & \cos\left( \frac{t}{\eps}\right)\left[ \frac{1}{4} - \left( y_1(t)\cos\left(\frac{t}{\eps}\right) + y_2(t)\sin\left(\frac{t}{\eps}\right) \right)^2  \right]\left[ -y_1(t)\sin\left(\frac{t}{\eps}\right) + y_2(t)\cos\left( \frac{t}{\eps} \right) \right]
            \end{array} \right.
        \end{equation}
        
        \normalsize which is the classical form for highly oscillatory systems. The average field of the system in its canonical form is given, for all  $y\in\RR^2$, by

        \begin{equation}
            \avgf(y) = \frac{1}{8}(1-|y|^2)y.
        \end{equation}
    \end{enumerate}

    For a highly oscillatory differential equation of the form \eqref{equation:highly_oscillatory}, one can approximate the averaged field  $\avgdf$ and the highly oscillatory generator $\phi^\eps$ (defined by the formal power series \eqref{variable_change_averaged_field_formal}) using the following approximation sequences \cite{AVERAGING_chartier2020new}:

    \begin{equation}
        \left\{ \begin{array}{c c l}
            \phi_\tau^{[0]}(y) & = & y \\
            \phi_\tau^{[k+1]}(y) & = & y + \eps\displaystyle\int_0^\tau f\left(\sigma,\phi_\sigma^{[k]}(y)\right) - \frac{\pt \phi_\sigma^{[k]}}{\pt y}(y)F^{[k]}(y)\dd\sigma 
        \end{array} \right. \label{sequence_approximation_variable_change}
    \end{equation}
    
    \noindent and

    \begin{equation}
        \left\{ \begin{array}{c c l}
            F^{[0]}(y) & = & \avgf(y) \\
            F^{[k]}(y) & = & \left(\frac{\pt \langle \phi_\cdot^{[k]}\rangle} {\pt y}(y)\right)^{-1}\Big\langle f\left(\cdot , \phi_\cdot^{[k]}(y)\right) \Big\rangle.
        \end{array} \right. \label{sequence_approximation_averaged_field}
    \end{equation}

    These sequences can be used to obtain an approximation of arbitrary order for the formal power series \eqref{variable_change_averaged_field_formal}

    \begin{equation}
        \phi_\tau^{\eps}(y) = \phi_\tau^{[k]}(y) + \mathcal{O}\left( \eps^{k+1} \right) \text{  and } \avgdf(y) = F^{[k]}(y) + \mathcal{O}\left( \eps^{k+1} \right)
    \end{equation}

    \subsection{Approximation of the averaged field and high oscillation generator}

    In this subsection, we investigate the approximation error between the learned averaged field and the high-oscillation generator, compared to the theoretical averaged field and the high-oscillation generator, both at orders 0 and 1 for the inverted pendulum. We examine the learning error with respect to the high-oscillation parameter $\eps$. Specifically, we compute the values

    \begin{equation}
        \underset{y\in\Omega}{Max}\left| F_{\theta}\left(y,0,\eps\right) - F^{[k]}(y) \right| \text{  and } \underset{(\tau,y) \in [0,2\pi] \times \Omega}{Max}\left| \phi_{\theta,+}\left(\tau,y,\eps\right) - \phi^{[k]}_\tau(y) \right| \label{formula:learning error_epsilon}
    \end{equation}
    
    \noindent where $\phi^{[k]}$ and $F^{[k]}$ are computed using the formulas \eqref{sequence_approximation_variable_change} and \eqref{sequence_approximation_averaged_field} for $k=0,1$ respectively.\\

\begin{figure}[H]
		\centering
		\begin{minipage}{0.45\linewidth}
			\includegraphics[width=\linewidth]{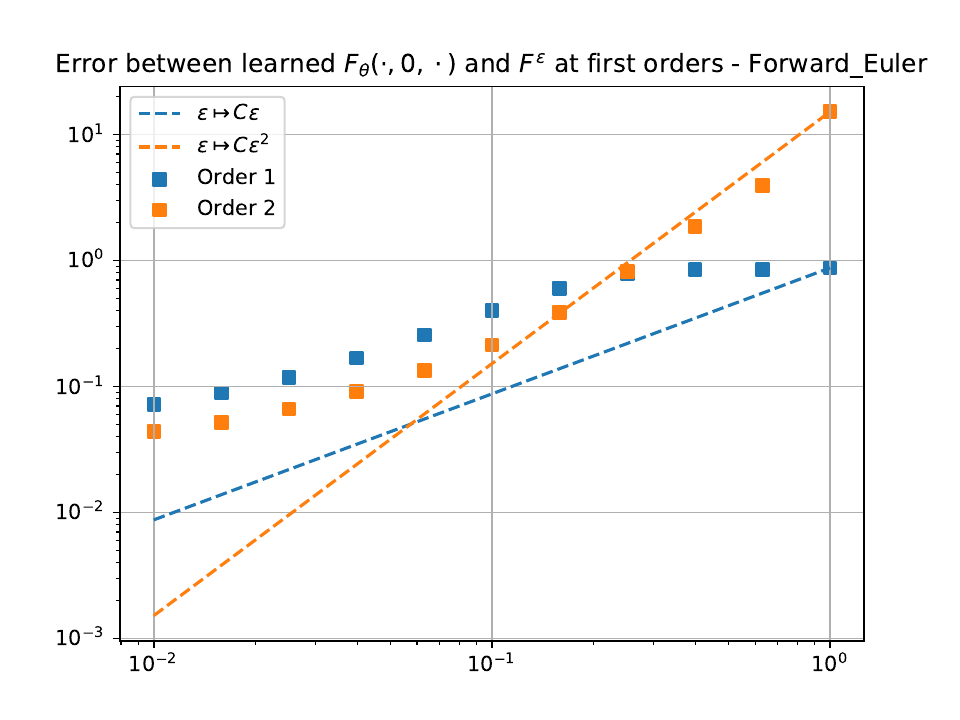}
		\end{minipage}
		\begin{minipage}{0.45\linewidth}
			\includegraphics[width=\linewidth]{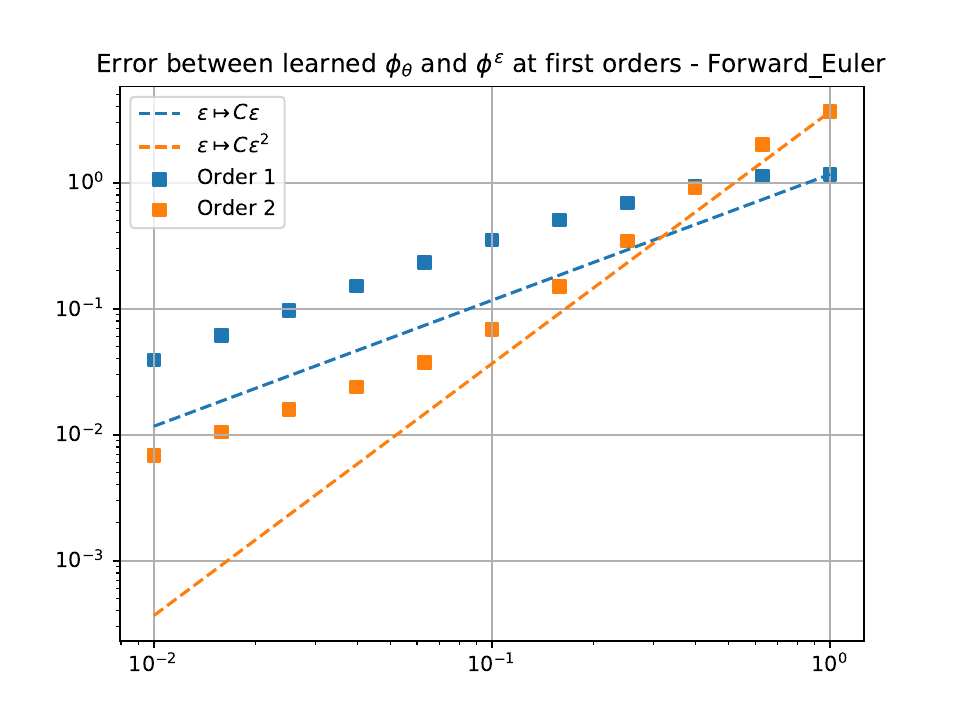}
		\end{minipage}
		\caption{Inverted pendulum with forward Euler method. Left: Error between $F^{[k]}$ and $F_{\theta}(\cdot,0,\eps)$ for $k = 0,1$ w.r.t. $\eps$. Right: Error between $\phi^{[k]}$ and $\phi_{\theta,+}(\cdot,\cdot,\eps)$ for $k = 0,1$ w.r.t. $\eps$. }
		\label{fig_learning_error_Inverted_Pendulum_Forward_Euler}
    \end{figure}

\begin{figure}[H]
		\centering
		\begin{minipage}{0.45\linewidth}
			\includegraphics[width=\linewidth]{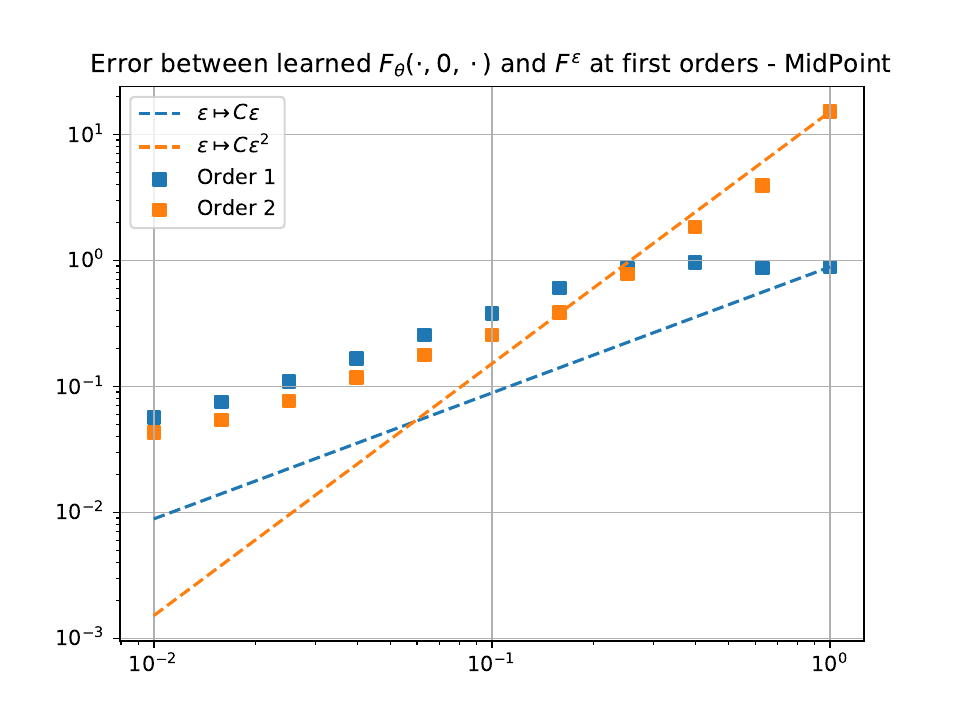}
		\end{minipage}
		\begin{minipage}{0.45\linewidth}
			\includegraphics[width=\linewidth]{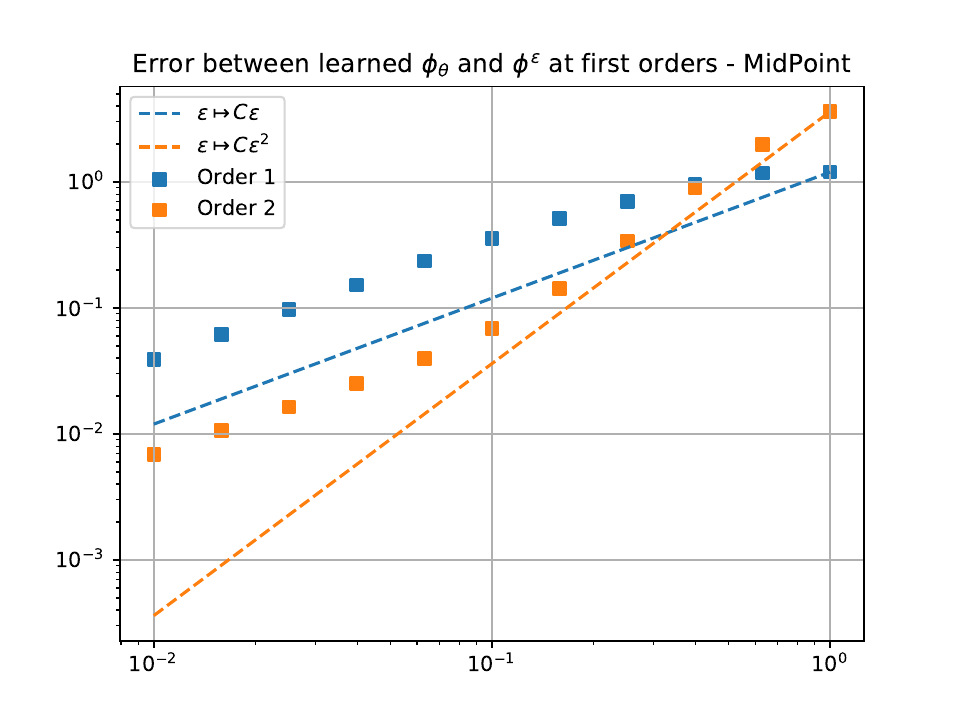}
		\end{minipage}
		\caption{Inverted pendulum with midpoint method. Left: Error between $F^{[k]}$ and $F_{\theta}(\cdot,0,\eps)$ for $k = 0,1$ w.r.t. $\eps$. Right: Error between $\phi^{[k]}$ and $\phi_{\theta,+}(\cdot,\cdot,\eps)$ for $k = 0,1$ w.r.t. $\eps$. }
		\label{fig_learning_error_Inverted_Pendulum_MidPoint}
    \end{figure}

\begin{figure}[H]
		\centering
		\begin{minipage}{0.45\linewidth}
			\includegraphics[width=\linewidth]{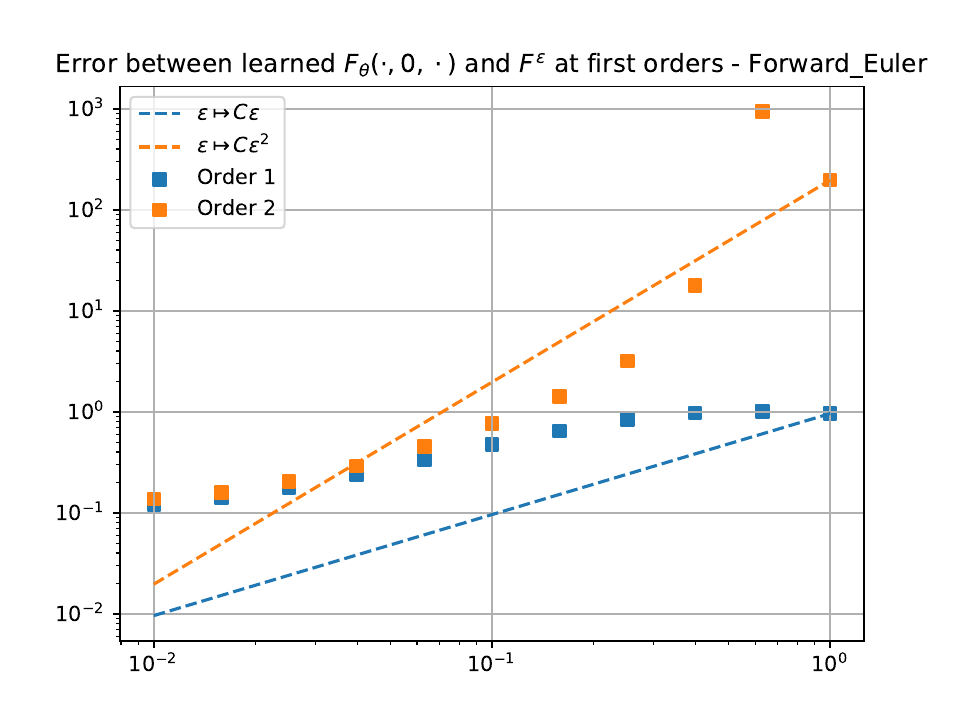}
		\end{minipage}
		\begin{minipage}{0.45\linewidth}
			\includegraphics[width=\linewidth]{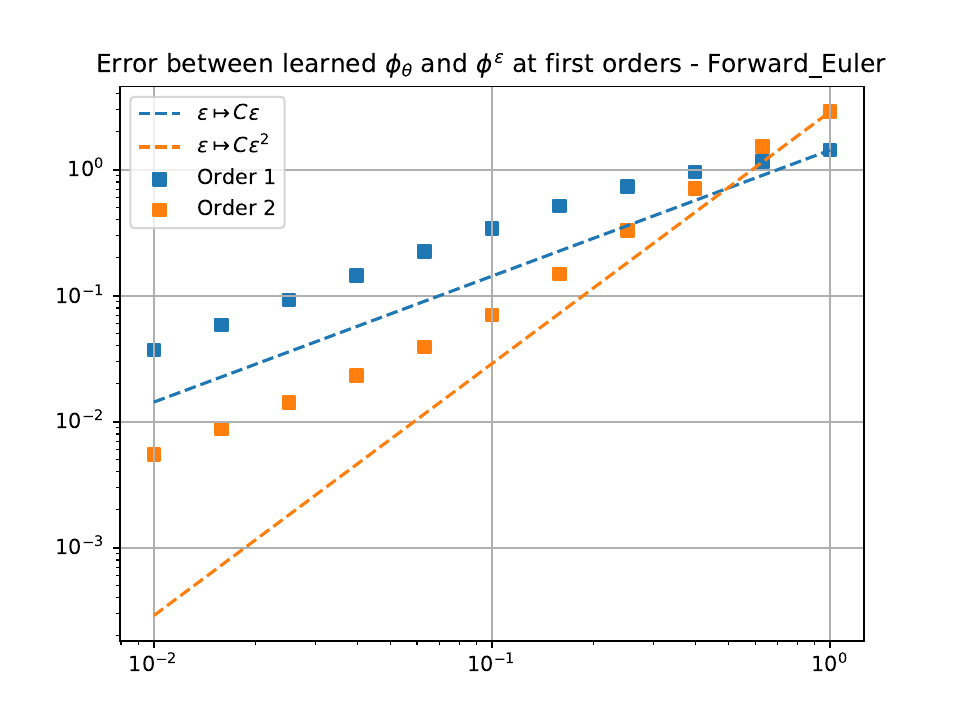}
		\end{minipage}
		\caption{Van der Pol oscillator with forward Euler method. Left: Error between $F^{[k]}$ and $F_{\theta}(\cdot,0,\eps)$ for $k = 0,1$ w.r.t. $\eps$. Right: Error between $\phi^{[k]}$ and $\phi_{\theta,+}(\cdot,\cdot,\eps)$ for $k = 0,1$ w.r.t. $\eps$. }
		\label{fig_learning_error_VDP_Forward_Euler}
    \end{figure}

    The figures \ref{fig_learning_error_Inverted_Pendulum_Forward_Euler}, \ref{fig_learning_error_Inverted_Pendulum_MidPoint} and \ref{fig_learning_error_VDP_Forward_Euler} confirm that the modified field can be effectively learned with our neural network. However, the learning error seems to slow the decay of the error bounds.
    
    \subsection{Loss decay and integration of ODE's}\label{subsection:Loss_decay_integration_ODE}

    To compare the integration of a dynamical system using the learned modified averaged field and the learned high-oscillation generator, we will solve the inverted pendulum problem using the Forward Euler and Midpoint methods. Before doing so, we examine the decay of the loss functions for the training and test data sets ($Loss_{Train}$ and $Loss_{Test}$). Their similarity provides a good indication that there is no overfitting, suggesting that the size of the training data set has been appropriately estimated. Since the MSE loss is used, it offers an insight into the square of the error for both the equation and the auto-encoder structure.\\

\begin{figure}[H]
        \centering
        \includegraphics[scale = 0.6]{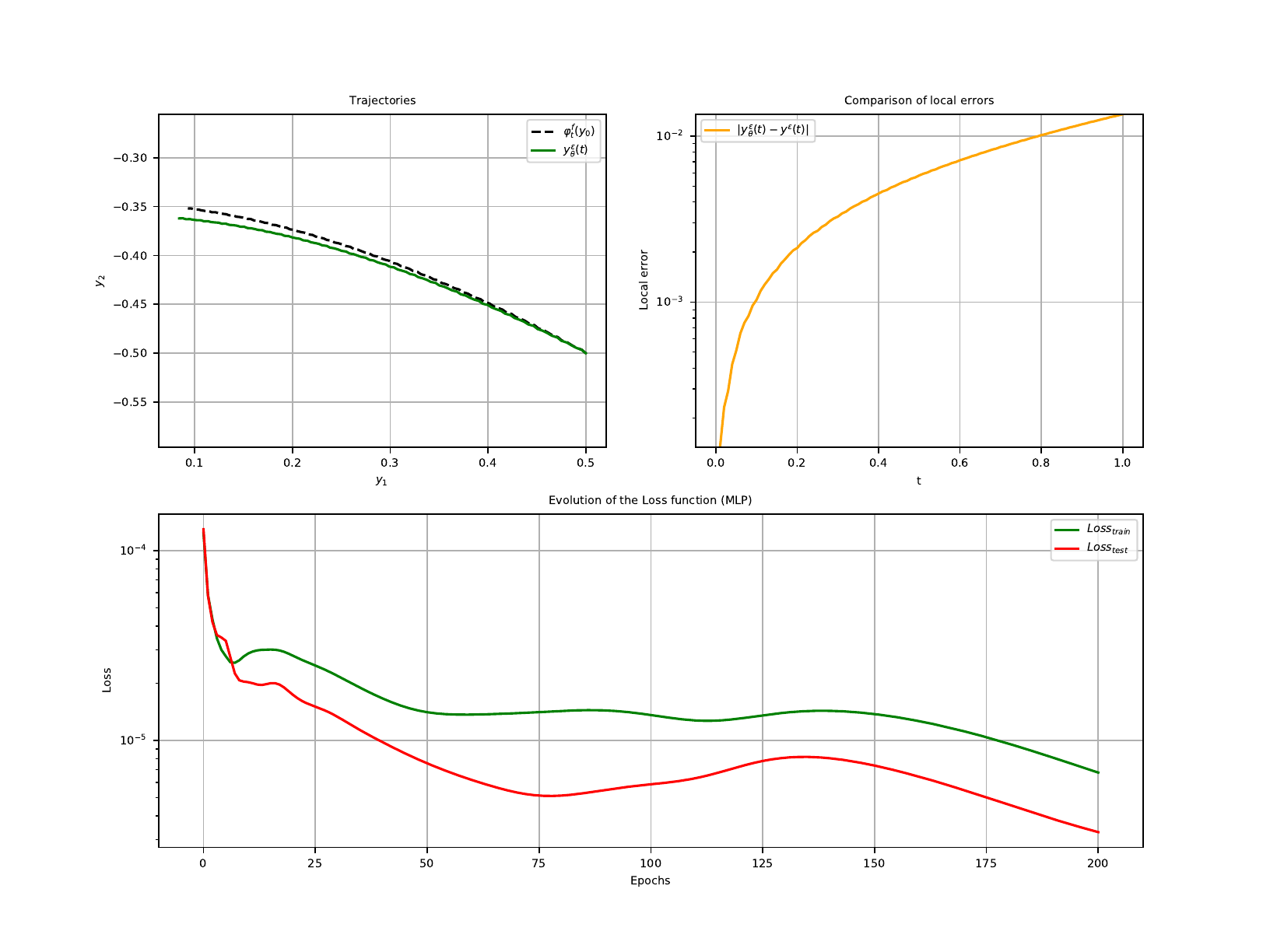}
        \caption{Comparison between $Loss$ decays (green: $Loss_{Train}$, red: $Loss_{Test}$), trajectories (dashed dark: exact flow, green: numerical flow with learned vector fields and local error (yellow) for the inverted pendulum with Forward Euler method in the case $\eps = 0.001$.}
        \label{fig_inverted_pendulum_FEuler_eps=0.001}
    \end{figure}

\begin{figure}[H]
        \centering
        \includegraphics[scale = 0.6]{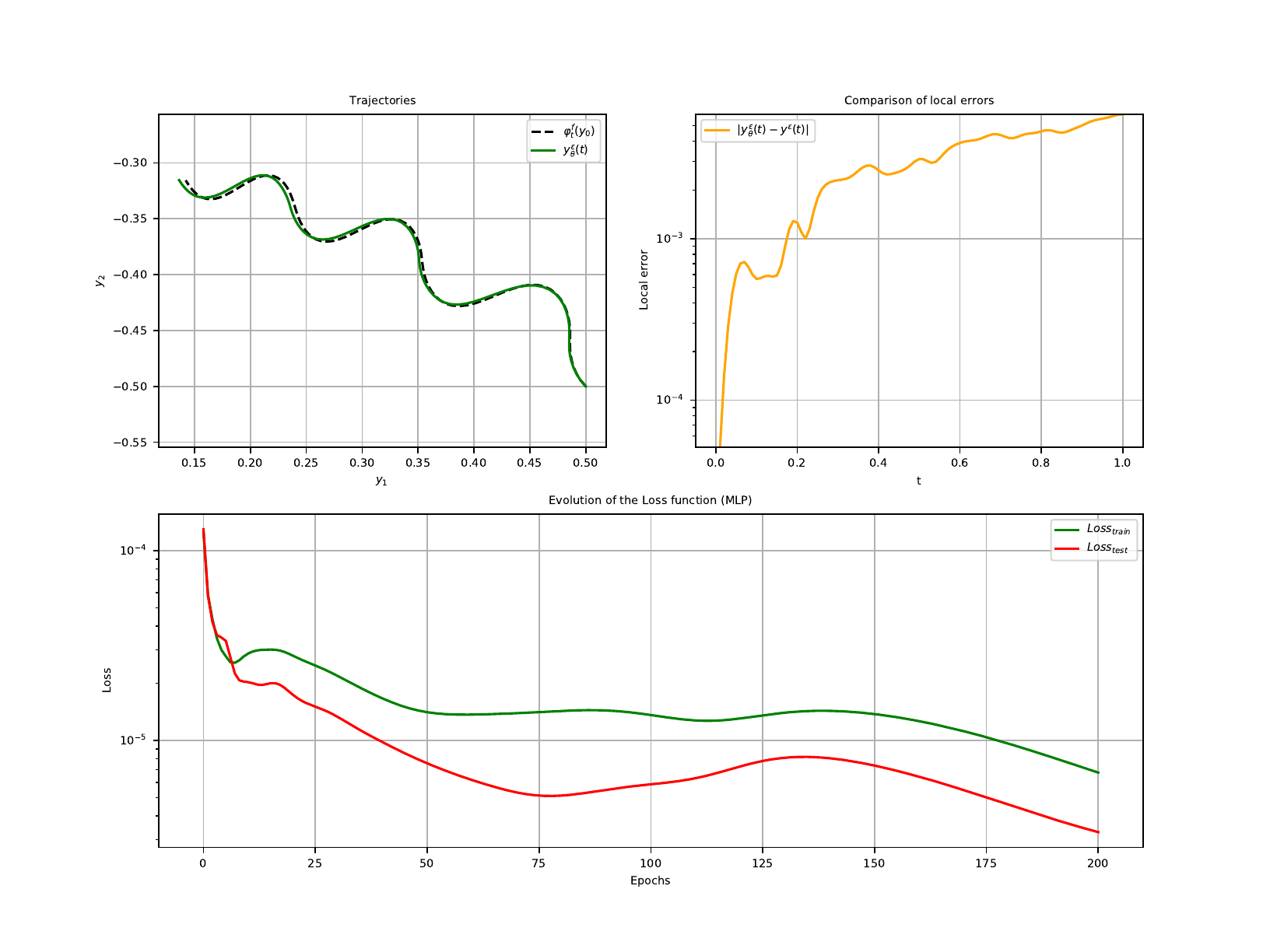}
        \caption{Comparison between $Loss$ decays (green: $Loss_{Train}$, red: $Loss_{Test}$), trajectories (dashed dark: exact flow, green: numerical flow with learned vector fields and local error (yellow) for the inverted pendulum with Forward Euler method in the case $\eps = 0.05$.}
        \label{fig_inverted_pendulum_FEuler_eps=0.05}
    \end{figure}

\begin{figure}[H]
        \centering
        \includegraphics[scale = 0.6]{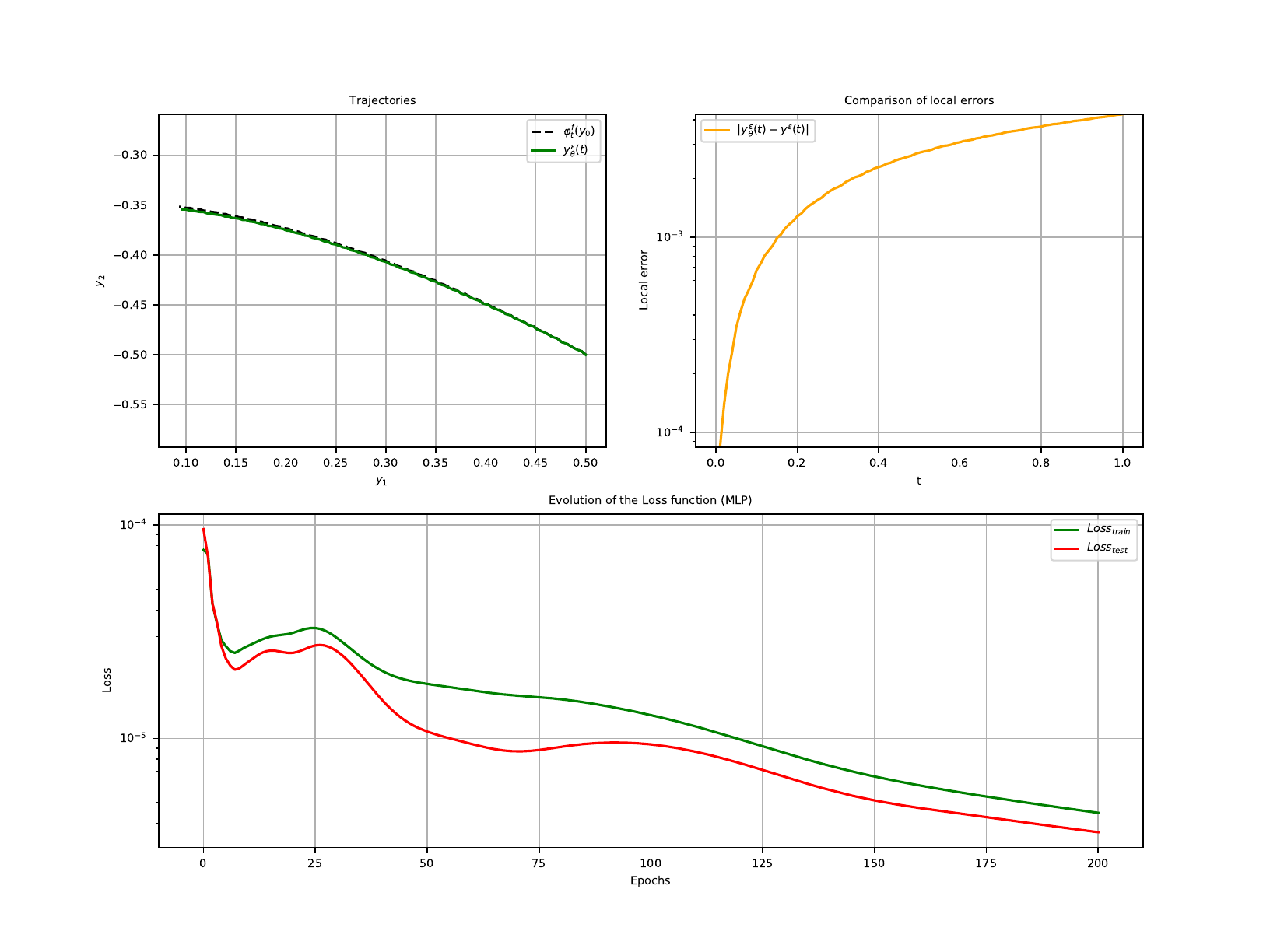}
        \caption{Comparison between $Loss$ decays (green: $Loss_{Train}$, red: $Loss_{Test}$), trajectories (dashed dark: exact flow, green: numerical flow with learned vector fields and local error (yellow) for the inverted pendulum with midpoint method in the case $\eps = 0.001$.}
        \label{fig_inverted_pendulum_MidPoint_eps=0.001}
    \end{figure}

\begin{figure}[H]
        \centering
        \includegraphics[scale = 0.6]{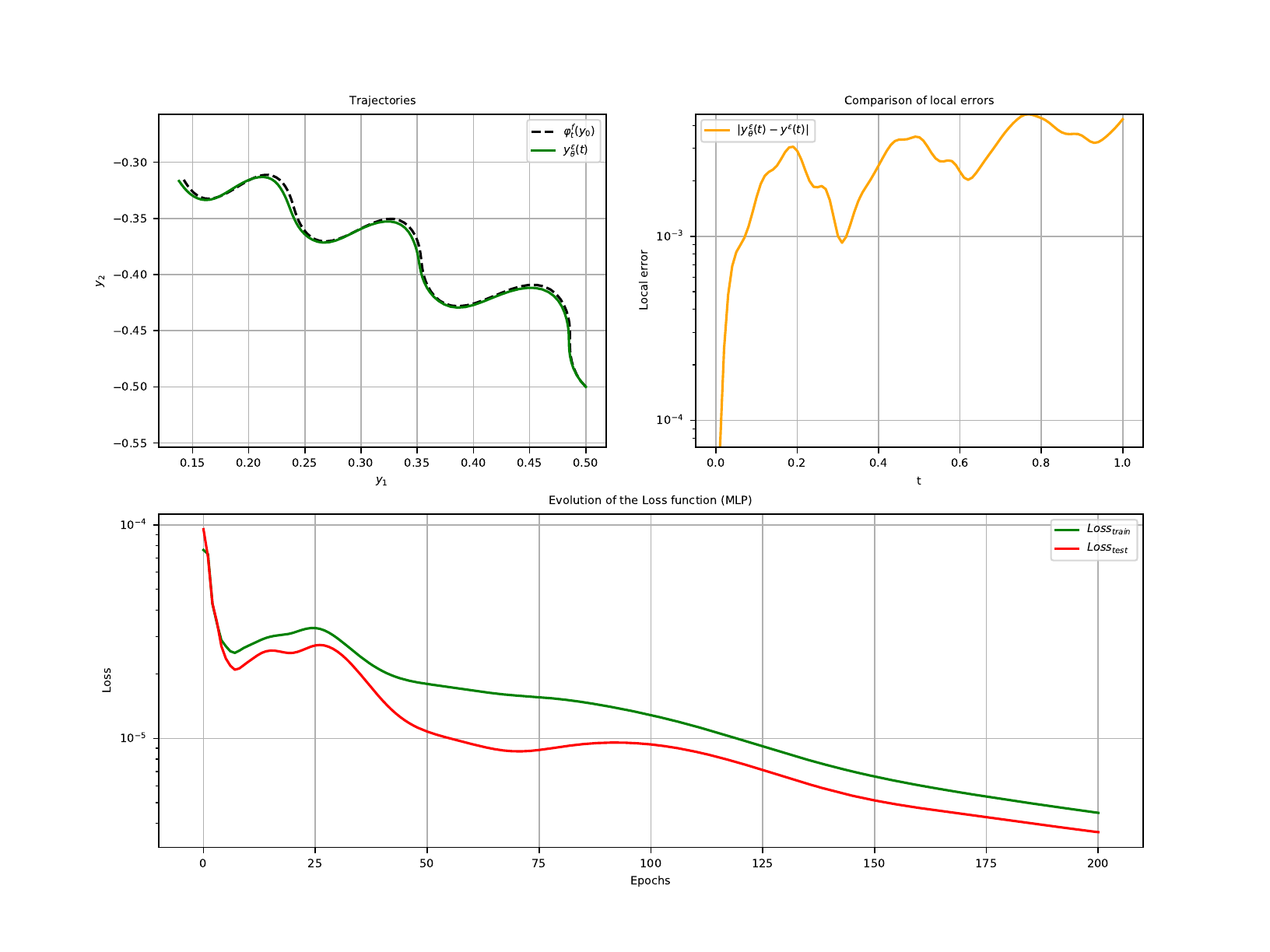}
        \caption{Comparison between $Loss$ decays (green: $Loss_{Train}$, red: $Loss_{Test}$), trajectories (dashed dark: exact flow, green: numerical flow with learned vector fields and local error (yellow) for the inverted pendulum with midpoint method in the case $\eps = 0.05$.}
        \label{fig_inverted_pendulum_MidPoint_eps=0.05}
    \end{figure}

\begin{figure}[H]
        \centering
        \includegraphics[scale = 0.6]{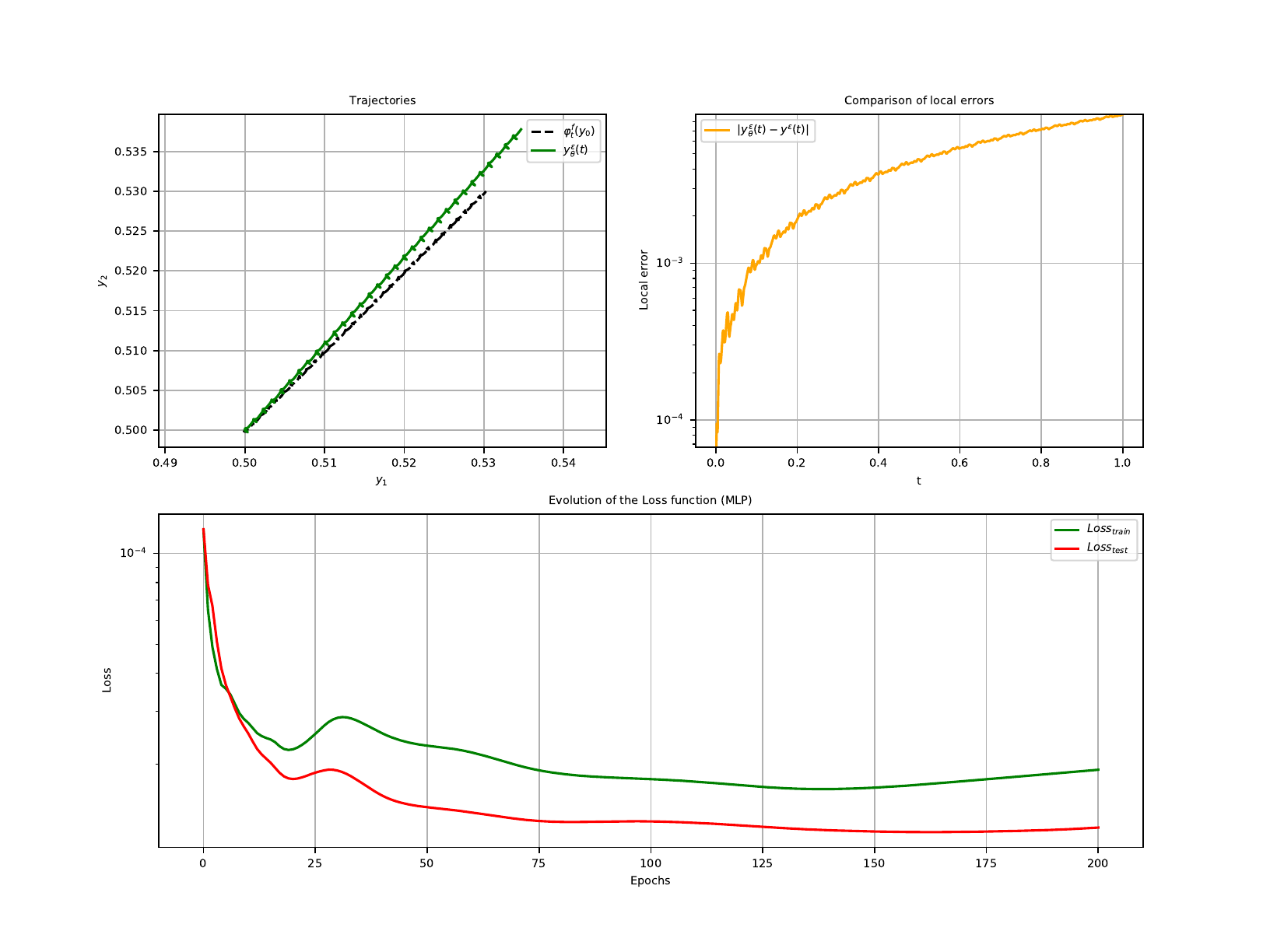}
        \caption{Comparison between $Loss$ decays (green: $Loss_{Train}$, red: $Loss_{Test}$), trajectories (dashed dark: exact flow, green: numerical flow with learned vector fields and local error (yellow) for the Van der Pol oscillator with Forward Euler method in the case $\eps = 0.01$.}
        \label{fig_VDP_FEuler_eps=0.01}
    \end{figure}

\begin{figure}[H]
        \centering
        \includegraphics[scale = 0.8]{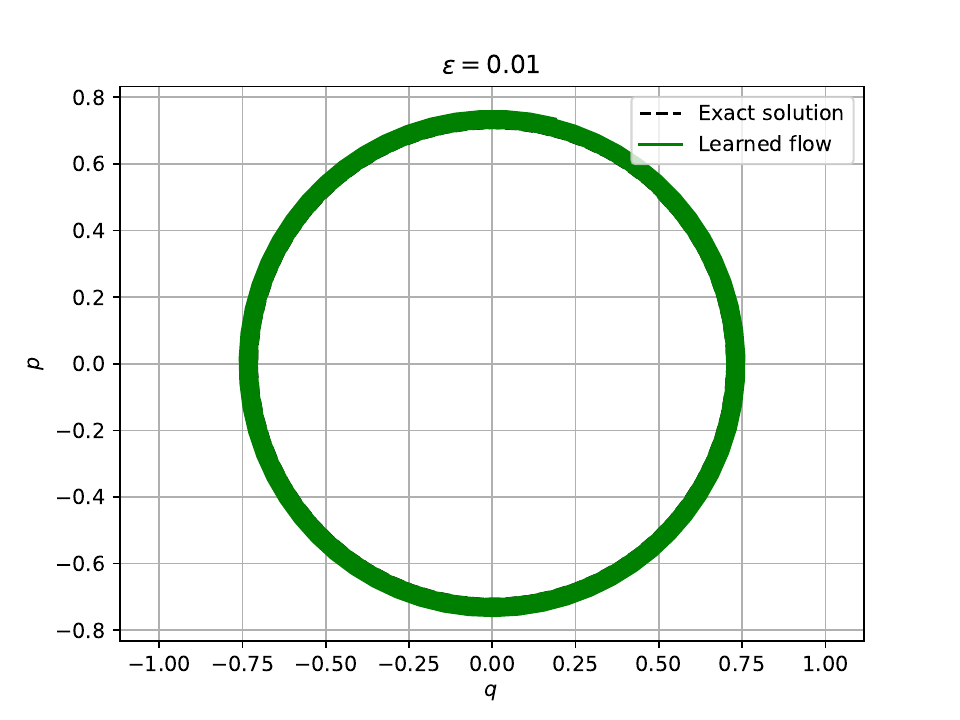}
        \caption{Comparison between trajectories (dashed dark: exact flow, green: numerical flow with learned vector field) for the Van der Pol oscillator with Forward Euler method in the case $\eps = 0.01$ after the inverse variable change \eqref{Variable_change_VDP}.}
        \label{fig_VDP_FEuler_Variable_change_eps=0.01}
    \end{figure}

\begin{figure}[H]
        \centering
        \includegraphics[scale = 0.6]{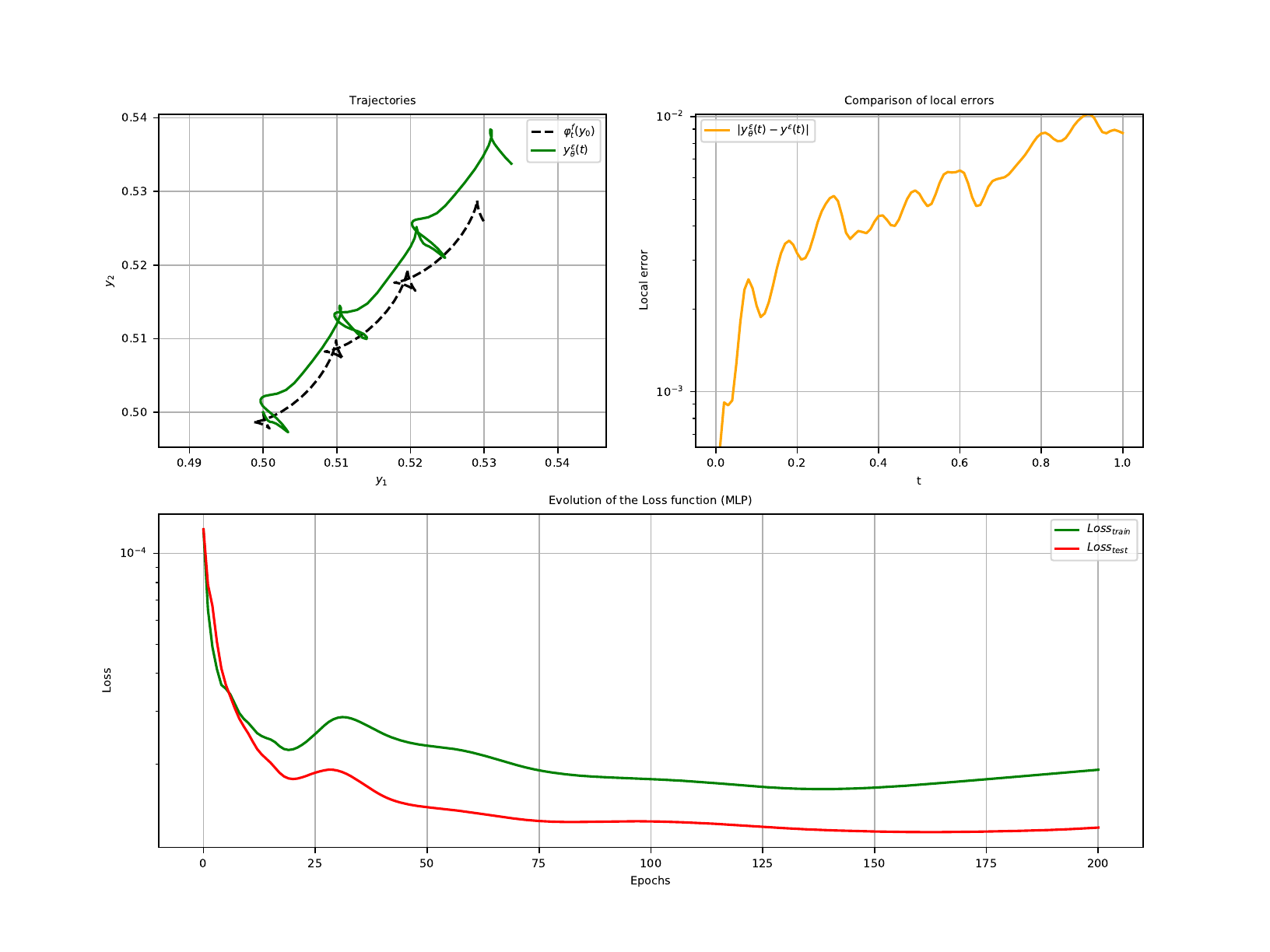}
        \caption{Comparison between $Loss$ decays (green: $Loss_{Train}$, red: $Loss_{Test}$), trajectories (dashed dark: exact flow, green: numerical flow with learned vector fields and local error (yellow) for the Van der Pol oscillator with Forward Euler method in the case $\eps = 0.1$.}
        \label{fig_VDP_FEuler_eps=0.1}
    \end{figure}

\begin{figure}[H]
        \centering
        \includegraphics[scale = 0.8]{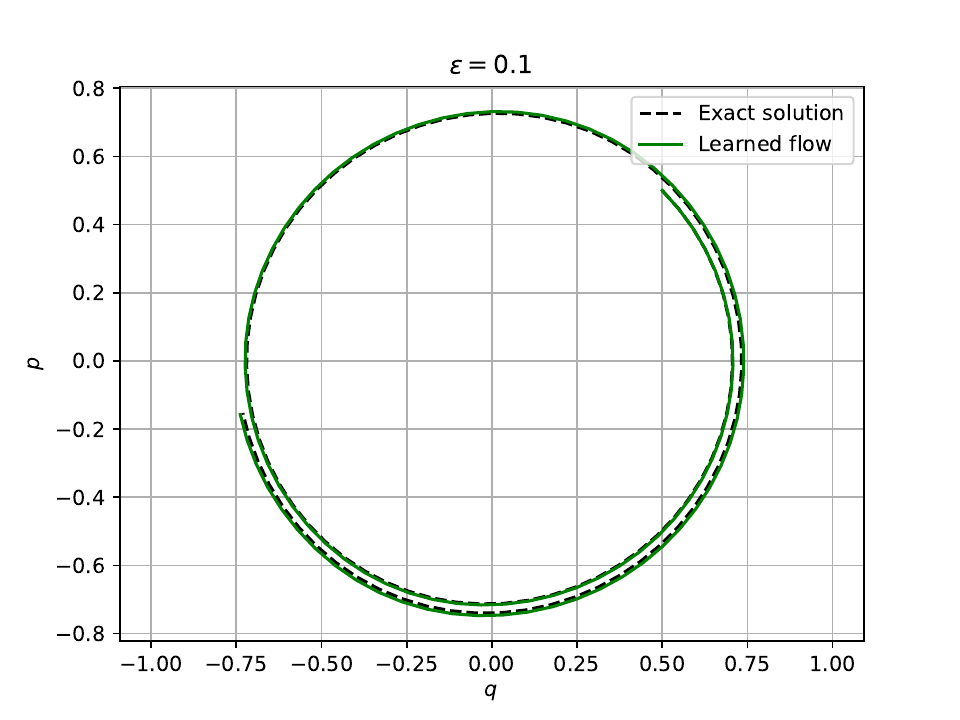}
        \caption{Comparison between trajectories (dashed dark: exact flow, green: numerical flow with learned vector field) for the Van der Pol oscillator with Forward Euler method in the case $\eps = 0.1$ after the inverse variable change \eqref{Variable_change_VDP}.}
        \label{fig_VDP_FEuler_Variable_change_eps=0.1}
    \end{figure}

    Figures \ref{fig_inverted_pendulum_FEuler_eps=0.001}, \ref{fig_inverted_pendulum_FEuler_eps=0.05}, \ref{fig_inverted_pendulum_MidPoint_eps=0.001}, and \ref{fig_inverted_pendulum_MidPoint_eps=0.05} demonstrate accurate numerical integration using the slow-fast decomposition-based method, applying the corresponding learned vector field for both the Forward Euler and Midpoint methods for the inverted pendulum. For the Van der Pol oscillator, although the learning for the transformed dynamical system appears to be less accurate (figures \ref{fig_VDP_FEuler_eps=0.01} and \ref{fig_VDP_FEuler_eps=0.1}), the integration after the inverse variable change \eqref{Variable_change_VDP} appears to be correct, as shown in figures \ref{fig_VDP_FEuler_Variable_change_eps=0.01} and \ref{fig_VDP_FEuler_Variable_change_eps=0.1}.\\

\begin{figure}[H]
        \centering
        \includegraphics[scale = 0.6]{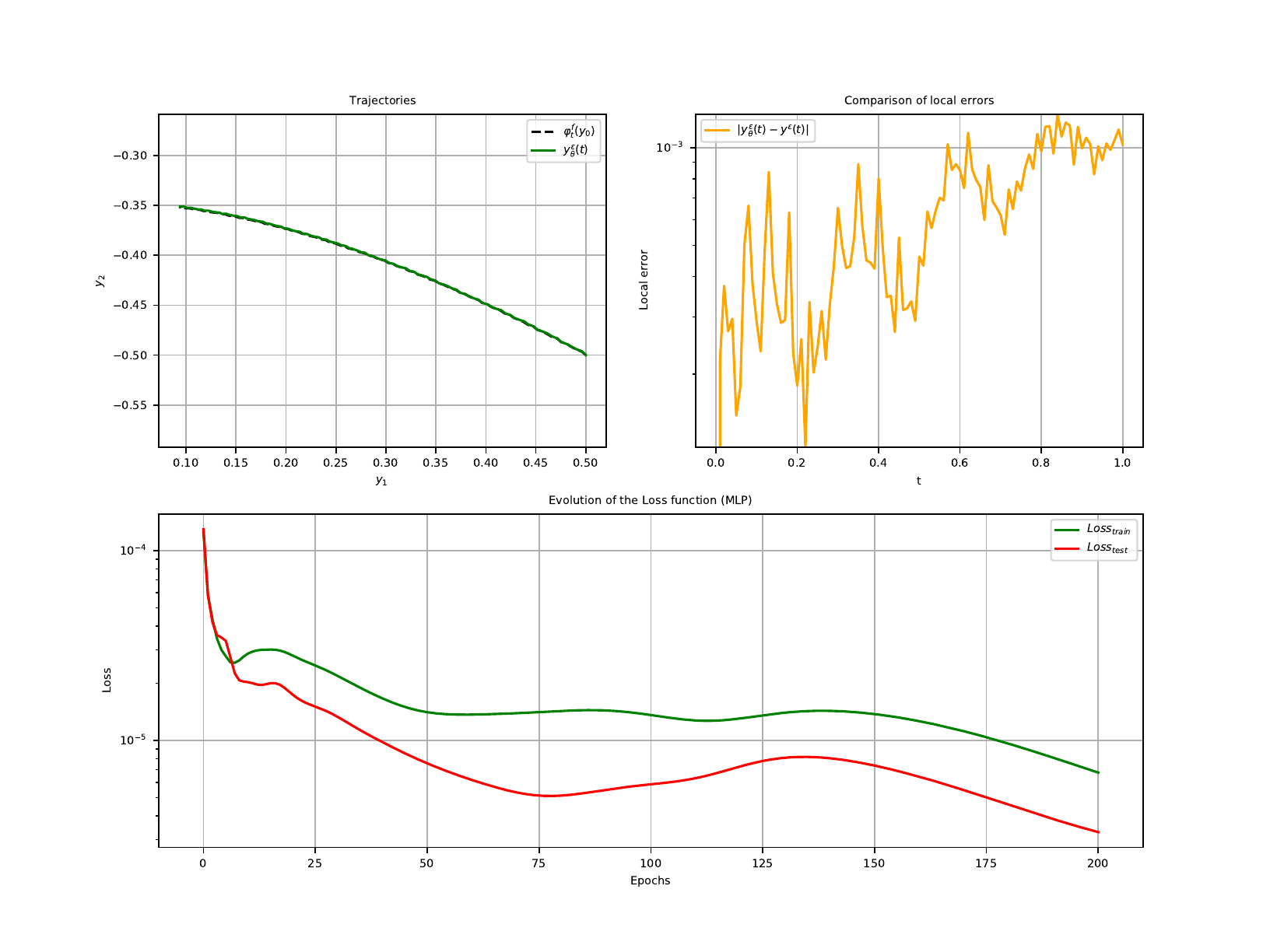}
        \caption{Comparison between $Loss$ decays (green: $Loss_{Train}$, red: $Loss_{Test}$), trajectories (dashed dark: exact flow, green: numerical flow with learned vector fields and local error (yellow) for the inverted pendulum with Forward Euler method with Micro-Macro correction in the case $\eps = 0.001$.}
        \label{fig_inverted_pendulum_FEuler_MM_eps=0.001}
    \end{figure}

\begin{figure}[H]
        \centering
        \includegraphics[scale = 0.6]{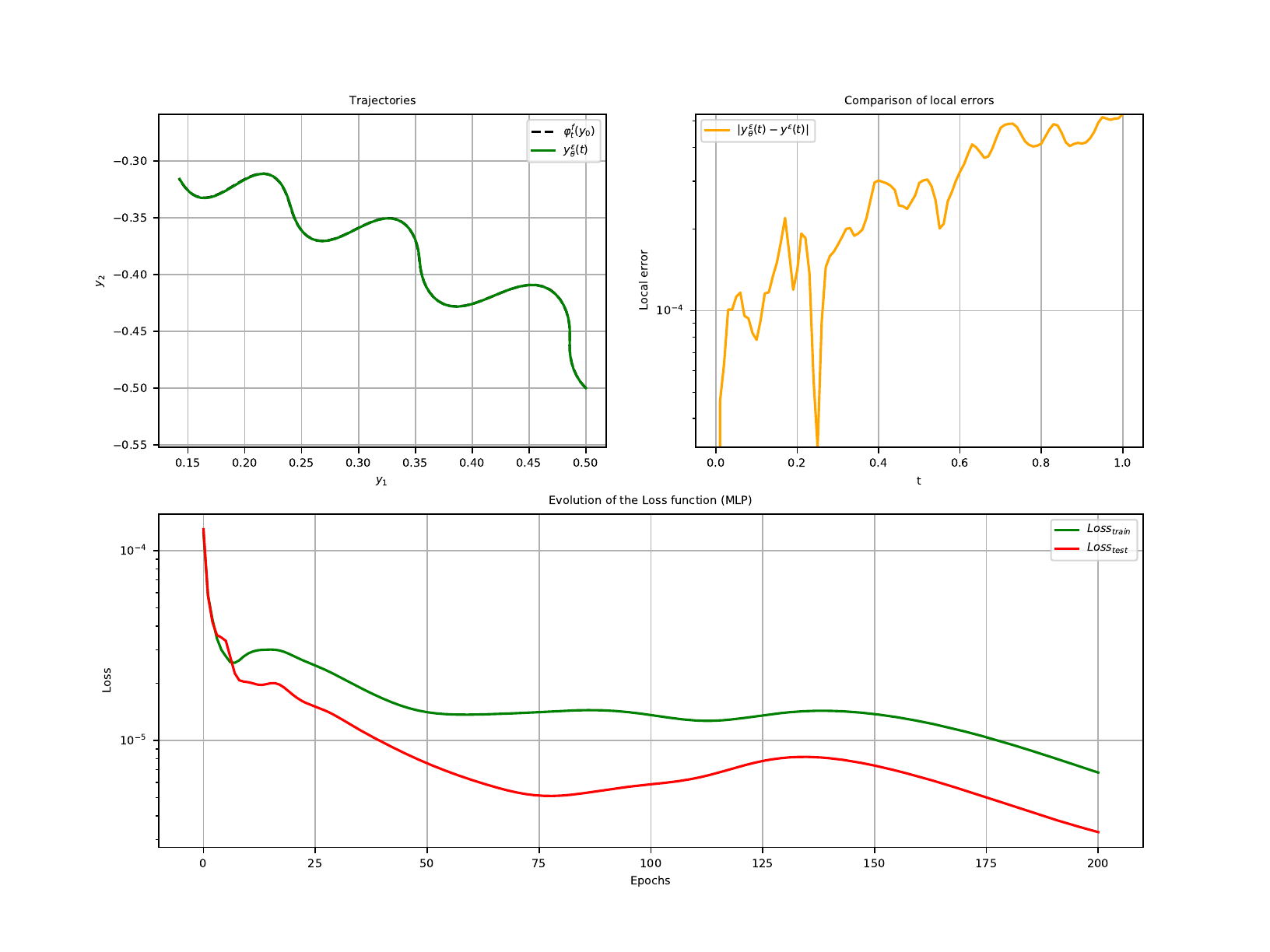}
        \caption{Comparison between $Loss$ decays (green: $Loss_{Train}$, red: $Loss_{Test}$), trajectories (dashed dark: exact flow, green: numerical flow with learned vector fields and local error (yellow) for the inverted pendulum with Forward Euler method with Micro-Macro correction in the case $\eps = 0.05$.}
        \label{fig_inverted_pendulum_FEuler_MM_eps=0.05}
    \end{figure}

\begin{figure}[H]
        \centering
        \includegraphics[scale = 0.6]{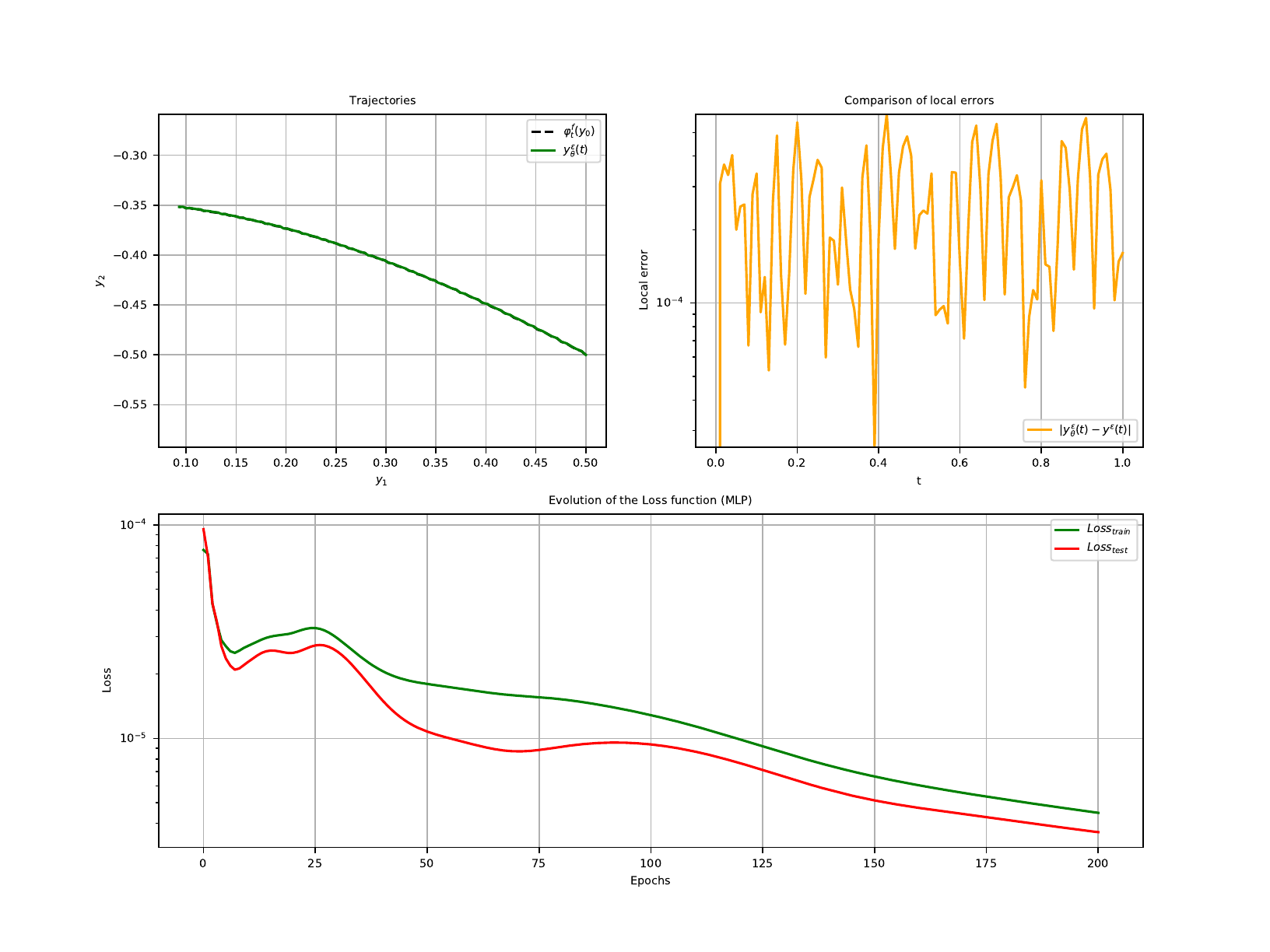}
        \caption{Comparison between $Loss$ decays (green: $Loss_{Train}$, red: $Loss_{Test}$), trajectories (dashed dark: exact flow, green: numerical flow with learned vector fields and local error (yellow) for the inverted pendulum with midpoint method with Micro-Macro correction in the case $\eps = 0.001$.}
        \label{fig_inverted_pendulum_MidPoint_MM_eps=0.001}
    \end{figure}

\begin{figure}[H]
        \centering
        \includegraphics[scale = 0.6]{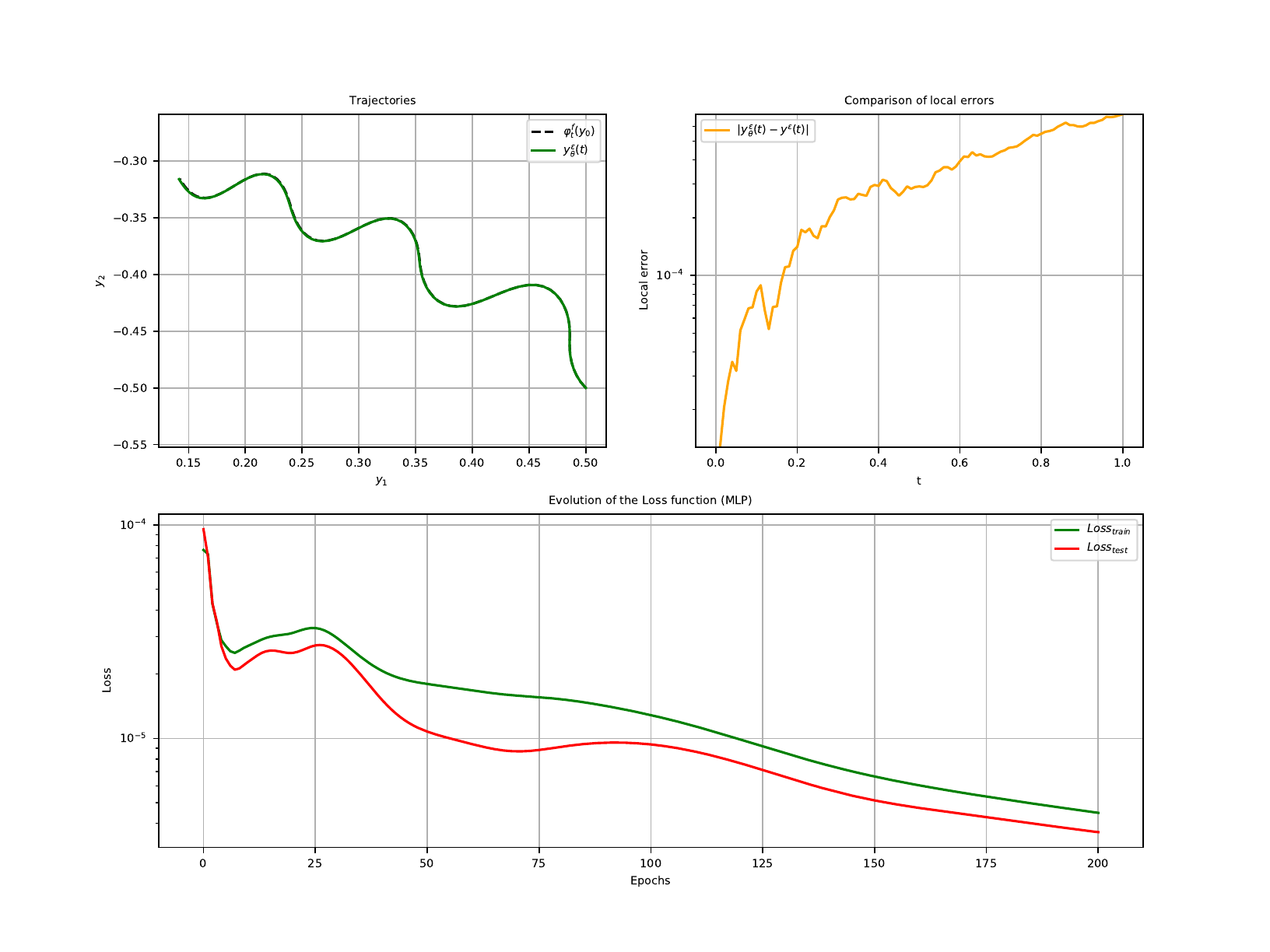}
        \caption{Comparison between $Loss$ decays (green: $Loss_{Train}$, red: $Loss_{Test}$), trajectories (dashed dark: exact flow, green: numerical flow with learned vector fields and local error (yellow) for the inverted pendulum with midpoint method with Micro-Macro correction in the case $\eps = 0.05$.}
        \label{fig_inverted_pendulum_MidPoint_MM_eps=0.05}
    \end{figure}

\begin{figure}[H]
        \centering
        \includegraphics[scale = 0.6]{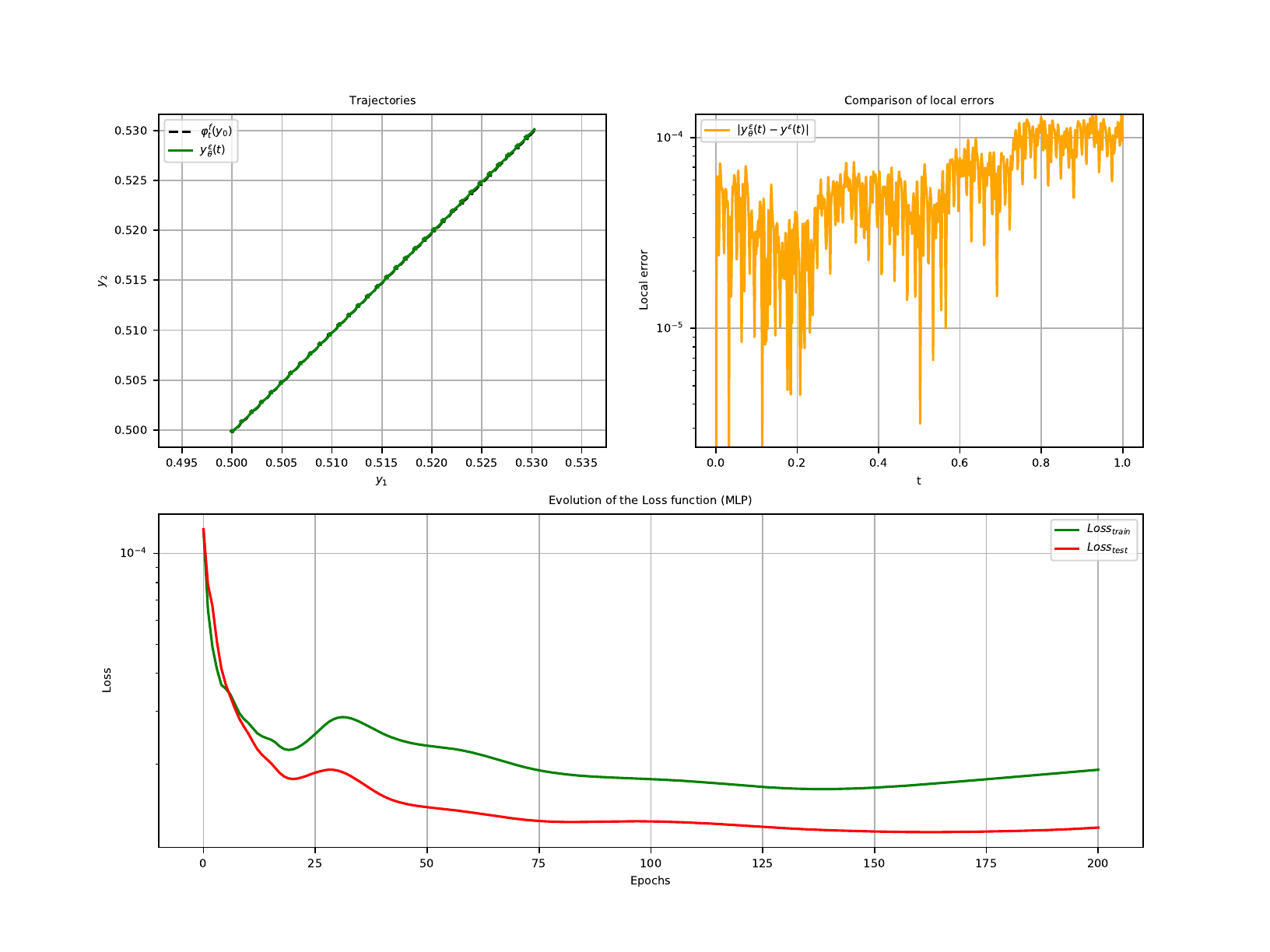}
        \caption{Comparison between $Loss$ decays (green: $Loss_{Train}$, red: $Loss_{Test}$), trajectories (dashed dark: exact flow, green: numerical flow with learned vector fields and local error (yellow) for the Van der Pol oscillator with Forward Euler method with Micro-Macro correction in the case $\eps = 0.01$.}
        \label{fig_VDP_FEuler_MM_eps=0.01}
    \end{figure}

\begin{figure}[H]
        \centering
        \includegraphics[scale = 0.8]{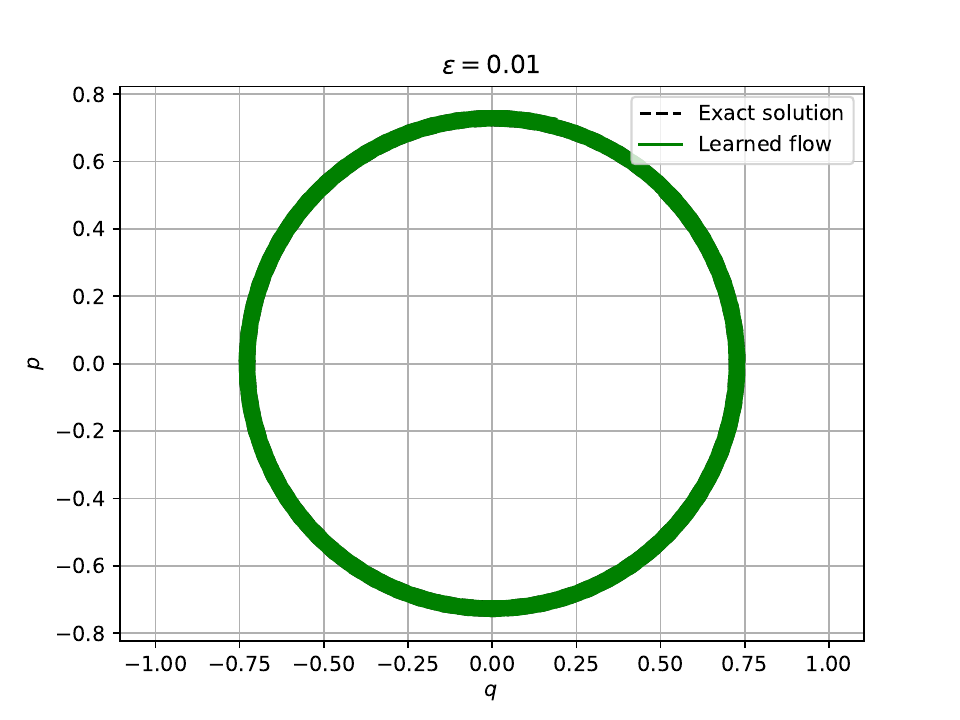}
        \caption{Comparison between trajectories (dashed dark: exact flow, green: numerical flow with learned vector field) for the Van der Pol oscillator with Forward Euler method with Micro-Macro correction in the case $\eps = 0.01$ after the inverse variable change \eqref{Variable_change_VDP}.}
        \label{fig_VDP_FEuler_MM_Variable_change_eps=0.01}
    \end{figure}

\begin{figure}[H]
        \centering
        \includegraphics[scale = 0.6]{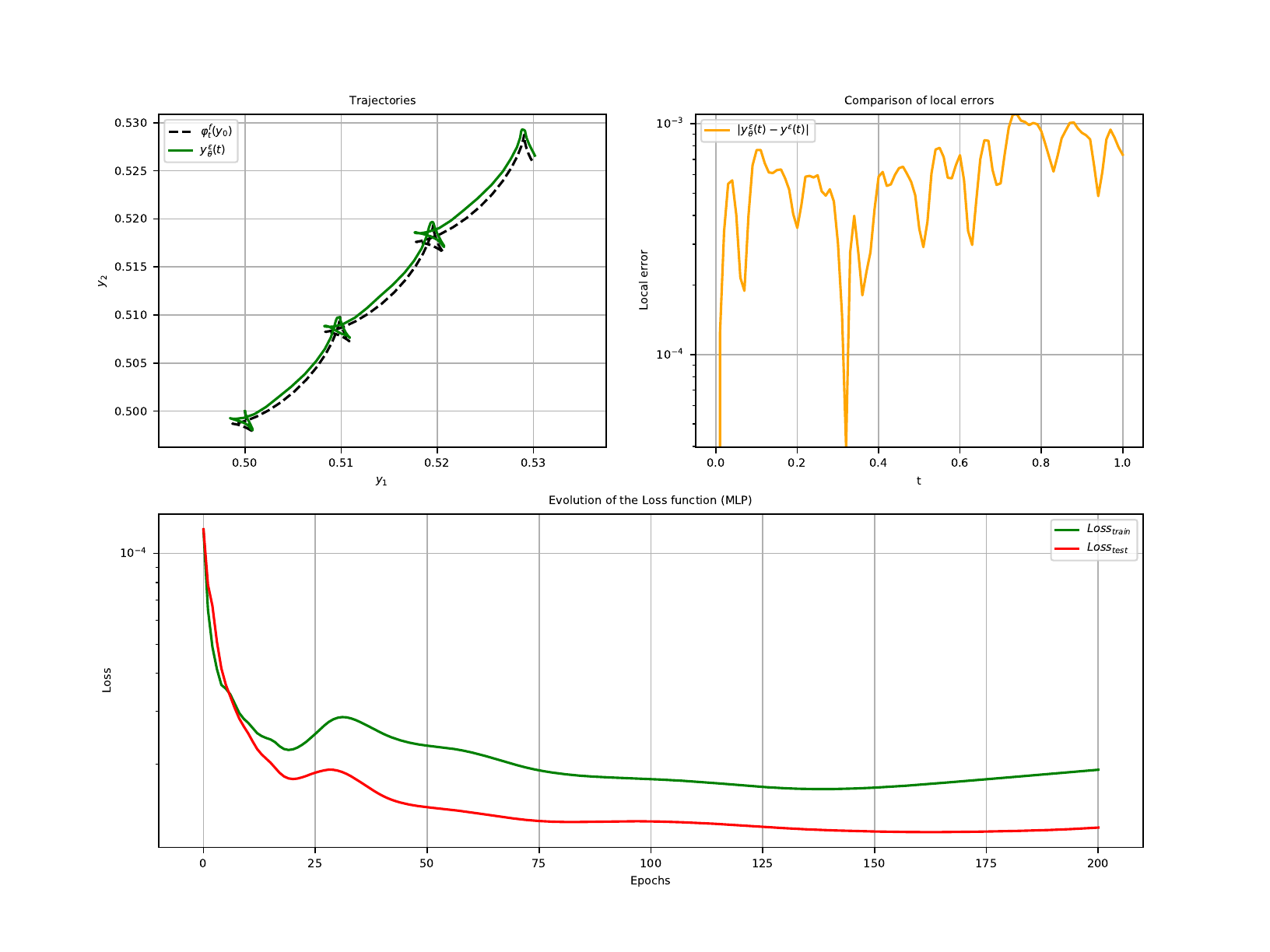}
        \caption{Comparison between $Loss$ decays (green: $Loss_{Train}$, red: $Loss_{Test}$), trajectories (dashed dark: exact flow, green: numerical flow with learned vector fields and local error (yellow) for the Van der Pol oscillator with Forward Euler method with Micro-Macro correction in the case $\eps = 0.1$.}
        \label{fig_VDP_FEuler_MM_eps=0.1}
    \end{figure}

\begin{figure}[H]
        \centering
        \includegraphics[scale = 0.8]{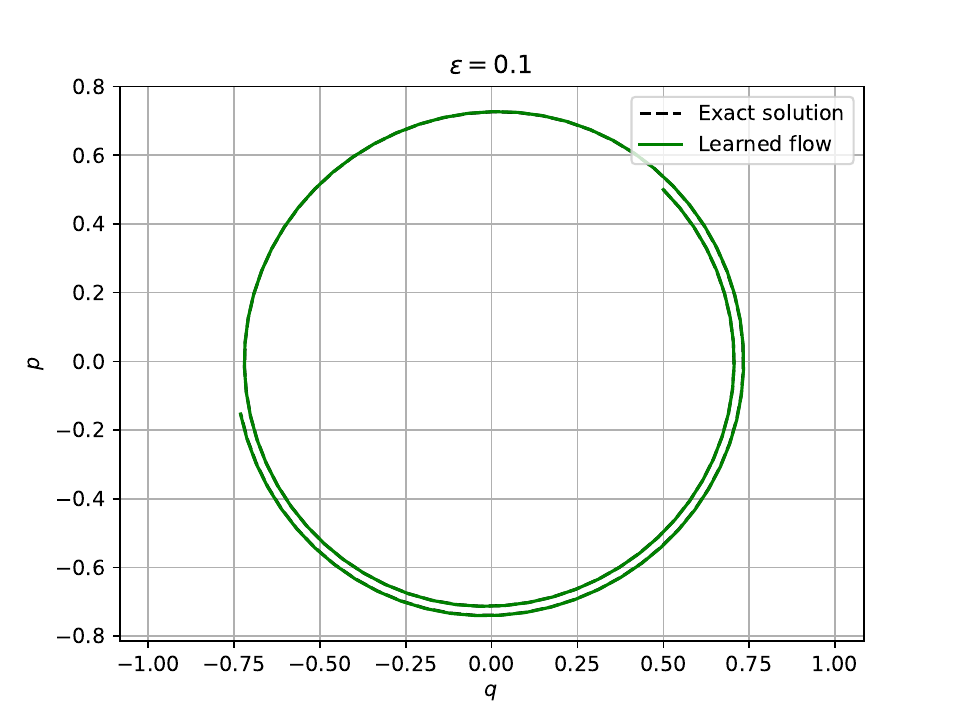}
        \caption{Comparison between trajectories (dashed dark: exact flow, green: numerical flow with learned vector field) for the Van der Pol oscillator with Forward Euler method with Micro-Macro correction in the case $\eps = 0.1$ after the inverse variable change \eqref{Variable_change_VDP}.}
        \label{fig_VDP_FEuler_MM_Variable_change_eps=0.1}
    \end{figure}

    However, incorporating the micro-macro correction significantly enhances the accuracy, especially for larger values of $\eps$. In particular, the numerical approximation of the solution for the inverted pendulum becomes more accurate, as observed in figures \ref{fig_inverted_pendulum_FEuler_MM_eps=0.001}, \ref{fig_inverted_pendulum_FEuler_MM_eps=0.05}, \ref{fig_inverted_pendulum_MidPoint_MM_eps=0.001}, and \ref{fig_inverted_pendulum_MidPoint_MM_eps=0.05}. For the Van der Pol oscillator, it can be seen that the system's solution is accurately approximated using the micro-macro correction, as shown in figures \ref{fig_VDP_FEuler_MM_eps=0.01} and \ref{fig_VDP_FEuler_MM_eps=0.1}, even before the variable change, as can be seen in figures \ref{fig_VDP_FEuler_MM_Variable_change_eps=0.01} and \ref{fig_VDP_FEuler_MM_Variable_change_eps=0.1}.


    \subsubsection{Inverted Pendulum - forward Euler method}

    \subsubsection{Inverted Pendulum - midpoint method}

    \subsubsection{Van der Pol oscillator - forward Euler method}

    
    \subsection{Error curves w.r.t. step size}

\begin{figure}[H]
		\centering
		\begin{minipage}{0.45\linewidth}
			\includegraphics[width=\linewidth]{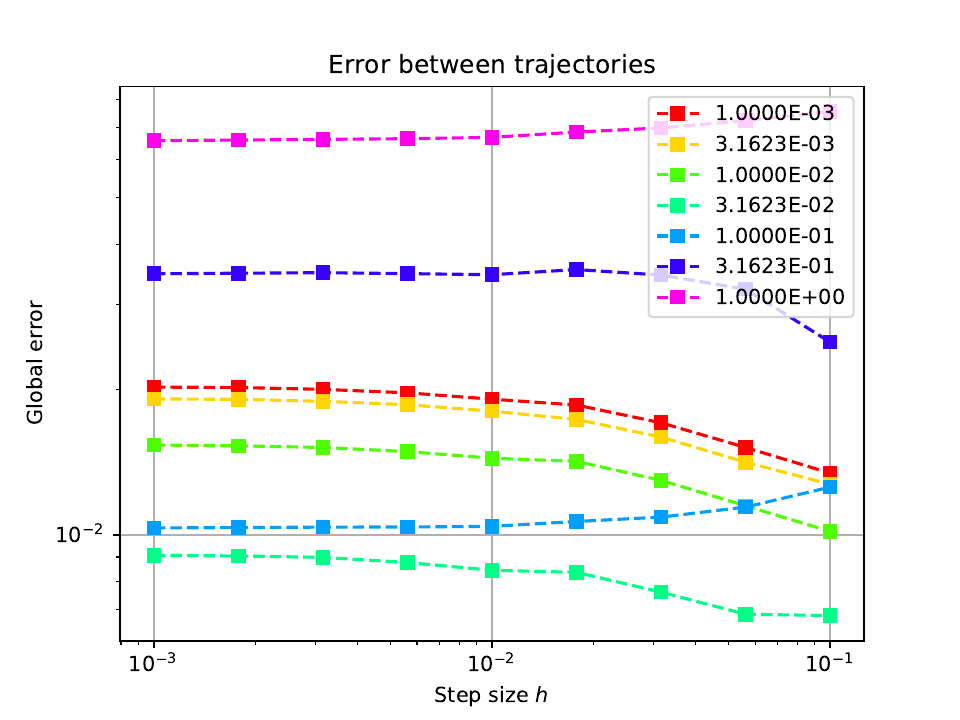}
		\end{minipage}
		\begin{minipage}{0.45\linewidth}
			\includegraphics[width=\linewidth]{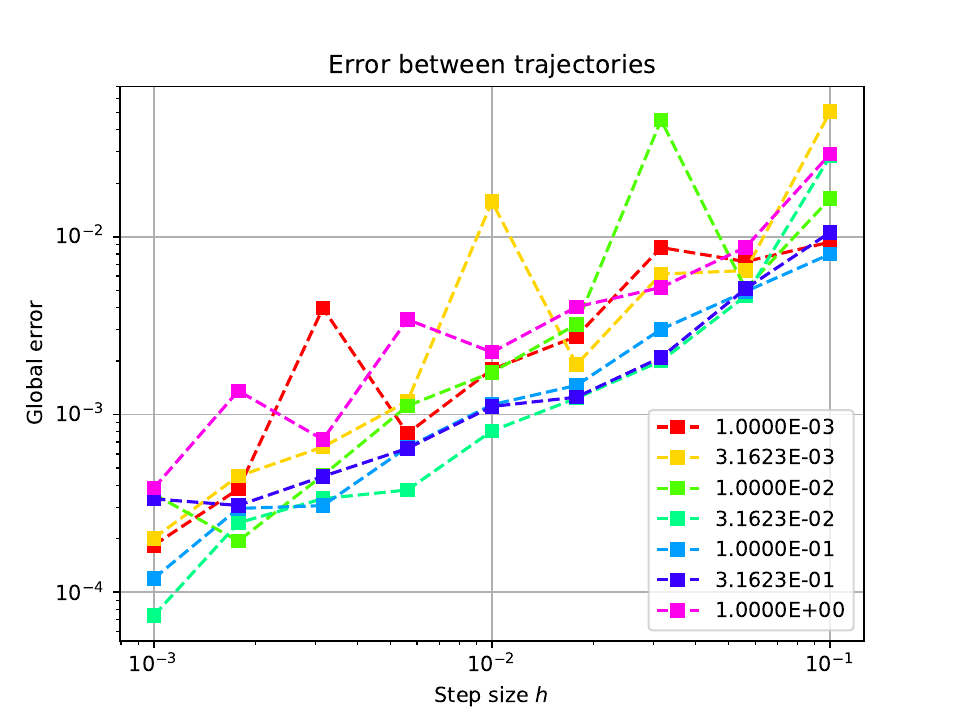}
		\end{minipage}
		\caption{Integration errors (each color corresponds to a high oscillation parameter $\eps$) of Inverted Pendulum with Forward Euler. Left: \textit{Slow-fast} decomposition-based method. Right: With Micro-Macro correction.}
		\label{fig_global_error_pendulum_FEuler}
    \end{figure}

\begin{figure}[H]
		\centering
		\begin{minipage}{0.45\linewidth}
			\includegraphics[width=\linewidth]{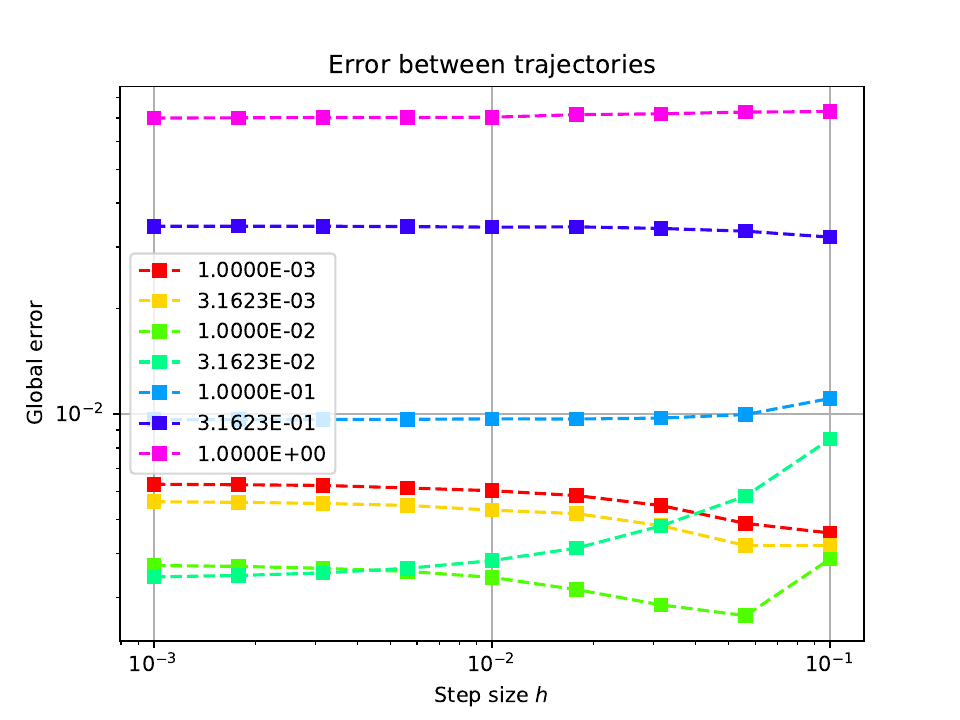}
		\end{minipage}
		\begin{minipage}{0.45\linewidth}
			\includegraphics[width=\linewidth]{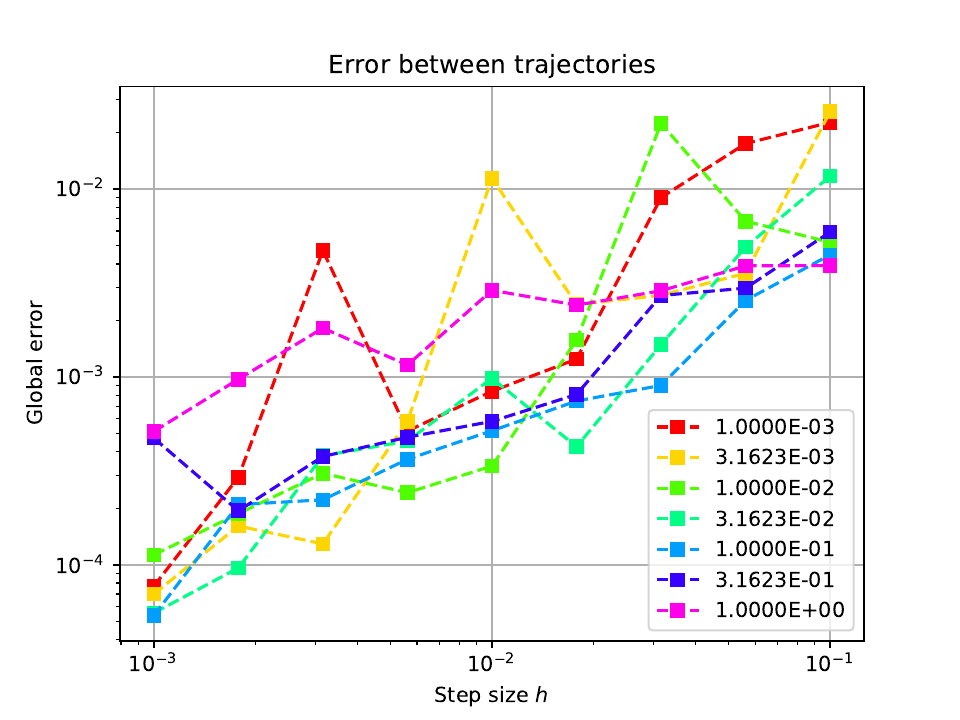}
		\end{minipage}
		\caption{Integration errors (each color corresponds to a high oscillation parameter $\eps$) of Inverted Pendulum with midpoint. Left: \textit{Slow-fast} decomposition-based method. Right: With Micro-Macro correction.}
		\label{fig_global_error_pendulum_MidPoint}
    \end{figure}

\begin{figure}[H]
		\centering
		\begin{minipage}{0.45\linewidth}
			\includegraphics[width=\linewidth]{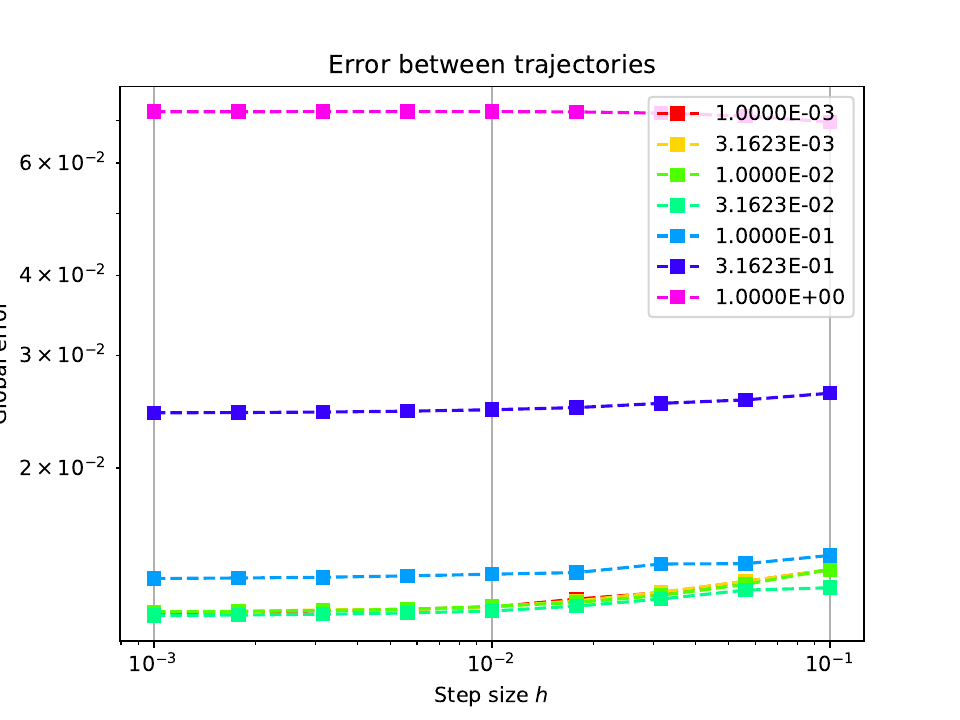}
		\end{minipage}
		\begin{minipage}{0.45\linewidth}
			\includegraphics[width=\linewidth]{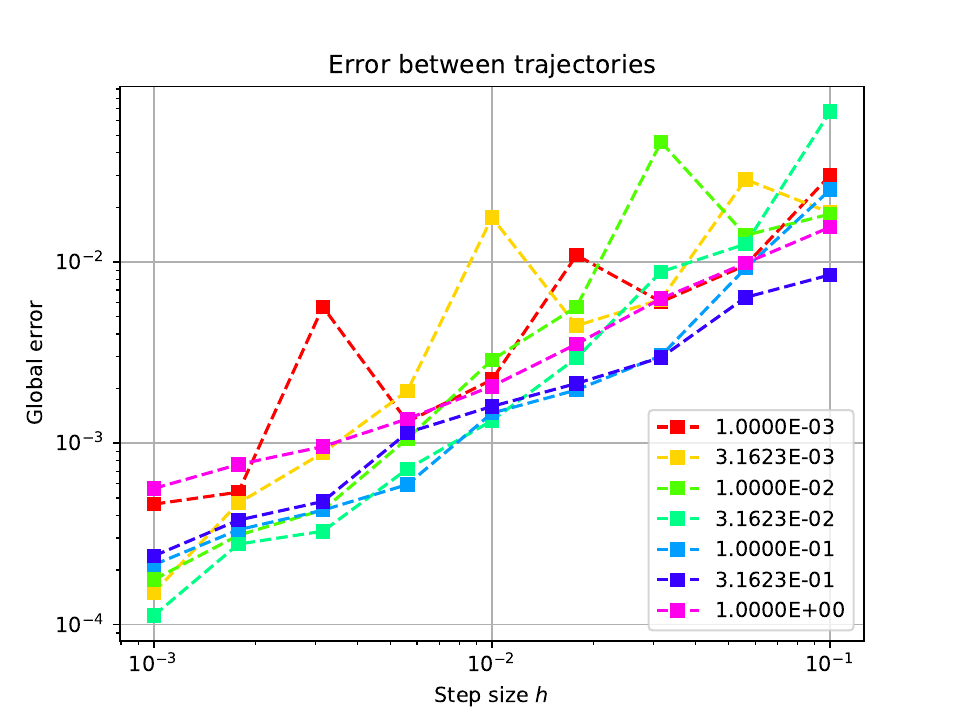}
		\end{minipage}
		\caption{Integration errors (each color corresponds to a high oscillation parameter $\eps$) of Van der Pol oscillator with Forward Euler. Left: \textit{Slow-fast} decomposition-based method. Right: With Micro-Macro correction.}
		\label{fig_global_error_VDP_FEuler}
    \end{figure}

    We also investigate the global error between the exact flow and the numerical flow obtained from the learned vector fields. The errors are plotted as a function of the step size, with the curves showing perfect agreement with the estimates from the previous theorems (for the Forward Euler and Midpoint methods). Figures \ref{fig_global_error_pendulum_FEuler}, \ref{fig_global_error_pendulum_MidPoint}, and \ref{fig_global_error_VDP_FEuler} demonstrate that numerical integration using the micro-macro correction is more accurate than the \textit{slow-fast} decomposition-based method, confirming the results obtained in subsection \ref{subsection:Loss_decay_integration_ODE}.

    \subsection{Uniform accuracy test}

    Since the micro-macro method is presented as a \textit{uniformly accurate} (UA) method \cite{AVERAGING_chartier2020new,AVERAGING_chartier2022derivative},  we can check whether the micro-macro correction with machine learning retains uniform accuracy by plotting the global errors against the parameter $\eps$ for various step sizes $h$.\\

\begin{figure}[H]
		\centering
		\begin{minipage}{0.45\linewidth}
			\includegraphics[width=\linewidth]{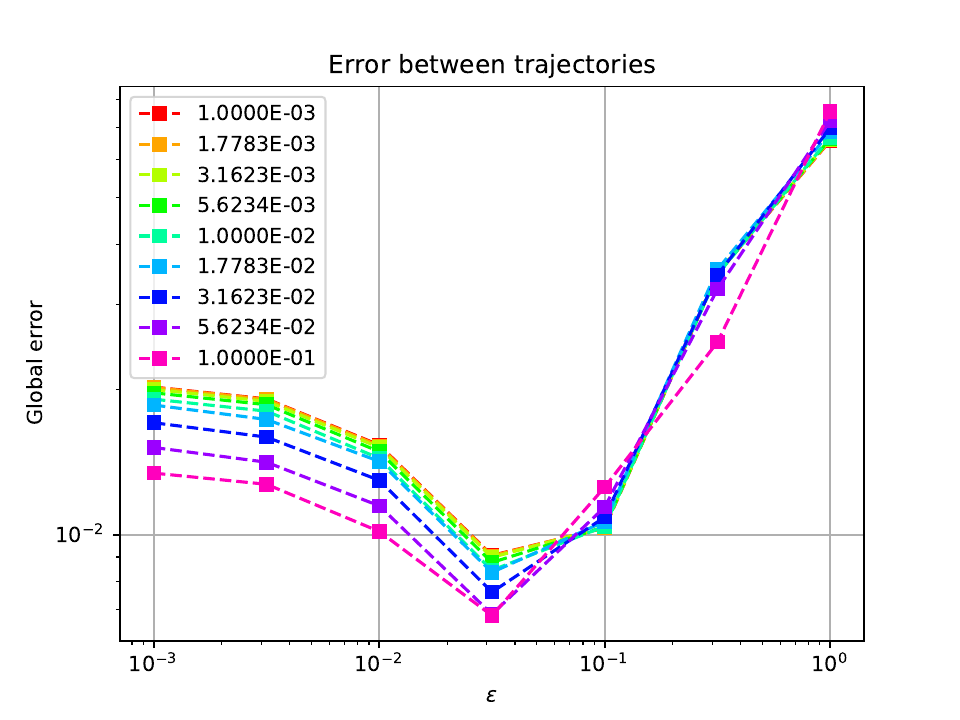}
		\end{minipage}
		\begin{minipage}{0.45\linewidth}
			\includegraphics[width=\linewidth]{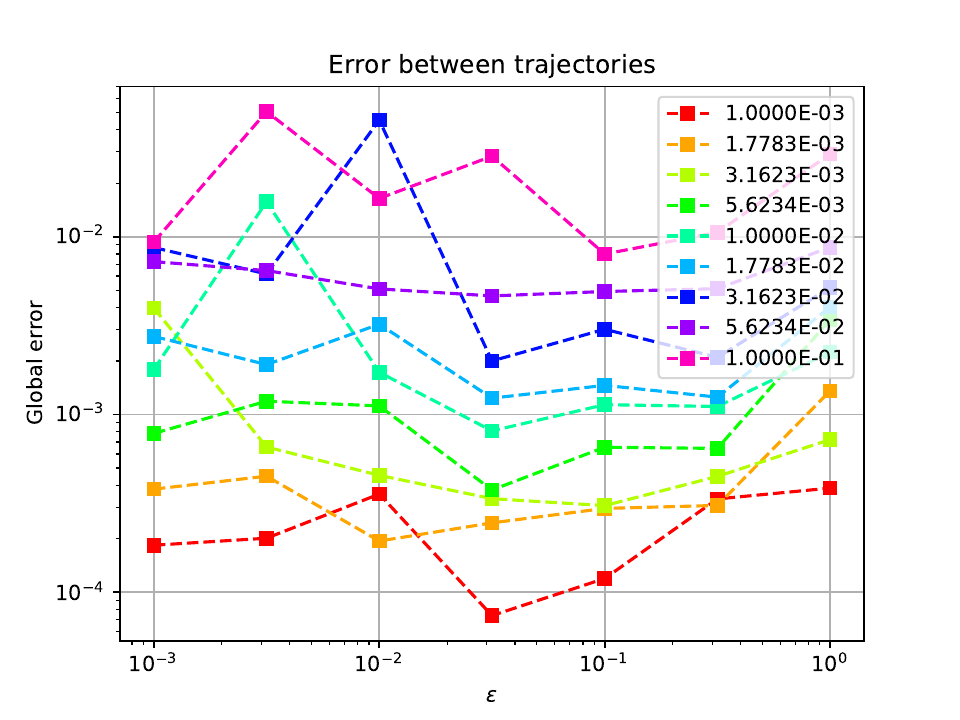}
		\end{minipage}
		\caption{Uniform accuracy test (each color corresponds to a step size $h$) of Inverted Pendulum system with Forward Euler. Left: \textit{Slow-fast} decomposition-based method. Right: With Micro-Macro correction.}
		\label{fig_global_error_UA_pendulum_FEuler}
    \end{figure}

\begin{figure}[H]
		\centering
		\begin{minipage}{0.45\linewidth}
			\includegraphics[width=\linewidth]{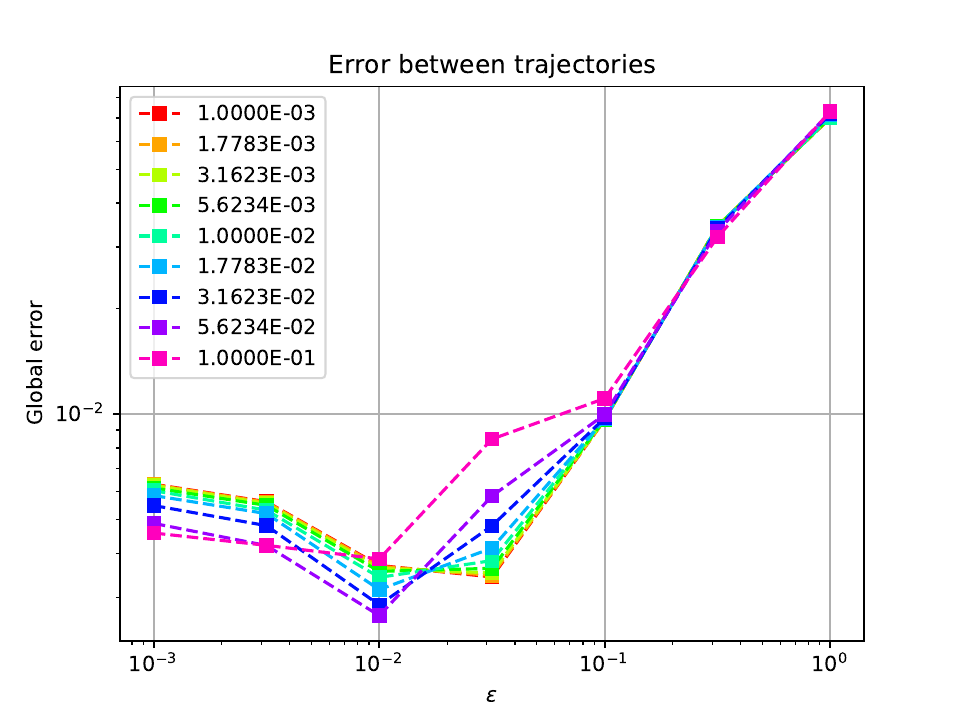}
		\end{minipage}
		\begin{minipage}{0.45\linewidth}
			\includegraphics[width=\linewidth]{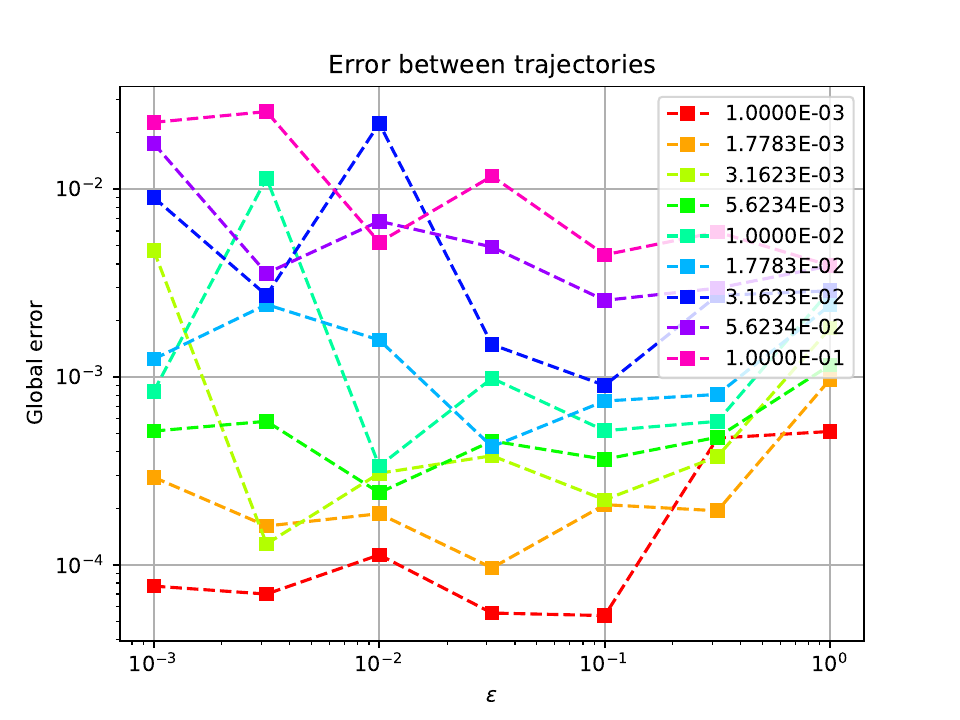}
		\end{minipage}
		\caption{Uniform accuracy test (each color corresponds to a step size $h$) of Inverted Pendulum system with midpoint. Left: \textit{Slow-fast} decomposition-based method. Right: With Micro-Macro correction.}
		\label{fig_global_error_UA_pendulum_Midpoint}
    \end{figure}

\begin{figure}[H]
		\centering
		\begin{minipage}{0.45\linewidth}
			\includegraphics[width=\linewidth]{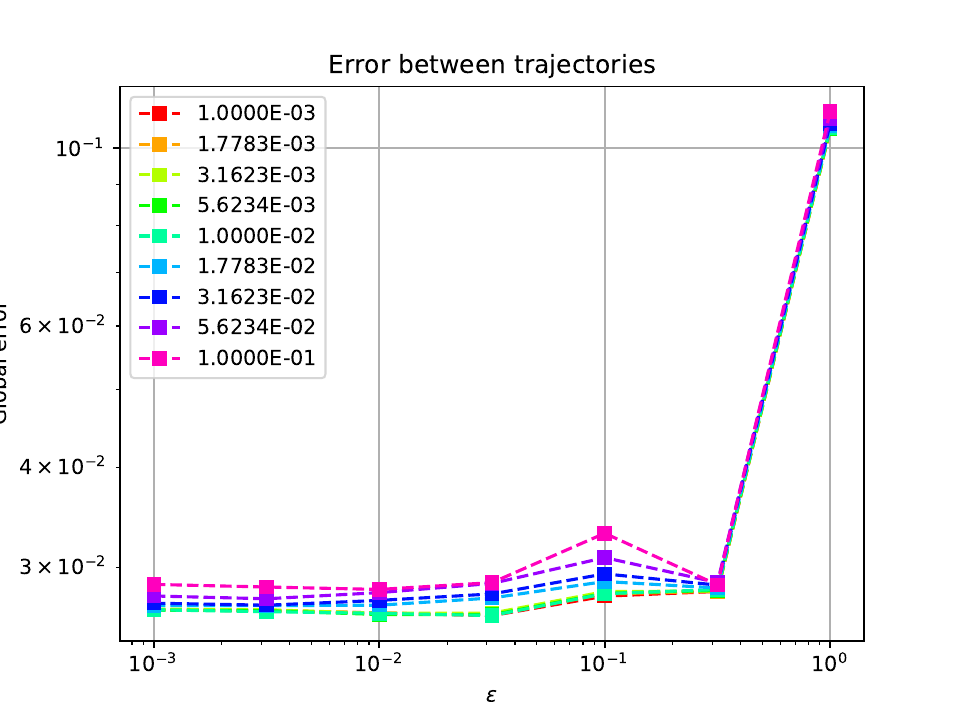}
		\end{minipage}
		\begin{minipage}{0.45\linewidth}
			\includegraphics[width=\linewidth]{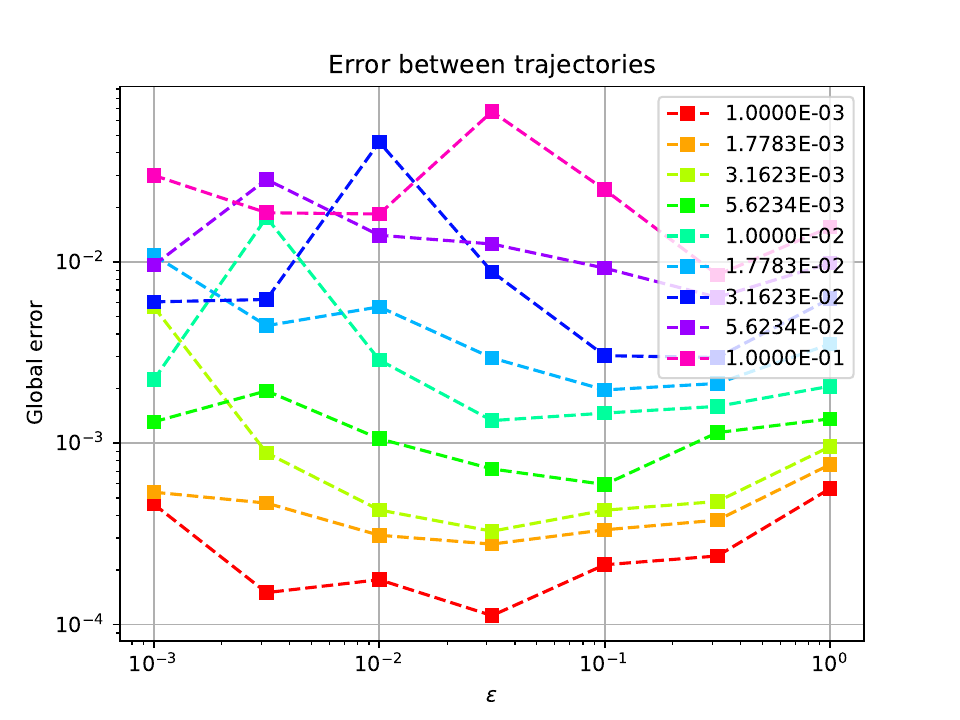}
		\end{minipage}
		\caption{Uniform accuracy test (each color corresponds to a step size $h$) of Van der Pol oscillator with Forward Euler. Left: \textit{Slow-fast} decomposition-based method. Right: With Micro-Macro correction.}
		\label{fig_global_error_UA_VDP_FEuler}
    \end{figure}

    Figures \ref{fig_global_error_UA_pendulum_FEuler}, \ref{fig_global_error_UA_pendulum_Midpoint}, and \ref{fig_global_error_UA_VDP_FEuler} show that uniform accuracy is nearly verified with the micro-macro correction. However, this property is not observed with the \textit{slow-fast} decomposition-based method, particularly for larger values of $\eps$. This behavior can be attributed to the exponential remainder in the formula \eqref{estimate:exponential_error}.

    \subsection{Evaluation of alternative method}

\begin{figure}[H]
        \centering
        \includegraphics[scale = 0.6]{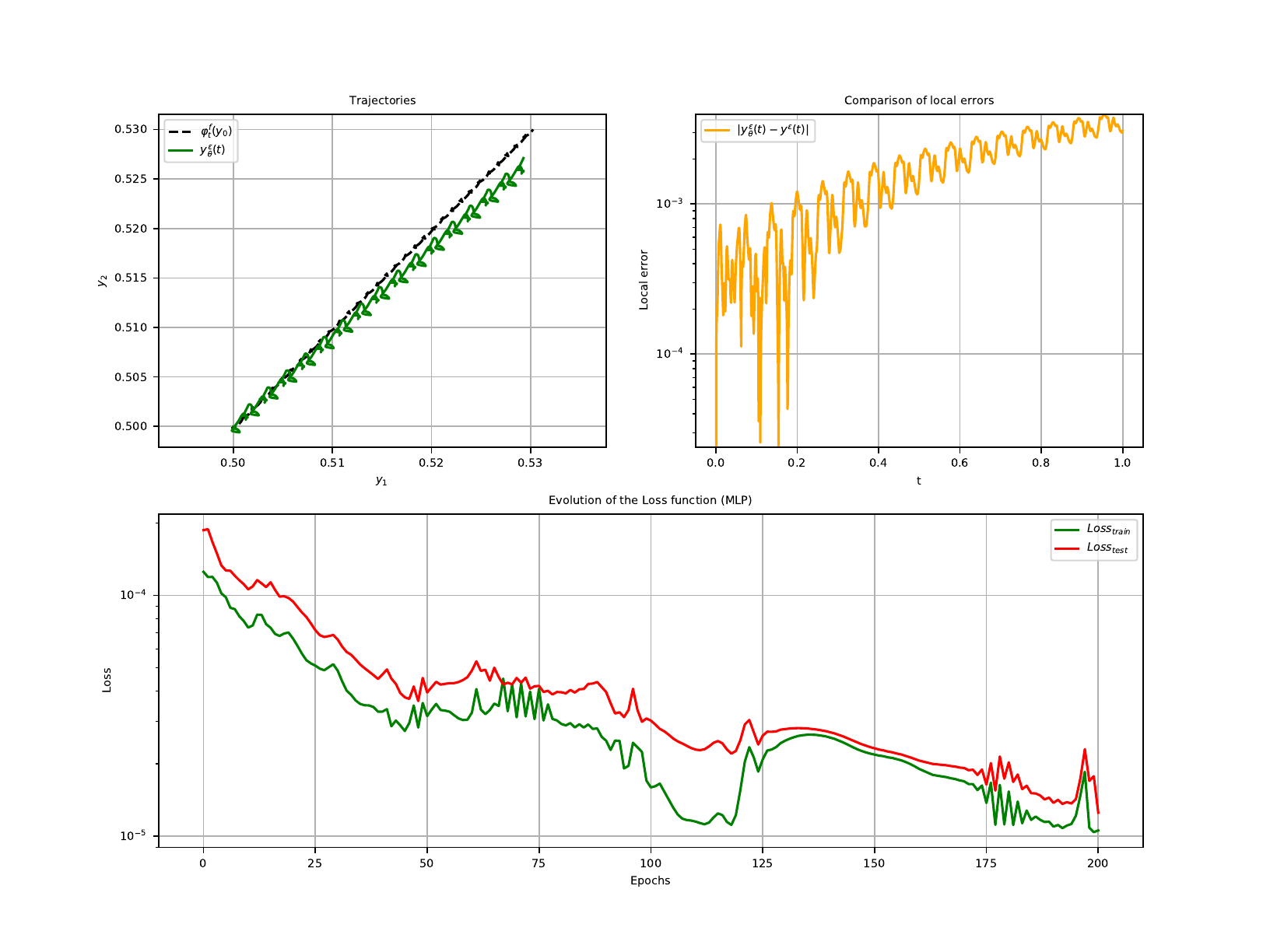}
        \caption{Comparison between $Loss$ decays (green: $Loss_{Train}$, red: $Loss_{Test}$), trajectories (dashed dark: exact flow, green: numerical flow with learned vector fields and local error (yellow) for the Van der Pol oscillator with Forward Euler method in the case $\eps = 0.01$.}
        \label{fig_Comparison_VDP_FEuler_eps=0.01}
    \end{figure}

\begin{figure}[H]
        \centering
        \includegraphics[scale = 0.8]{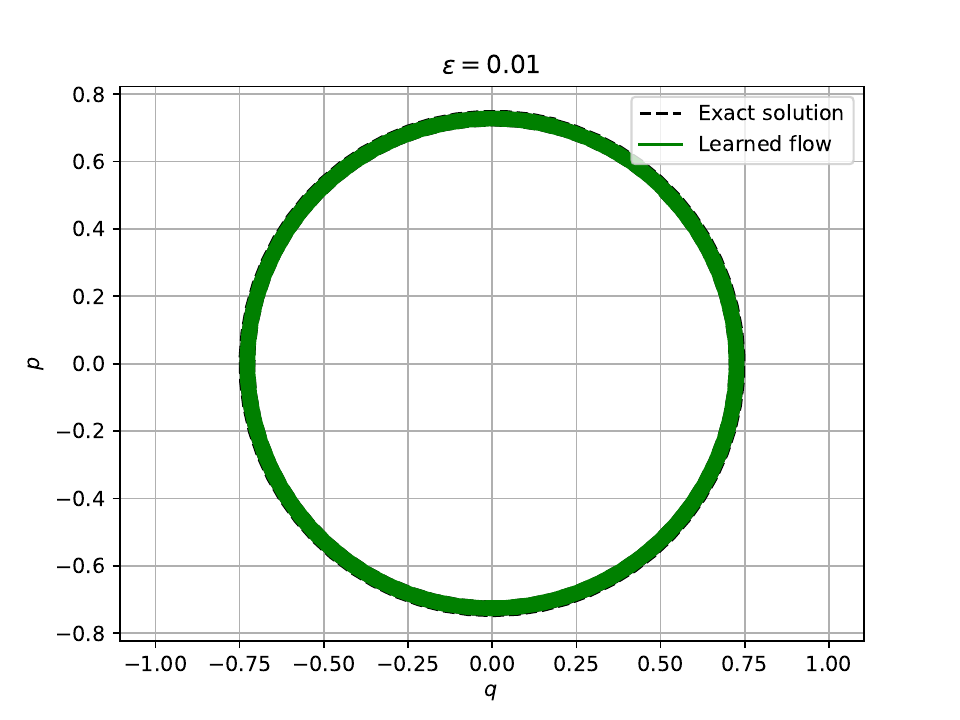}
        \caption{Comparison between trajectories (dashed dark: exact flow, green: numerical flow with learned vector field) for the Van der Pol oscillator with Forward Euler method in the case $\eps = 0.01$ after the inverse variable change \eqref{Variable_change_VDP}.}
        \label{fig_Comparison_VDP_FEuler_Variable_change_eps=0.01}
    \end{figure}

\begin{figure}[H]
        \centering
        \includegraphics[scale = 0.6]{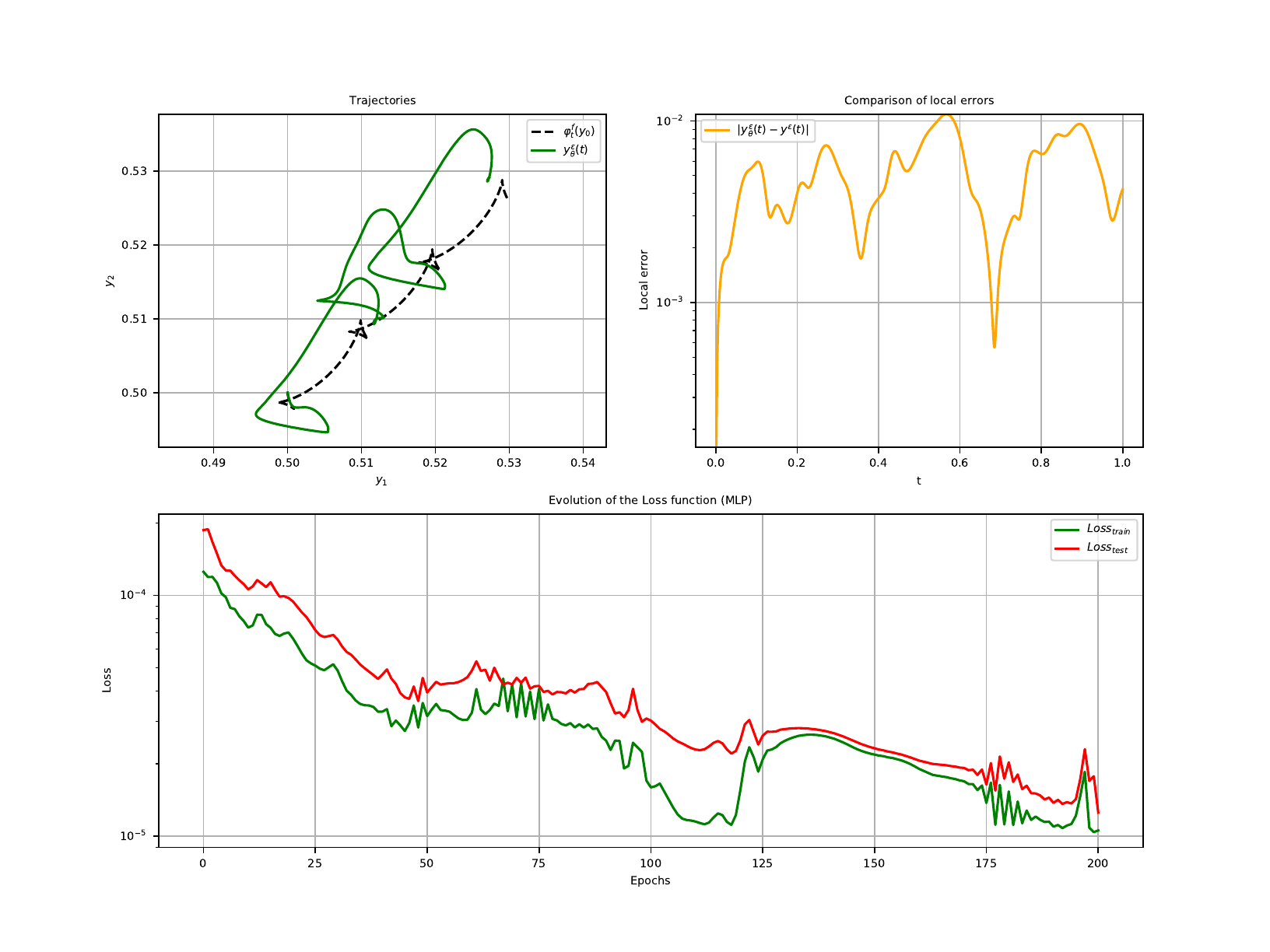}
        \caption{Comparison between $Loss$ decays (green: $Loss_{Train}$, red: $Loss_{Test}$), trajectories (dashed dark: exact flow, green: numerical flow with learned vector fields and local error (yellow) for the Van der Pol oscillator with Forward Euler method in the case $\eps = 0.1$.}
        \label{fig_Comparison_VDP_FEuler_eps=0.1}
    \end{figure}

\begin{figure}[H]
        \centering
        \includegraphics[scale = 0.8]{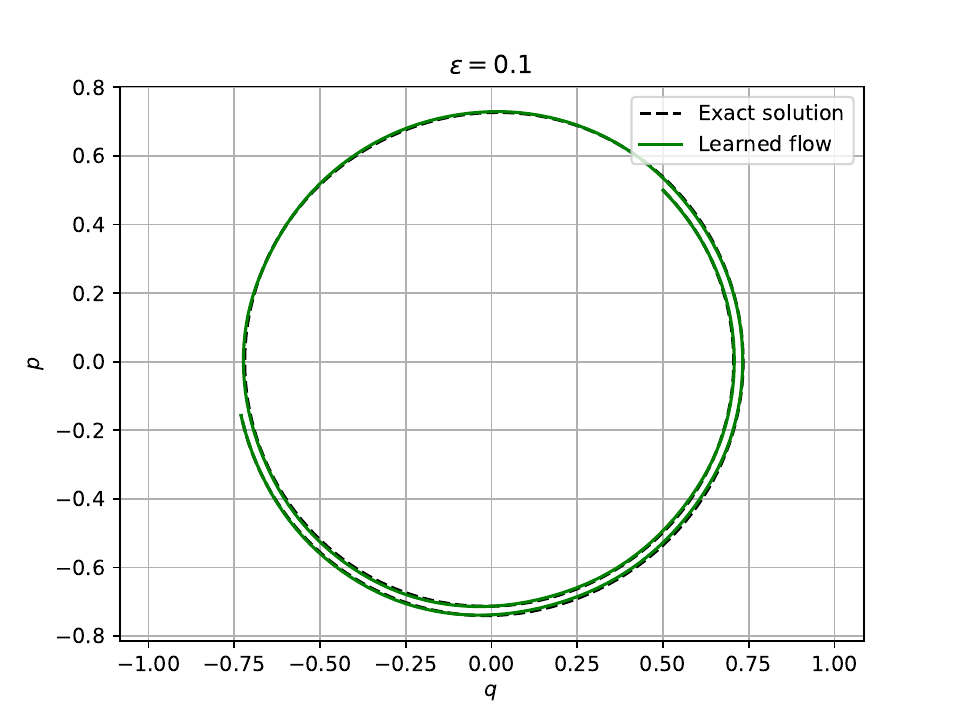}
        \caption{Comparison between trajectories (dashed dark: exact flow, green: numerical flow with learned vector field) for the Van der Pol oscillator with Forward Euler method in the case $\eps = 0.1$ after the inverse variable change \eqref{Variable_change_VDP}.}
        \label{fig_Comparison_VDP_FEuler_Variable_change_eps=0.1}
    \end{figure}

\begin{figure}[H]
        \centering
        \includegraphics[scale = 0.6]{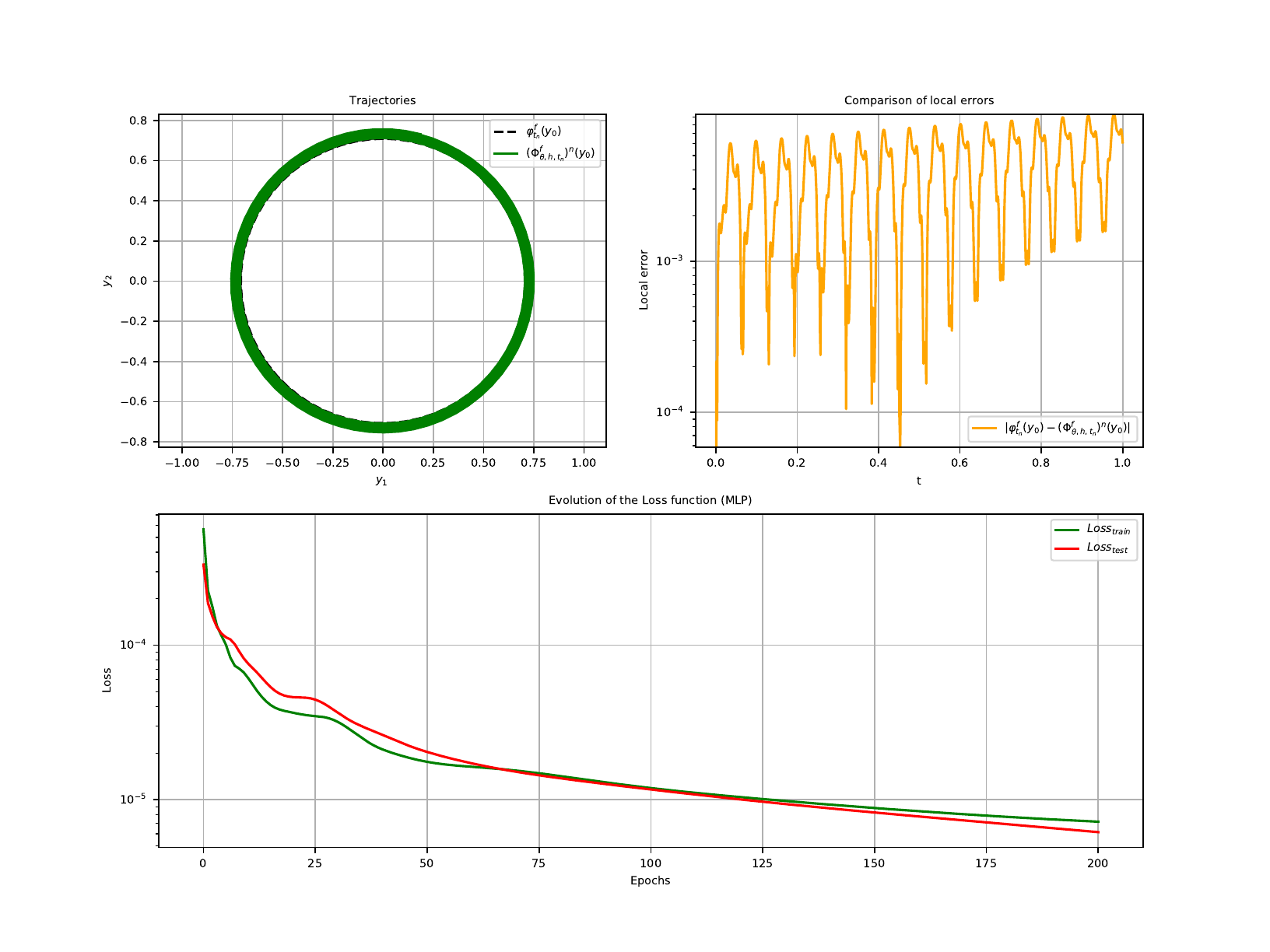}
        \caption{Comparison between $Loss$ decays (green: $Loss_{Train}$, red: $Loss_{Test}$), trajectories (dashed dark: exact flow, green: numerical flow with learned vector fields and local error (yellow) for the Van der Pol oscillator with alternative method in the case $\eps = 0.01$.}
        \label{fig_Comparison_VDP_FEuler_eps=0.01_Autonomous}
    \end{figure}

\begin{figure}[H]
        \centering
        \includegraphics[scale = 0.6]{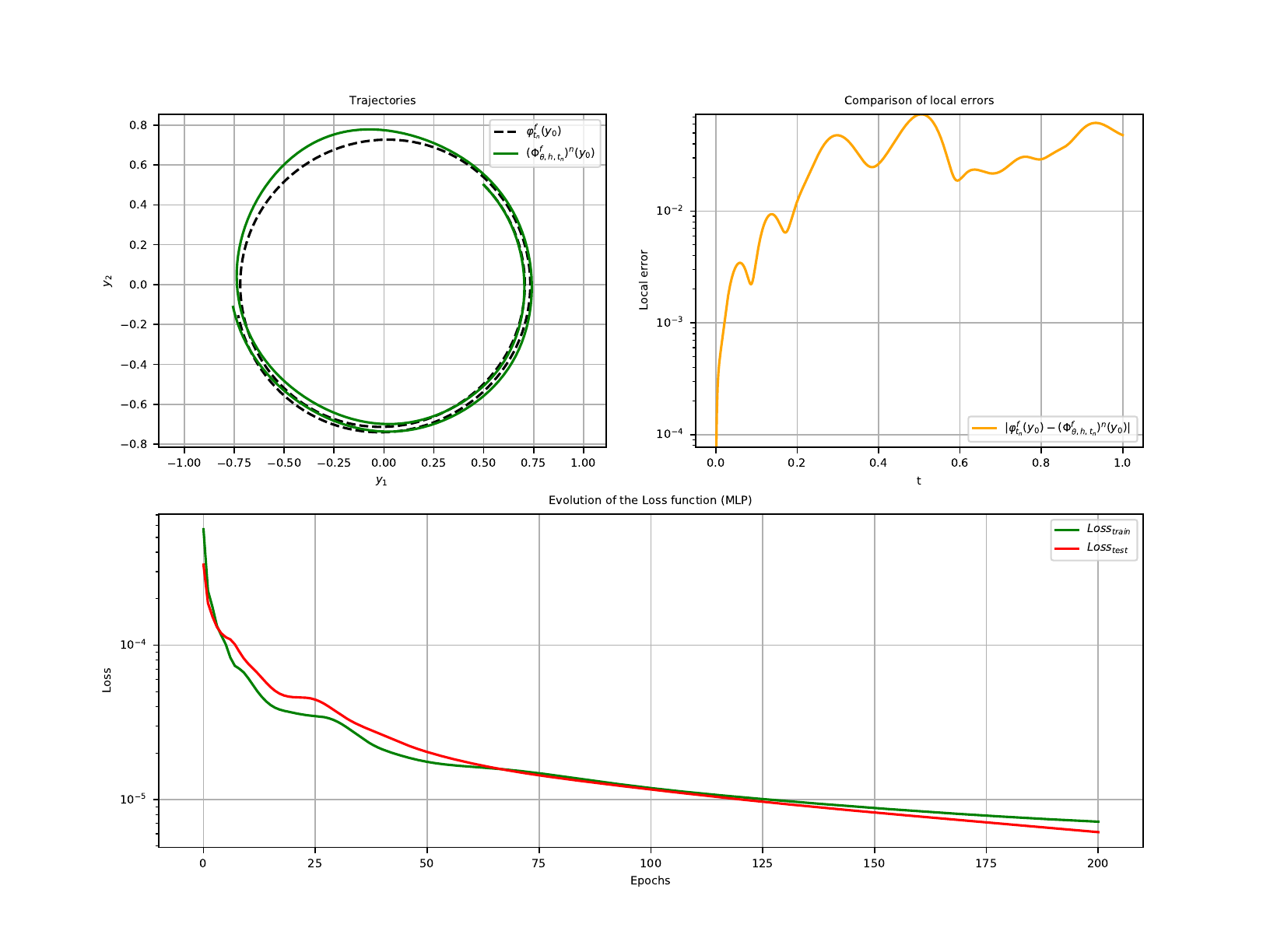}
        \caption{Comparison between $Loss$ decays (green: $Loss_{Train}$, red: $Loss_{Test}$), trajectories (dashed dark: exact flow, green: numerical flow with learned vector fields and local error (yellow) for the Van der Pol oscillator with alternative method in the case $\eps = 0.1$.}
        \label{fig_Comparison_VDP_FEuler_eps=0.1_Autonomous}
    \end{figure}

\begin{figure}[H]
		\centering
		\begin{minipage}{0.45\linewidth}
			\includegraphics[width=\linewidth]{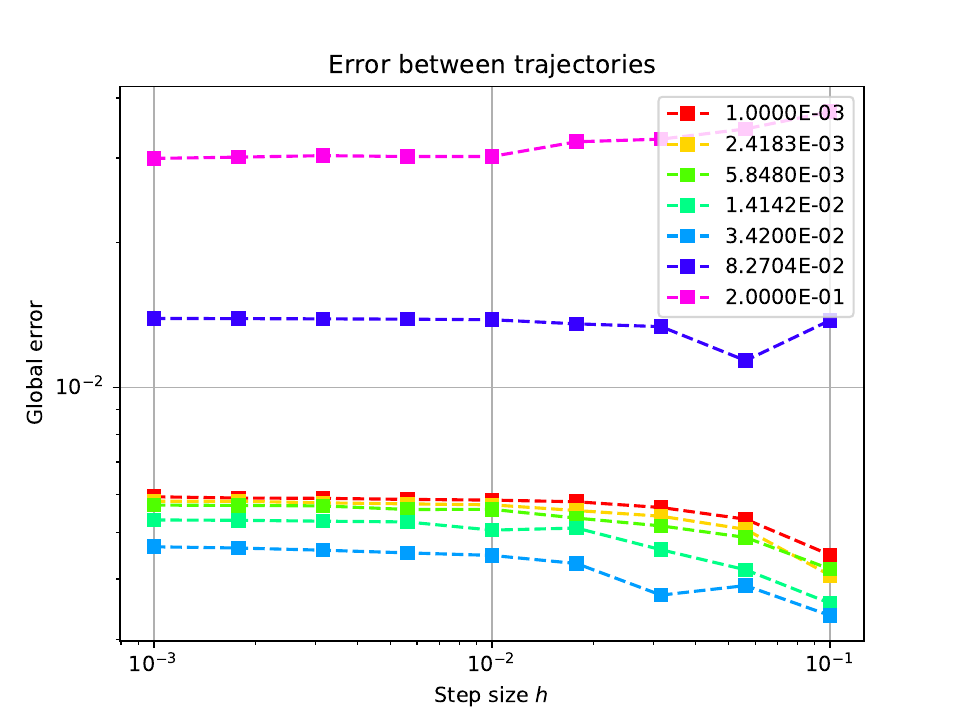}
		\end{minipage}
		\begin{minipage}{0.45\linewidth}
			\includegraphics[width=\linewidth]{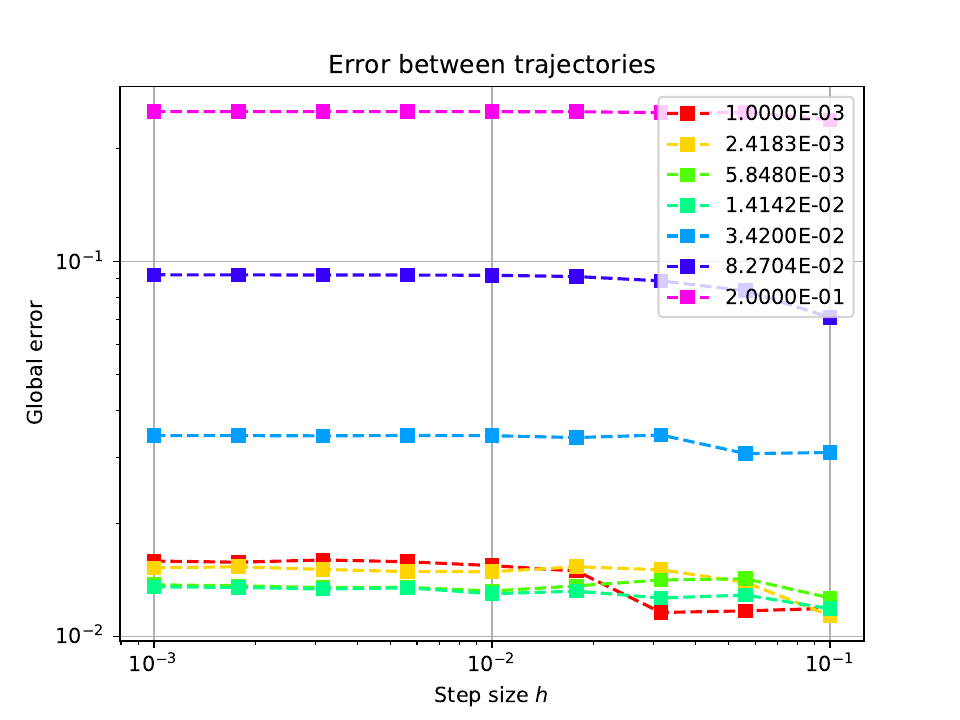}
		\end{minipage}
		\caption{Integration errors (each color corresponds to a high oscillation parameter $\eps$) of Van der Pol oscillator. Left: \textit{Slow-fast} decomposition-based method. Right: Alternative method for autonomous case.}
		\label{fig_Comparison_global_error_VDP}
    \end{figure}

    In this subsection, we compare the classical method using \textit{slow-fast} decomposition with an alternative method in the autonomous case, using the Van der Pol oscillator as an example. Despite the absence of an auto-encoder, the classical method appears to yield more accurate solutions (as shown in Figures \ref{fig_Comparison_VDP_FEuler_eps=0.01}, \ref{fig_Comparison_VDP_FEuler_Variable_change_eps=0.01}, \ref{fig_Comparison_VDP_FEuler_eps=0.1}, and \ref{fig_Comparison_VDP_FEuler_Variable_change_eps=0.1}) than the alternative method (Figures \ref{fig_Comparison_VDP_FEuler_eps=0.01_Autonomous} and \ref{fig_Comparison_VDP_FEuler_eps=0.1_Autonomous}). Furthermore, Figure \ref{fig_Comparison_global_error_VDP} confirms the difference in accuracy.

    
    \section{Conclusions}

    The numerical experiments presented in this paper demonstrate the feasibility of learning both the highly oscillatory generator and the averaged field using neural networks. Additionally, numerical integration without pre-computation provides accurate approximations of the exact solutions to highly oscillatory differential equations. While the method based on \textit{slow-fast} decomposition is simpler and computationally less demanding than the micro-macro correction, it still achieves nearly uniform accuracy and consistency—though perfect uniform accuracy and consistency are limited by learning errors. Given the challenges of applying machine learning to high-dimensional problems, our methods are particularly well-suited for low-dimensional systems.
    
    \section*{Acknoledgements}
    
    The author would like to thank Philippe Chartier, Mohamed Lemou and Florian Méhats for their guidance, valuable advice and careful reading of this work.

    \bibliographystyle{plainnat}
    \bibliography{./biblio/biblio}

    \appendix

    \section{Proof of theorem \ref{theorem:highly_oscillatory_general_slow-fast}}

    Let $t \mapsto \psi_t^\eps$ represent the exact solution of the autonomous differential equation associated with the vector field $F^{\eps,[n_\eps]}$. Additionally, let $\psi_{\theta,n}^\eps$ denote the corresponding numerical flow:

    \begin{equation}
        \left\{\begin{array}{c c l}
            \psi_{\theta,0}^\eps & = & \yeps(0) \\
            \psi_{\theta,n+1}^\eps & = & \Phi_h^{F_{\theta}(\cdot,h,\eps)}(\psi_{\theta,n}^\eps) 
        \end{array}\right.
    \end{equation}
    
    First, we examine the estimates over this autonomous ODE. Then a general error estimate is provided.\\

    \begin{enumerate}[label=\textbf{(\roman*).}]
        \item \textbf{Consistency error (autonomous ODE):} Consistency error is given by

        \begin{eqnarray}
            \eps_{\theta,n} & := & \psi_{t_{n+1}}^\eps - \Phi_h^{F_{\theta}(\cdot,h,\eps)}(\psi_{t_n}^\eps) \\
            & = & \pphi_h^{F^{\eps,[n_\eps]}}(\psi_{t_n}^\eps) - \Phi_h^{\widetilde{F_h^{\eps,[n_\eps],[q]}}}(\psi_{t_n}^\eps) + \Phi_h^{\widetilde{F_h^{\eps,[n_\eps],[q]}}}(\psi_{t_n}^\eps) - \Phi_h^{F_{\theta}(\cdot,h,\eps)}(\psi_{t_n}^\eps) \nonumber
        \end{eqnarray}

        According to the first hypothesis of theorem $\ref{theorem:highly_oscillatory_general_slow-fast}$, we obtain

        \begin{equation}
            |\eps_{\theta,n}| \leqslant \overline{C}h^{p+q} + C\delta_F h
        \end{equation}

        \item \textbf{Local truncation error (autonomous ODE):} Local truncation error is given by

        \begin{equation}
            e_{\theta,n} := \psi_{\theta,n}^\eps - \psi_{t_n}^\eps.
        \end{equation}

        Thus we have

        \begin{equation}
            e_{\theta,n+1} = \Phi_h^{F_{\theta}(\cdot,h,\eps)}(\psi_{\theta,n}^\eps) - \Phi_h^{F_{\theta}(\cdot,h,\eps)}(\psi_{t_n}^\eps) - \eps_{\theta,n}.
        \end{equation}

        Thus there exists a constant $\lambda$ (the Lipschitz constant of $F_{\theta}$ with respect to its space variable) such that

        \begin{equation}
            |e_{\theta,n+1}| \leqslant (1+\lambda h)|e_{\theta,n}| + \overline{C}h^{p+q} + C\delta_F h.
        \end{equation}

        Using the dicrete Grönwall lemma we obtain

        \begin{equation}
            \underset{0 \leqslant n \leqslant N}{Max}|e_{\theta,n}| \leqslant \frac{e^{\lambda T}-1}{\lambda}\left[\overline{C}h^{p+q-1} + C\delta_F\right]
        \end{equation}

        \item \textbf{General error estimate:} Let us consider and denote $\eta_{\theta,n}$ as the following error

        \begin{eqnarray}
            \eta_{\theta,n} & := & \yeps(t_n) - y_{\theta,n}^\eps \\
            & = & \yeps(t_n) - \phi_{\theta,+}\left( \frac{t_n}{\eps} , \psi_{\theta,n}^\eps , \eps \right). \nonumber
        \end{eqnarray}

        By employing the following decomposition:

        \begin{eqnarray}
            \eta_{\theta,n} & = & \yeps(t_n) - \phi_{\frac{t_n}{\eps}}^{\eps,[n_\eps]}(\psi_{t_n}^\eps) + \phi_{\frac{t_n}{\eps}}^{\eps,[n_\eps]}(\psi_{t_n}^\eps) - \phi_{\theta,+}\left(\frac{t_n}{\eps},\psi_{t_n}^\eps , \eps\right) \\
            & + & \phi_{\theta,+}\left(\frac{t_n}{\eps},\psi_{t_n}^\eps , \eps\right) - \phi_{\theta,+}\left(\frac{t_n}{\eps},\psi_{\theta,n}^\eps , \eps\right), \nonumber
        \end{eqnarray}
        
        we bound the first, second and third terms from above using $(\ref{estimate:exponential_error})$, the learning error $\delta_{\phi,+}$ and the estimate for $e_{\theta,n}$. This will provide the desired estimate.
        
    \end{enumerate}

    \section{Proof of theorem \ref{theorem:highly_oscillatory_general_Micro-Macro}}

    \begin{enumerate}[label  = \textbf{\arabic*.}]
        \item \textbf{Error estimate for $v_{\theta,n}$:}

        For the estimate of $v$, we apply a classical error estimation technique commonly used for autonomous differential equations, similar to the estimation of the autonomous part in the proof of \ref{theorem:highly_oscillatory_general_slow-fast}.

        \begin{enumerate}[label = \textbf{(\roman*).}]
            \item \textbf{Consistency Error:} Consistency error is given by

            \begin{eqnarray}
                \eps_{\theta,v,n} & := & v(t_{n+1}) - \Phi_h^{F_{\theta}(\cdot,h,\eps)}(v(t_n))  \\
                & = & \pphi_h^{F^{[p]}}(v(t_n)) - \Phi_h^{\widetilde{F^{[p],[q]}_h}} (v(t_n)) + \Phi_h^{\widetilde{F^{[p],[q]}_h}} (v(t_n)) - \Phi_h^{F_{\theta}(\cdot,h,\eps)}(v(t_n)) \nonumber
            \end{eqnarray}
            
            according to the first hypothesis of Theorem \ref{theorem:highly_oscillatory_general_Micro-Macro} we obtain

            \begin{eqnarray}
                |\eps_{\theta,v,n}| & \leqslant & \overline{C}h^{p+q} + C\delta_Fh.
            \end{eqnarray}
    
            \item \textbf{Local Truncation Error:} Local truncation error is given by

            \begin{eqnarray}
                e_{\theta,v,n} & := & v(t_n) - v_{\theta,n}.
            \end{eqnarray}

            Therefore, we have

            \begin{eqnarray}
                e_{\theta,v,n+1} & := & v(t_{n+1}) - v_{\theta,n+1} \nonumber \\
                & = & \Phi_h^{F_{\theta}(\cdot,h,\eps)}(v(t_n)) - \Phi_h^{F_{\theta}(\cdot,h,\eps)}(v_{\theta,n}) + \eps_{\theta,v,n}.
            \end{eqnarray}

            Thus, there exists a constant $\lambda$ (the Lipschitz constant of $F_\theta$ w.r.t. space variable) such that:

            \begin{eqnarray}
                |e_{\theta,v,n+1}| & \leqslant & (1+\lambda h)|e_{v,\theta,n}| + \overline{C}h^{p+q} + C\delta_Fh.
            \end{eqnarray}

            By applying the discrete Grönwall lemma, we obtain

            \begin{eqnarray}
                \underset{0 \leqslant n \leqslant N}{Max}|e_{\theta,v,n}| & \leqslant & \frac{e^{\lambda T}-1}{\lambda}\left[\overline{C}h^{p+q-1} + C\delta_F\right]
            \end{eqnarray}
        \end{enumerate}

        \item \textbf{Error estimate for $w_{\theta,n}$:}

        For the estimate of $w$, we treat $v$ as a source term and perform an error estimation similar to the approach used for autonomous differential equations.

        \begin{enumerate}[label = \textbf{(\roman*)}.]
            \item \textbf{Consistency Error:} Consistency error is given by

            \begin{eqnarray}
                \eps_{\theta,w,n} & := & w(t_{n+1}) - \Phi_h^{g_\theta\left(\frac{t_n}{\eps},\cdot,v(t_n)\right)}w(t_n) \\
                & = & \underbrace{w(t_{n+1}) - \Phi_{t_n,h}^{g\left(\frac{t_n}{\eps},\cdot,v(t_n)\right)}w(t_n)}_{\text{Consistency error for classical scheme}}          + \Phi_{t_n,h}^{g\left(\frac{t_n}{\eps},\cdot,v(t_n)\right)}w(t_n) - \Phi_{t_n,h}^{g_\theta\left(\frac{t_n}{\eps},\cdot,v(t_n)\right)}w(t_n). \nonumber
            \end{eqnarray}

            Since the first term can be bounded from above by $Mh^{p+1}$ (where $M$ is a constant independent of $h$ and $\eps$), and the second term can be bounded from above using the first hypothesis of theorem \ref{theorem:highly_oscillatory_general_Micro-Macro}, we obtain

            \begin{eqnarray}
                |e_{\theta,w,n}| & \leqslant & M'h^{p+1} + C\delta_gh
            \end{eqnarray}

            \item \textbf{Local Truncation Error:} Local truncation error is given by

            \begin{eqnarray}
                e_{\theta,w,n} & = & w(t_n) - w_{\theta,n}.
            \end{eqnarray}

            Therefore, we have

            \begin{eqnarray}
                e_{\theta,w,n+1} & = & w(t_{n+1}) - w_{\theta,n+1} \nonumber \\
                & = & \Phi_{t_n,h}^{g_\theta\left(\frac{t_n}{\eps},\cdot,v(t_n)\right)}(w(t_n)) - \Phi_{t_n,h}^{g_\theta\left(\frac{t_n}{\eps},\cdot,v_{\theta,n}\right)}(w_{\theta,n}) + \eps_{\theta,w,n} \nonumber \\
                & = & \Phi_{t_n,h}^{g_\theta\left(\frac{t_n}{\eps},\cdot,v(t_n)\right)}(w(t_n)) - \Phi_{t_n,h}^{g_\theta\left(\frac{t_n}{\eps},\cdot,v_{\theta,n}\right)}(w(t_n)) \\
                & + & \Phi_{t_n,h}^{g_\theta\left(\frac{t_n}{\eps},\cdot,v_{\theta,n}\right)}(w(t_n)) - \Phi_{t_n,h}^{g_\theta\left(\frac{t_n}{\eps},\cdot,v_{\theta,n}\right)}(w_{\theta,n}) \nonumber \\
                & + & \eps_{\theta,w,n} \nonumber
            \end{eqnarray}

            The first term can be bounded from above using the first hypothesis of theorem \ref{theorem:highly_oscillatory_general_Micro-Macro}, while the second term can be estimated by introducing $\mu$, the Lipschitz constant of $g_\theta$, and applying the second hypothesis of theorem \ref{theorem:highly_oscillatory_general_Micro-Macro}. Thus we obtain

            \begin{eqnarray}
                |e_{\theta,w,n+1}| & \leqslant & Ch\lnorm g_\theta\left(\frac{t_n}{\eps},\cdot,v(t_n)\right) - g_\theta\left(\frac{t_n}{\eps},\cdot,v_{\theta,n}\right) \rnorm_{\Linfty(\Omega)} \nonumber \\
                & + & (1+\mu h)|e_{\theta,w,n}| + Mh^{p+1} + Q\delta_gh  \\
                & \leqslant & C\beta h|e_{\theta,v,n}| + (1+\mu h)|e_{\theta,w,n}| + M'h^{p+1} + C\delta_gh, \nonumber
            \end{eqnarray}
            
            where $\beta$ is the Lipschitz constant of $g_\theta$ with respect to $v$. Since $|e_{\theta,w,n}|$ can be bounded from above by using the error estimate for $v$, we obtain
            
            \begin{eqnarray}
                |e_{\theta,w,n+1}| & \leqslant & (1+\mu h)|e_{\theta,w,n}| + M'h^{p+1} \\
                & + & h\left[ \frac{e^{\lambda T}-1}{\lambda}\beta(\overline{C}h^{p+q-1}+C\delta_F) + C\delta_g \right]. \nonumber
            \end{eqnarray}

            Using the discrete Grönwall lemma, we obtain

            \begin{eqnarray}
                \underset{0 \leqslant n \leqslant N}{Max}|e_{\theta,w,n}| & \leqslant & \frac{e^{\mu T}-1}{\mu}\left[ M'h^p + \frac{e^{\lambda T}-1}{\lambda}\beta(\overline{C}h^{p+q-1}+C\delta_F) + C\delta_g \right]
            \end{eqnarray}
        \end{enumerate}

        \item \textbf{Error estimate for $\yeps_{\theta,n}$:} Finally, we perform an error estimate over $y$ using the formula \eqref{formula:Micro-Macro_ML} to describe our numerical method. Since we have

        \begin{equation}
            \yeps_{\theta,n} = \phi_{\theta,+}\left( \frac{t_n}{\eps} , v_{\theta,n} , \eps \right) + w_{\theta,n},
        \end{equation}

        The local truncation errors can be expressed as follows:

        \begin{eqnarray}
            e_{\theta,n} & := & \yeps(t_n) - \yeps_{\theta,n} \nonumber \\
            & = & \phi^{[p]}_{\frac{t_n}{\eps}}(v(t_n)) + w(t_n) - \phi_{\theta,+}\left(\frac{t_n}{\eps},v_{\theta,n},\eps\right) - w_{\theta,n} \\
            & = & \phi^{[p]}_{\frac{t_n}{\eps}}(v(t_n)) - \phi^{[p]}_{\frac{t_n}{\eps}}(v_{\theta,n}) + \phi^{[p]}_{\frac{t_n}{\eps}}(v_{\theta,n}) - \phi_{\theta,+}\left(\frac{t_n}{\eps},v_{\theta,n},\eps\right) + w(t_n) - w_{\theta,n}. \nonumber
        \end{eqnarray}

        The first difference term can be bounded from above using $\alpha_\phi$, the Lipschitz constant of $\phi^{[p]}$, while the second difference term can be estimated using the learning error for $\phi_{\theta,+}$. Therefore we obtain

        \begin{eqnarray}
            |e_{\theta,n}| & \leqslant & \alpha_{\phi}|e_{\theta,v,n}| + \delta_{\phi,+}\eps + |e_{\theta,w,n}|
        \end{eqnarray}
        
        and the desired estimate follows.
        
    \end{enumerate}

    \section{Proof of theorem \ref{theorem:autonomous_highly_oscillatory}}

    A standard proof used for autonomous ODEs (consistency and local truncation error) gives the following result. If we denote $t \longmapsto \psi_t^\eps$, the solution associated with $g^{\eps}$ and define $\displaystyle \psi_{\theta,n}^\eps = \left(\pphi_\theta(\cdot,h,\eps)\right)^n(\yeps(0))$, then:

    \begin{equation}
        \left| \psi_{t_n}^\eps - \psi_{\theta,n}^\eps \right| \leqslant \frac{e^{\lambda T}-1}{\lambda}\delta\pphi,
    \end{equation}

    \noindent where $\lambda$ is the Lipschitz constant of $R_{\theta,\pphi}$ with respect to $y$.\\

    Then, using the following decomposition:

    \begin{eqnarray}
        \yeps_{\theta,n} - \yeps(t_n) & = & \phi_\theta\left( \frac{t_n}{\eps} , \psi_{\theta,n}^\eps , \eps \right) - \Phi_{\frac{t_n}{\eps}}^\eps\left(\psi_{\theta,n}^\eps\right) \nonumber \\
        & + & \Phi_{\frac{t_n}{\eps}}^\eps\left(\psi_{\theta,n}^\eps\right) - \Phi_{\frac{t_n}{\eps}}^\eps\left( \pphi_{t_n}^{g^\eps}(\yeps(0)) \right) \\
        & + & \phi_{\frac{t_n}{\eps}}^\eps\left( \pphi_{t_n}^{g^\eps}(\yeps(0)) \right) - \yeps(t_n), \nonumber
    \end{eqnarray}

    Since the first term is bounded from above by $\delta_\phi$, the second term is bounded from above by $L\left| \psi_{t_n}^\eps - \psi_{\theta,n}^\eps \right|$, where $L$ is the Lipschitz constant of $\Phi^\eps_{\frac{t_n}{\eps}}$ w.r.t. space variable, and the third term is estimated by the exponential remainder, we get the desired estimate.

    \section{Implementation of implicit methods}

    In the formula $(\ref{Loss_train})$, the numerical flow takes the input $\phi_{\theta,-}\left(\frac{t_0}{\eps},y_0,\eps\right)$. For implicit methods, however, the numerical flow is a function of both the input and the output. $\phi_{\theta,-}\left(\frac{t_0+h}{\eps},y_1,\eps\right)$ is considered as the output, and if we consider, for example, the midpoint rule, we have

    \begin{equation}
        \Phi_h^{F_{\theta}(\cdot,h,\eps)}\left( \phi_{\theta,-}\left(\frac{t_0}{\eps},y_0,\eps\right) \right) = \phi_{\theta,-}\left(\frac{t_0}{\eps},y_0,\eps\right) + hF_{\theta}\left(\frac{\phi_{\theta,-}\left(\frac{t_0}{\eps},y_0,\eps\right) + \phi_{\theta,-}\left(\frac{t_0+h}{\eps},y_1,\eps\right)}{2} , h , \eps\right)
    \end{equation}

    \section{Computation of learning errors}

    \subsection{Space and time discretizations}

    To compute the learning error with respect to $\eps$, we discretise the space and time domains. In the formulas \eqref{formula:learning error_epsilon} we make the following approximations:

    \begin{equation}
        \underset{y\in\Omega}{Max}\left| F_{\theta}\left(y,0,\eps\right) - F^{[k]}(y) \right| \approx \underset{0 \leqslant i \leqslant I}{Max}\left|F_\theta(y_i,\eps) - F^{[k]}(y_i)\right|
    \end{equation}
    and

    \begin{equation}
        \underset{(\tau,y) \in [0,2\pi] \times \Omega}{Max}\left| \phi_{\theta,+}\left(\tau,y,\eps\right) - \phi^{[k]}_\tau(y) \right| \approx \underset{0 \leqslant i \leqslant I,0 \leqslant j \leqslant J}{Max}\left| \phi_{\theta,+}\left( \tau_j , y_i , \eps \right) \right| .
    \end{equation}
    
    \noindent where $\left\{ y_i \right\}_{0 \leqslant i \leqslant I}$ and $\left\{ \tau_j \right\}_{0 \leqslant j \leqslant J}$ are discretisations of $\Omega$ and $[0,2\pi]$ respectively. In our simulations we set $I = 960 = 31^2-1$ (with $\Omega = [-2,2]^2$ represented by the points $\left\{(y_{1,l_1},y_{2,l_2})\right\}_{0 \leqslant l_1 , l_2 \leqslant 30}$) and $J = 30$.\\

    \subsection{Intergals and derivatives representation}

    We also set $\phi^{[0]} = Id$ and $F^{[0]} = \avgf$ for order 0. For order 1 we take, by the following formulas, \eqref{sequence_approximation_variable_change} and \eqref{sequence_approximation_averaged_field}:

    \begin{equation}
        \phi^{[1]}_\tau(y) = y + \eps\int_0^\tau f(\sigma , y)\dd\sigma
    \end{equation}
    
    \noindent and

    \begin{equation}
        F^{[1]}(y) = \left( \frac{\pt \langle \phi_\cdot^{[1]} \rangle }{\pt y}(y) \right)^{-1}\Big\langle f\left(\cdot,\phi_\cdot^{[1]}(y)\right) \Big\rangle.
    \end{equation}

    To compute an integral of the form $\displaystyle\int_a^b g(\sigma)\dd\sigma$, we use a Gauss quadrature with 10 points. \\
    
    Additionally, to approximate the space derivative of a function $g : \RR^2 \mapsto \RR^2$ (Jacobian matrix), we use the following finite difference approximation:

    \begin{equation}
        \frac{\pt g}{\pt y}(y) = \frac{1}{2\eta}\begin{bmatrix}
            g(y+\eta\cdot e_1) - g(y-\eta\cdot e_1) & g(y+\eta\cdot e_2) - g(y-\eta\cdot e_2) 
        \end{bmatrix} + \mathcal{O}\left(\eta^2\right),
    \end{equation}
    
    \noindent where $e_1,e_2$ are the two vectors of the canonical base of $\RR^2$.\\

    Furthermore, if $g,F=\RR^d \longrightarrow \RR^d$, we can approximate the directional derivative $\displaystyle{\frac{\pt g}{\pt y}(y)F(y)}$ using this finite difference approximation:

    \begin{equation}
        \frac{\pt g}{\pt y}(y)F(y) = \frac{1}{2\eta}\big[ g(y+\eta F(y) - g(y-\eta F(y) \big] + \mathcal{O}(\eta^2).
    \end{equation}
    
   We take $\eta = 10^{-5}$ in our case.

   \section{Influence of learning error}

\begin{figure}[H]
		\centering
		\begin{minipage}{0.45\linewidth}
			\includegraphics[width=\linewidth]{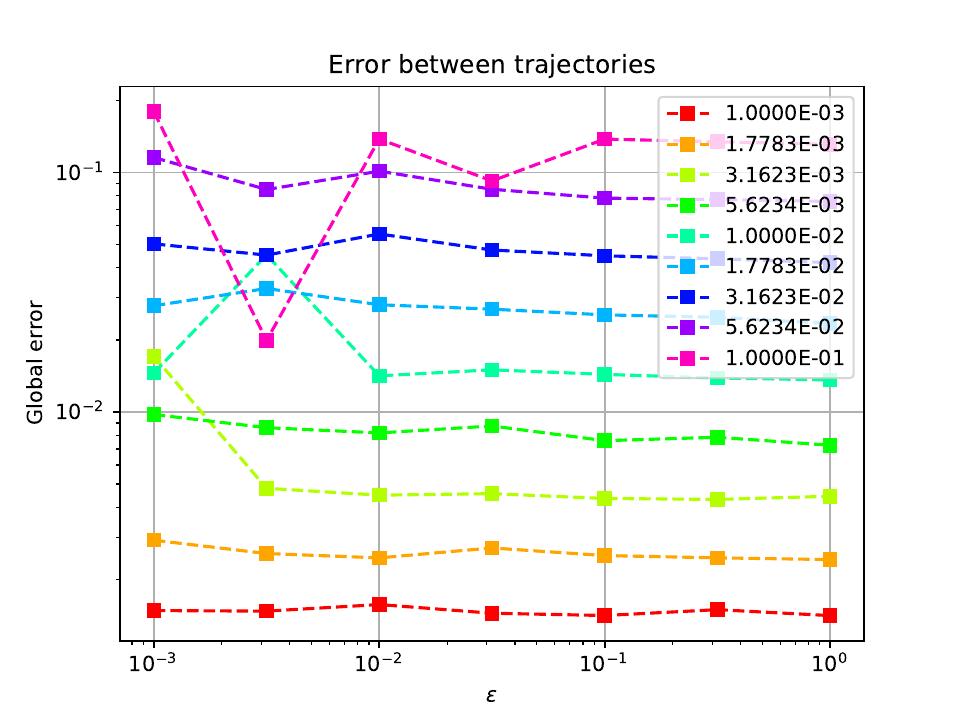}
		\end{minipage}
		\begin{minipage}{0.45\linewidth}
			\includegraphics[width=\linewidth]{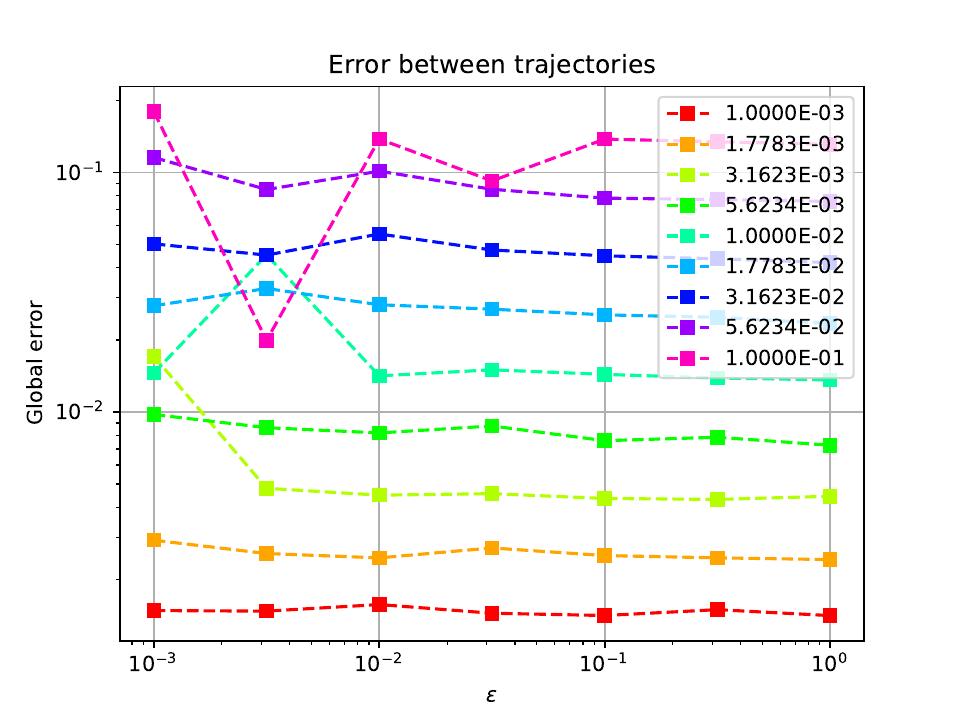}
		\end{minipage}
		\caption{Influence of learning error over integration error (UA test) for Van der Pol oscillator. Left: $K=\numprint{800}$ data, $25$ neurons and $1$ hidden layer per neural network. Right: $K=\numprint{500000}$ data, $150$ neurons and $2$ hidden layers per neural network.}
		\label{fig_Coarse_train_learning_error_influence}
    \end{figure}

   We assess the impact of learning error on integration error, as discussed in Theorems \ref{theorem:highly_oscillatory_general_slow-fast}, \ref{theorem:highly_oscillatory_general_Micro-Macro}, and \ref{theorem:autonomous_highly_oscillatory}. To illustrate this property, we compare Uniform Accuracy (UA) tests after two trainings sessions, using different numbers of hidden layers and neurons with Micro-macro correction method. Figure \ref{fig_Coarse_train_learning_error_influence} demonstrates that efficient learning leads to a reduction in integration error.

    \section{Choice of the parameters}

    \subsection{Inverted Pendulum - Forward Euler method}

    \begin{center}
	{\tiny 
		\begin{tabular}{| l | l |}
			\hline
			\multicolumn{2}{|c|}{\textcolor{red}{\textbf{Parameters}}} \\
			\hline
			\multicolumn{2}{|l|}{\textcolor{DarkGreen}{\textbf{\# Math Parameters:}}} \\
			\hline
			\textbf{Dynamical system:} & Inverted Pendulum \\
            \textbf{Numerical method:} & Forward Euler \\
			\textbf{Interval where step sizes are selected:} &  $[h_-,h_+]=[10^{-3},10^{-1}]$ \\
            \textbf{Interval where small parameters are selected:} &  $[\eps_-,\eps_+]=[10^{-3},1]$ \\
			\textbf{Time for ODE simulation:} &  $T=1$ \\
			\textbf{step size for ODE simulation:} & $h=0.01$ \\
            \textbf{High oscillation parameter for ODE simulation:} & $\eps \in \left\{ 5\cdot10^{-2} , 10^{-3} \right\}$ \\
            \textbf{Initial datum:} & $y^{\eps}(0) = (0.5,-0.5)$ \\
			\hline
			\multicolumn{2}{|l|}{\textcolor{DarkGreen}{\textbf{\# Machine Learning Parameters:}}} \\
			\hline
			\textbf{Domain where initial data are selected:} &  $\Omega = [-2,2]^2$ \\
			\textbf{Number of data:} & $K=\numprint{1000000}$ \\
			\textbf{Proportion of data for training:} & $80 \%$ - $K_0 = \numprint{800000}$ \\
                \textbf{Batch size:} & $B = \numprint{100}$ \\
			\textbf{Hidden layers per MLP:} & $2$ \\
			\textbf{Neurons on each hidden layer:} & $200$ \\
			\textbf{Learning rate:} & $2\cdot 10^{-3}$ \\
			\textbf{Weight decay:} & $1\cdot 10^{-9}$ \\
			\textbf{Epochs:} & $200$ \\
			\hline
		\end{tabular}
		}
    \end{center}

    \noindent\textbf{Computational time for data creation:} 1 h 44 min 19 s\\
    \textbf{Computational time for training:} 9 h 33 min 24s

    \subsection{Inverted Pendululm - midpoint method}

    \begin{center}
	{\tiny 
		\begin{tabular}{| l | l |}
			\hline
			\multicolumn{2}{|c|}{\textcolor{red}{\textbf{Parameters}}} \\
			\hline
			\multicolumn{2}{|l|}{\textcolor{DarkGreen}{\textbf{\# Math Parameters:}}} \\
			\hline
			\textbf{Dynamical system:} & Inverted Pendulum \\
            \textbf{Numerical method:} & midpoint \\
			\textbf{Interval where step sizes are selected:} &  $[h_-,h_+]=[10^{-3},10^{-1}]$ \\
            \textbf{Interval where small parameters are selected:} &  $[\eps_-,\eps_+]=[10^{-3},1]$ \\
			\textbf{Time for ODE simulation:} &  $T=1$ \\
			\textbf{step size for ODE simulation:} & $h=0.01$ \\
            \textbf{High oscillation parameter for ODE simulation:} & $\eps \in \left\{ 5\cdot10^{-2} , 10^{-3} \right\}$ \\
            \textbf{Initial datum:} & $y^{\eps}(0) = (0.5,-0.5)$ \\
			\hline
			\multicolumn{2}{|l|}{\textcolor{DarkGreen}{\textbf{\# Machine Learning Parameters:}}} \\
			\hline
			\textbf{Domain where initial data are selected:} &  $\Omega = [-2,2]^2$ \\
			\textbf{Number of data:} & $K=\numprint{1000000}$ \\
			\textbf{Proportion of data for training:} & $80 \%$ - $K_0 = \numprint{800000}$ \\
                \textbf{Batch size:} & $B = \numprint{100}$ \\
			\textbf{Hidden layers per MLP:} & $2$ \\
			\textbf{Neurons on each hidden layer:} & $200$ \\
			\textbf{Learning rate:} & $2\cdot 10^{-3}$ \\
			\textbf{Weight decay:} & $1\cdot 10^{-9}$ \\
			\textbf{Epochs:} & $200$ \\
			\hline
		\end{tabular}
		}
    \end{center}

    \noindent\textbf{Computational time for data creation:} 1 h 44 min 19 s \footnote{data set used for Inverted Pendulum is the same for both Forward Euler and midpoint.}\\
    \textbf{Computational time for training:} 9 h 11 min 8 s

    \subsection{Van der Pol oscillator - Forward Euler method}

    \begin{center}
	{\tiny 
		\begin{tabular}{| l | l |}
			\hline
			\multicolumn{2}{|c|}{\textcolor{red}{\textbf{Parameters}}} \\
			\hline
			\multicolumn{2}{|l|}{\textcolor{DarkGreen}{\textbf{\# Math Parameters:}}} \\
			\hline
			\textbf{Dynamical system:} & Van der Pol \\
            \textbf{Numerical method:} & Forward Euler \\
			\textbf{Interval where step sizes are selected:} &  $[h_-,h_+]=[10^{-3},10^{-1}]$ \\
            \textbf{Interval where small parameters are selected:} &  $[\eps_-,\eps_+]=[10^{-3},1]$ \\
			\textbf{Time for ODE simulation:} &  $T=1$ \\
			\textbf{step size for ODE simulation:} & $h=0.01$ (for $\eps = 0.1$) and $h=0.001$ (for $\eps = 0.01$)  \\
            \textbf{High oscillation parameter for ODE simulation:} & $\eps \in \left\{ \cdot10^{-1} , 10^{-2} \right\}$ \\
            \textbf{Initial datum:} & $y^{\eps}(0) = (0.5,0.5)$ \\
			\hline
			\multicolumn{2}{|l|}{\textcolor{DarkGreen}{\textbf{\# Machine Learning Parameters:}}} \\
			\hline
			\textbf{Domain where initial data are selected:} &  $\Omega = [-2,2]^2$ \\
			\textbf{Number of data:} & $K=\numprint{20000000}$ \\
			\textbf{Proportion of data for training:} & $80 \%$ - $K_0 = \numprint{16000000}$ \\
                \textbf{Batch size:} & $B = \numprint{200}$ \\
			\textbf{Hidden layers per MLP:} & $2$ \\
			\textbf{Neurons on each hidden layer:} & $200$ \\
			\textbf{Learning rate:} & $2\cdot 10^{-3}$ \\
			\textbf{Weight decay:} & $1\cdot 10^{-9}$ \\
			\textbf{Epochs:} & $200$ \\
			\hline
		\end{tabular}
		}
    \end{center}

    \noindent\textbf{Computational time for data creation:} 2 Days 8 h 2 min 30 s\\
    \textbf{Computational time for training:} 3 Days 13 h 50 min 30 s

    \subsection{Van der Pol oscillator: Comparison between classical and alternative method}

    \begin{center}
	{\tiny 
		\begin{tabular}{| l | l |}
			\hline
			\multicolumn{2}{|c|}{\textcolor{red}{\textbf{Parameters}}} \\
			\hline
			\multicolumn{2}{|l|}{\textcolor{DarkGreen}{\textbf{\# Math Parameters:}}} \\
			\hline
			\textbf{Dynamical system:} & VDP \\
			\textbf{Interval where step sizes are selected:} &  $[h_-,h_+]=[10^{-3},10^{-1}]$ \\
            \textbf{Interval where small parameters are selected:} &  $[\eps_-,\eps_+]=[10^{-3},0.2]$ \\
			\textbf{Time for ODE simulation:} &  $T=1$ \\
			\textbf{step size for ODE simulation:} & $h=0.001$ \\
            \textbf{Small parameter for ODE simulation:} & $\eps\in\left\{10^{-3},10^{-2},10^{-1}\right\}$ \\
            \textbf{Initial datum:} & $y^{\eps}(0) = (0.5,0.5)$ \\
			\hline
			\multicolumn{2}{|l|}{\textcolor{DarkGreen}{\textbf{\# Machine Learning Parameters:}}} \\
			\hline
			\textbf{Domain where initial data are selected:} &  $\Omega = [-2,2]^2$ \\
			\textbf{Number of data:} & $K=\numprint{100000}$ \\
			\textbf{Proportion of data for training:} & $80 \%$ - $K_0 = \numprint{80000}$ \\
                \textbf{Batch size:} & $B = \numprint{100}$ \\
			\textbf{Hidden layers per MLP:} & $2$ \\
			\textbf{Neurons on each hidden layer:} & $200$ \\
			\textbf{Learning rate:} & $2\cdot 10^{-3}$ \\
			\textbf{Weight decay:} & $1\cdot 10^{-9}$ \\
			\textbf{Epochs:} & $200$ \\
			\hline
		\end{tabular}
		}
    \end{center}

    \begin{enumerate}[label = -]
        \item \textbf{Classical method (with auto-encoder):}\\
        \textbf{Computational time for data creation:} 10 min 48 s\\
        \textbf{Computational time for training:} 57 min 03 s\\
        \item \textbf{Alternative method adapted for autonomous case:}\\
        \textbf{Computational time for data creation:} 10 min 52 s\\
        \textbf{Computational time for training:} 48 min 47 s
    \end{enumerate}

\end{document}